\documentclass{article}
\usepackage{mathtools}
\usepackage{mathrsfs,enumerate}
\usepackage{amsmath}
\usepackage{amssymb}
\usepackage[numbers]{natbib}
\usepackage{hyperref}
\usepackage{dsfont}
\usepackage[draft]{todonotes}
\usepackage[a4paper,top=2cm,bottom=2cm,left=2cm,right=2cm,bindingoffset=5mm]{geometry}
\newcommand{\R}{\mathbb{R}}
\newcommand{\T}{\mathbb{T}}

\newcommand{\Z}{\mathbb{Z}}
\newcommand{\cD}{{\ensuremath{\mathcal D}}}

\newcommand{\cP}{{\ensuremath{\mathcal P}}}

\newcommand{\ee}{{\mbox{\boldmath$e$}}}

\newcommand{\oo}{{\mbox{\boldmath$o$}}}

\newcommand{\Kliminf}{K\kern-3pt-\kern-2pt\mathop{\rm
lim\,inf}\limits}  % Kuratowski liminf di insiemi
\newcommand{\Klimsup}{K\kern-3pt-\kern-2pt\mathop{\rm lim\,sup}\limits}  % Kuratowski liminf di insiemi
\renewcommand{\d}{{\mathrm d}}

\newcommand{\restr}[1]{\lower3pt\hbox{$|_{#1}$}}

\newcommand{\nchi}{{\raise.3ex\hbox{$\chi$}}}
\def\qed{\ifmmode % if math mode, assume display: omit penalty etc.
  \else \leavevmode\unskip\penalty9999 \hbox{}\nobreak\hfill
  \fi               
    \qquad           \hbox{\hskip.5em $\square$
                \hskip.1em}}
\def\endproofsym{\qed}

\newenvironment{proof}[1][Proof]{\def\endproofsym{\qed}\trivlist\item[\hskip\labelsep{%
\noindent{\normalfont\emph{#1}.}\hskip .321429\parindent}]\ignorespaces}% Added colon after proof title 1/15/97
{\endproofsym\endtrivlist}

\numberwithin{equation}{section}

\newtheorem{theorem}{Theorem}[section]
\newtheorem{proposition}[theorem]{Proposition}
\newtheorem {lemma}[theorem]{Lemma}

\newtheorem {definition}[theorem]{Definition}

\newtheorem{remark}[theorem]{Remark}

\newenvironment{system}
{\left\lbrace\begin{aligned}}
{\end{aligned}\right.}

\newcommand{\smallplus}{{\raise.3ex\hbox{$\scriptstyle+$}}}
\newcommand{\dive}{\operatorname{div}}

\newcommand{\KL}{K\kern-2pt L}

\def\vep{\varepsilon}
\def\vfi{\varphi}
\def\de{\delta}
\def\lam{\lambda}
\def\be{\begin{equation}}
\def\ee{\end{equation}}
\def\rife#1{\eqref{#1}}
\def\m{\noalign{\medskip}}
\def\mA{{\mathcal A}}
\def\mD{{\mathcal D}}
\def\mO{{\mathcal O}}
\def\into{\int_{\Omega}}
\def\elle#1{L^{#1}(\Omega)}
\def\oo{\overline \Omega}
\def\oQ{\overline Q}

\title{Regularizing effects  of the entropy functional in optimal transport and planning problems}

\author{
Alessio Porretta\thanks{Dipartimento di Matematica, Universit\`a di Roma Tor Vergata. 
Via della Ricerca Scientifica 1, 00133 Roma, Italy.    Email: \texttt{porretta@mat.uniroma2.it}. The author is member of GNAMPA research group of Indam.} 
 }

\date{\today}

\begin{document}

\maketitle
\begin{abstract}
We analyze optimal transport problems with additional entropic cost evaluated along curves in the Wasserstein space which join two probability measures $m_0,m_1$. 
The effect of the additional entropy functional results into an elliptic regularization for the (so-called) Kantorovich potentials of the dual problem. Assuming  the initial and terminal measures to be positive and smooth, we prove that the optimal curve remains smooth for all time.
We focus on the case that the transport problem is set on a convex bounded domain in the $d$-dimensional Euclidean space (with no-flux condition on the boundary), but we also mention the case of Gaussian-like measures in the whole space. The approach follows ideas introduced by P.-L. Lions in the theory  of mean-field games \cite{L-college}. The result provides with a  smooth approximation of minimizers in optimization problems with penalizing congestion terms, which appear in mean-field control or mean-field   planning problems. This allows us to exploit new estimates for this kind of problems by using displacement convexity properties in  the Eulerian approach.
\end{abstract}

%\tableofcontents

\tableofcontents

\section{Introduction}

The aim of this article is to analyze a regularized version of the classical mass optimal transport problem between given measures in the Euclidean space $\R^d$. 
As is well-known, the Kantorovich formulation of this problem reads as follows: given two probability measures $m_0, m_1$ in $\R^d$, find 
\be\label{Ka}
W_2(m_0, m_1):= \min \left\lbrace \int\!\!\!\int_{\R^d \times \R^d}
  |x -y|^2\,\d \gamma(x,y)   : \gamma \in \cP(\R^d \times \R^d), \
  \pi^1_\sharp \gamma= m_0, \ \pi^2_\sharp \gamma= m_1  \right\rbrace,
\ee
where $\gamma$ is a probability measure in $\R^d \times \R^d$, $\pi^i: \R^d \times \R^d \to \R^d$, $i = 1,2$, stand, respectively, for the
$x$ and $y$ projections, and where $\pi^i_\sharp \gamma$ denotes the push forward of the measure $\gamma$ through the mapping $\pi^i$. The dynamic formulation of the problem, due to Benamou and Brenier \cite{BB}, shows that the optimal value is realized by the energy-minimization problem
\be\label{bebr}
W_2(m_0, m_1)= \min \,\,  %\left\{  
\int_0^1\!\! \int_{\R^d}  \frac 12\, |v|^2 dm\,, \qquad (m,v)\,:\quad \begin{cases} m_t - \dive(vm)=0 & \\ m(0)=m_0\,, m(1)=m_1 & \end{cases}
%\right\}
\ee  
where the minimum is meant on all couples $(m,v)$ which satisfy the continuity equation in a suitable sense. The equivalence between \rife{Ka} and \rife{bebr} has a clear geometrical meaning in terms of the Wasserstein space of probability measures
(with finite second moments), since $W_2$ defines  a metric  in  this space (so-called Kantorovich-Wasserstein distance) and the minimum curve $m(t)$ in \rife{bebr}
is the corresponding geodesic connecting $m_0$ and $m_1$. The optimal  velocity field $v$ is actually uniquely associated to this geodesic and the energy term in \rife{bebr} can be interpreted as a classical kinetic energy in terms of the metric derivative of the curve $m(t)$ in the Wasserstein space (\cite{AGS}). It is also known
that the geodesic $m(t)$ coincides with the so-called McCann's displacement interpolation of the optimal plan \cite{McCann}, i.e.  $m(t) =((1-t)\pi^1+t\pi^2)_\sharp \gamma$. For those and many other classical results in  optimal transport theory, we refer to the books  \cite{AGS}, \cite{Sa}, \cite{Vi}. 
\vskip0.5em
In this paper we study a natural regularization of the Wasserstein geodesic and, more generally, of curves $m(t)$ which are optimal for similar transportation costs.
{  It is well-known that the regularity of the geodesic curve, which is minimal  in \rife{bebr},  depends not only on the regularity of the marginals but also on the positivity of the density; however, the positivity set of $m$ may shrink in time,  for $t\in (0,1)$, see e.g.   \cite[Thm 1]{SaWa}. }
%even if the marginals $m_0, m_1$ are smooth, then the geodesic curve which is minimal  in \rife{bebr} may develop singularities. Precisely, this  may happen when the positivity set of $m$ shrinks for $t\in (0,1)$; we refer e.g. to \cite[Thm 1]{SaWa}. 
 This issue motivated 
the suggestion of possible regularizations of the classical optimal transport problem, the most popular being the so-called entropic regularization  which has been intensively investigated in the recent times for numerical efficiency (see \cite{entropic} and references therein). We also refer to \cite{Chizat+al}, \cite{Gentil+al}, \cite{GiTa}, \cite{LMS} for other related entropic perturbations, or regularizations,  of the classical problem.

What we analyze here is the modification of the classical optimal transport functional by addition of  an $\vep$-entropy along the curve $m(t)$. In the model case, this amounts to study
\be\label{eps-func}
\min \,\,  \int_0^T\!\! \!\int_{\R^d}  \frac 12\, |v|^2 dm+ \vep  \int_0^T\!\! \!\int_{\R^d}  \log\left(\frac{dm}{d\nu}\right)dm \,,\qquad (m,v)\,:\quad \begin{cases} m_t - \dive(vm)=0 & \\ m(0)=m_0\,, m(T)=m_1 & \end{cases} 
%\right\}
\ee
where the additional entropy is computed with respect to  a reference measure $\nu$. Typical choices for $\nu$ include the classical Lebesgue measure, in the case that the transport problem is restricted to a bounded domain $\Omega$, or the case of Gaussian measures, for problems in the whole space.
\vskip1em

There are   several motivations  for the interest in this problem. 
First of all, we will see that the additional entropy term produces a regularization of the geodesics of the classical optimal transport problem. Indeed, the
addition of  $\vep$-entropy in the functional yields an elliptic regularization which is reminiscent of  the vanishing viscosity approximation in first order problems. This elliptic regularization is observed in terms of the Kantorovich potential provided by the dual problem; exploiting this approach may enrich the viewpoints
on optimal transport problems.

Secondly, this kind of regularization is natural from the intrinsic geometry of   optimal transport, since it can be readily interpreted as a deformation of the Wasserstein metric which regularizes the corresponding  geodesics.  In particular, even if we develop here this analysis in the Euclidean setting, the extension to Riemannian manifolds will be natural. We plan to exploit this issue in future work.

As a further motivation, this approach provides with a natural setting where Eulerian calculus can be fully justified. In particular, we will prove  
estimates related to displacement convexity inequalities which are naturally robust in  this kind of approximation.

%\begin{itemize}
%\item The additional entropy term produces a regularization of the geodesics of the classical optimal transport problem. In particular, for smooth and positive marginals $m_0,m_1$, the minimizer $m(t)$ remains smooth for all times.
%
%\item The addition of  $\vep$-entropy in the functional yields an elliptic regularization which is reminiscent of  the vanishing viscosity approximation in first order problems. This elliptic regularization is observed in terms of the Kantorovich potential given by the dual problem, exploiting this approach may enrich the viewpoints
%on optimal transport problems.
%
%
%\item This kind of regularization is natural from the intrinsic geometrical viewpoint of the optimal transport, since it only amounts to add an entropy cost criterion in the selection of optimal curves in the Wasserstein space. It can therefore be readily interpreted as a deformation of the Wasserstein metric which regularizes the corresponding  geodesics.  In particular, even if we develop here this analysis in the Euclidean setting, the extension to Riemannian manifolds will be natural (see also Remark \ref{bochner}).
%
%\item This kind of approximation provides with a natural setting where Eulerian calculus can be fully justified. In particular, many estimates related to displacement convexity are naturally robust in this kind of approximation and can be fully developed.  
%
%\item Last but not least, this kind of approximation is natural in the framework of  mean-field theories of optimal control and differential games
%
%\end{itemize}  
\vskip0.5em
To put things into perspective, this kind of study was originally motivated by mean-field theories of optimal control and differential games. To this purpose, 
it is necessary to embed the classical optimal transport problem into a larger class of dynamical optimization problems, where 
the cost criterion involves both a  kinetic-type energy of the curve $m(t)$ and additional terms which minimize congestion effects:
\be\label{mfg-cost}
\min\,\, \int_0^T\!\!\! \int_{\R^d}  \frac12|v|^2\, dm+    \int_0^T\!\!\! \int_{\R^d} F(m(t))dxdt  \,:\qquad \begin{cases} m_t - \dive(vm)=0 & \\ m(0)=m_0\,, m(T)=m_1 & \end{cases}  
\ee
In \rife{mfg-cost}, we will suppose $F$ to be a nondecreasing function defined on the density of the curve $m(t)$, so that the functional will reduce concentration effects on the minimal configurations.

Minimization problems such as \rife{mfg-cost} have been extensively studied in the last decade, especially in connection with mean-field control theory. 
In particular, the introduction of mean field game theory by J.-M. Lasry and P.-L. Lions (\cite{LL1}, \cite{LL-japan}, \cite{L-college}) boosted the interest in this kind of problems, giving an interpretation of the minima in terms of Nash equilibria of large populations' differential games. In that context, any generic agent is represented as a dynamical state $x(t)$ and aims at minimizing an individual cost given by
\be\label{mfg-Lag}
\min \quad \int_0^1  \frac{|x'(t)|^2}2  dt +\int_0^1\!\! \into f(m(t)) dxdt
\ee
where $m(t)$ is a probability density which represents the distribution law of the  population. More precisely,  the agents consider $m(t)$ as an exogenous data  (an anticipated guess of the collective state) which at the Nash equilibrium coincides with the density of the population driven by the optimal strategies of the individuals. It turns out  that, when $f= F'$,  the couple $(u,m)$ given, respectively,  by the value function  and by  the density  of the agents, satisfy the system
\be\label{mfg-mod}
\begin{system}
&-u_t + \frac{|Du|^2}2  = F'(m)  &&\text{in }   (0,T)\times \Omega, \\
&m_t  - \dive (m\, Du)  = 0  &&\text{in } (0,T)\times \Omega,
\end{system}
\ee
which coincides with  the state-adjoint state system of the optimization problem \rife{mfg-cost}. 
%Here $H(p)$ is the Hamiltonian associated to $L$ and $H_p$ its derivative: $H_p(Du)$ denotes the optimal velocity field in the dynamic formulation of optimal transport.
%
%The optimality conditions satisfied by minima of problem \rife{mfg-cost} yield a couple $(u,m)$ solution of the system
 %In this context, 
 
The analysis of problem \rife{mfg-mod}, together with natural generalizations, was initiated   by P.-L. Lions in his courses at Coll\`ege de France devoted to  mean field games, and then pursued by many authors in different directions,  mostly relying on the study of the primal-dual optimization problems and/or on the Lagrangian formulation \rife{mfg-Lag}. For a non exhaustive list of contributions, see \cite{Carda1}, \cite{CG},  \cite{CGPT},  \cite{CMS}, \cite{GMST}, \cite{Gomes2}, \cite{LaSa}, \cite{OPS}. Let us stress that in this  literature the state space if often assumed to be  the flat torus (periodic case) and some final cost is often prescribed rather than prescribing the final marginal $m(T)$.
\vskip0.4em

In this article, we develop one of the pioneering ideas suggested by P.-L. Lions in this context (\cite[Lessons 27/11 and 4/12 2009]{L-college}), which relies on the reduction of   system \rife{mfg-mod} to a single elliptic equation on the function $u$. This equation is possibly degenerate in the set where $m$ vanishes, but adding an entropic term along the curve $m(t)$ preserves the ellipticity in a strong form. In the particular case that $F(m)= \vep m(\log m -1)$, which corresponds to the functional \rife{eps-func} with uniform Lebesgue measure $\nu$, this equation takes the following form:
\be\label{model2}
-u_{tt}   + 2 Du\cdot Du_t - D^2uDu \cdot Du - \vep \Delta u =0\,.
\ee
Let us recall that here $Du$ is the optimal velocity field for the functional \rife{eps-func}, and the corresponding $m$ can be reconstructed as
 $m:= \exp\left(\frac1\vep\left( \frac{|Du|^2}2-u_t\right)\right)$.

In this viewpoint, it is clear that the perturbed functional \rife{eps-func} yields a {\it vanishing viscosity} approximation of the first order case represented by optimal transport trajectories. It would actually be justified to call problem \rife{eps-func} as the {\it vanishing entropic} approximation of optimal transport. 
It is to be noticed that the optimal transportation problem corresponds to  equation \rife{model2} complemented with nonlinear Neumann boundary conditions at time $t=0$ and $t=T$:
\be\label{model-bou}
\begin{cases} 
u_t= \frac{|Du|^2}2 - \vep \log (m_0) & \hbox{at $t=0$,}\\
 u_t= \frac{|Du|^2}2 - \vep \log (m_1) & \hbox{at $t=T$.}
 \end{cases}
\ee
The purpose of this article is  to analyze the $\vep$-entropic problem \rife{eps-func} from this new perspective. Our  results include two main goals: 

\noindent (i) to give a suitable existence and regularity result for the nonlinear elliptic problem \rife{model2}--\rife{model-bou} in the case that  the marginals $m_0,m_1$ are positive and smooth. This problem is also identified as a singular limit of penalized optimal transportation problems.

%\noindent (ii) giving estimates for $t$ inside $(0,T)$ which may be uniform with respect to $\vep$, to be applied for the case that $m_0, m_1$ are not strictly positive.

\noindent (ii) to use this smooth setting to justify different type of a priori estimates, possibly independent of the positivity of $m$ (and possibly uniform as $\vep\to0$). Some of those estimates directly come from displacement convexity arguments developed in what is called the Eulerian viewpoint. This will  provide an extension of regularity results which were obtained for mean field game systems  \rife{mfg-mod} with different approaches (\cite{LaSa}).
\vskip0.5em

In order not to overlap with existing results (mostly restricted to dynamics on the flat torus), and in order to provide tools which may be useful in further applications, we develop our analysis in general convex subdomains $\Omega\subset \R^d$ where  a no-flux condition is prescribed on the continuity equation, which implies conservation of mass. In addition, it is convenient to add a potential energy term to the functional; this introduces some inhomogeneity in the problem and allows us to include the case where the entropy additional term is computed with respect to general reference measures $\nu$. 
%In  the specific case of the optimal transportation problem \rife{eps-func}, we also extend the result to the  case $\Omega= \R^d$ when $\nu$ is a Gaussian-type measure. This serves as a prototype of extension of our results to non compact manifolds, which seems of interest in the optimal transport literature.
\vskip0.5em

Let us now be more precise on the statement of our results, starting from the very model case of quadratic cost. 

\begin{theorem}\label{OT1} Let $\Omega$ be a convex bounded smooth domain, and assume that $V\in W^{2,\infty}(\Omega)$.
Let $m_0, m_1\in  {\cP} (\Omega)\cap W^{1,\infty}(\Omega)$ and suppose that $m_0, m_1>0 $ in $\oo$.  Then there exists $u\in C^2(Q)\cap C^1(\oQ)$ which is a classical solution of the problem
\be\label{elli-mod}
\begin{cases}
-u_{tt}   + 2 Du\cdot Du_t - D^2uDu \cdot Du - \vep \Delta u + DV\cdot Du=0 & \hbox{in $(0,T)\times \Omega$,}
\\
Du\cdot \vec\nu=0 & \hbox{on $(0,T)\times \partial \Omega$,}
\\
u_t= \frac{|Du|^2}2 - \vep \log (m_0) - V(x) & \hbox{at $t=0$, $x\in \Omega$}\\
 u_t= \frac{|Du|^2}2 - \vep \log (m_1) - V(x) & \hbox{at $t=T$, $x\in \Omega$.}
\end{cases}
\ee
Moreover, the function $u$ is unique up to a constant and  the function $m$ defined as
\be\label{mexp}
m:= \exp\left(\frac1\vep \left(  \frac{|Du|^2}2 -u_t-V(x)\right)\right)
\ee
is the unique minimizer of the functional 
\be\label{eps-func-neu}
B(m,v):= \int_0^T\!\! \!\int_{\R^d}  \frac 12\, |v|^2 dm+ \vep  \int_0^T\!\! \!\int_{\R^d}  \log\left(\frac{dm}{d\varrho}\right)dm \,,\quad (m,v)\,:\, \begin{cases} m_t - \dive(vm)=0 &\hbox{in $(0,T)\times \Omega$,} \\ 
v\cdot \vec \nu =0 & \hbox{on $(0,T)\times \partial \Omega$,}\\
m(0)=m_0\,, m(T)=m_1 & \end{cases} 
\ee
where $\varrho= e^{- \frac{V(x)}\vep}dx$. In particular, we have  that $m\in C^1(Q)\cap C^0(\oQ)$ and  $m>0$ in $\oQ$.
%Then problem \rife{OT}  admits  a solution $(u,m)$ such that $u\in C^2(Q), m\in C^1(Q)$, and this solution is unique (up to addition of a constant to $u$). Moreover we have  $m(t) >0 $ in $\oo$ for every $t\in [0,T]$. Finally, $m$ is the unique minimizer of the functional
%$$
%\inf \{ \int_0^T\into   m\, |v|^2dxdt + \int_0^T\into  \log\left( \frac {dm}{d\nu}\right) dm\}\,,\qquad m_t = \dive(mv)=0\,,\,\,\, m(0)=m_0, \, m(1)=m_1
%$$ 
%where $\nu= e^{- \frac{V(x)}\vep}dx$.
\end{theorem}

The  proof of Theorem \ref{OT1} relies on gradient estimates for quasilinear elliptic equations. To this respect, we follow Lions' approach but we also exploit 
some extension. In particular, even if we restrict here our analysis to positive measures, we give extra estimates which are local in time  and  independent of the positivity of the marginals $m_0$, $m_1$. Precisely, we prove that   for $t\in (0,T)$ there exists  a constant $K_t$ (depending only on the upper bounds of $u,m$, and on the distance of $t$ from $t=0,T$) such that
$$
\theta\, |Du(t)|^2 + \vep \log m(t) \leq K_t \qquad \forall t\in (0,T)\,,
$$ 
for some  $\theta\in (0,1)$. This gives  a quantitative estimate of the gradient of $u$   on the support of $m$, which we hope can be useful for further generalizations to non strictly positive measures.   

An extension of Theorem \ref{OT1} is also provided for the case of optimal transport in the whole space, relying on  the  dissipation property of the potential $V$, whenever the additional entropy is taken with respect to Gaussian-like measures. This will serve  as a prototype for a  development of similar analysis on non compact manifolds, which seems of interest in the optimal transport literature. For a more readable statement, here we rescale the potential $V$ according to the parameter $\vep$.

\begin{theorem}\label{OT2}   Assume that $V\in W^{2,\infty}(\R^d)$ satisfies $D^2 V(x) \xi\cdot \xi \geq \gamma\, |\xi|^2 $ for some $\gamma>0$, and for every $x,\xi\in \R^d$.
Let $m_0, m_1\in  C^1(\R^d)$ satisfy
$$
m_0 e^V, m_1 e^V \in W^{1,\infty}(\R^d)\qquad \hbox{and}\quad c_0 \, e^{-V(x)}\leq m_0\leq C_0 \, e^{-V(x)}\,,\qquad c_1 \, e^{-V(x)}\leq m_1\leq C_1 \, e^{-V(x)}
$$
for some positive constants $c_i, C_i>0$, $i=0,1$.
Then the problem \rife{eps-func} with $\nu= e^{- V(x) }dx$ admits a  unique minimizer $m$  given by $m= e^{-V} \, \exp\left(\frac1\vep \left(  \frac{|Du|^2}2 -u_t\right)\right)$, where $u$ is a solution of 
$$
-u_{tt}   + 2 Du\cdot Du_t - D^2uDu \cdot Du - \vep (\Delta u - DV\cdot Du) =0 \qquad \hbox{in $(0,T)\times \R^d$.}
$$
Moreover, we have  $m\in C^1([0,T]\times \R^d)$ and there exist constants $\gamma_0, \gamma_1$ (depending on $\vep, \gamma, c_0,c_1, C_0, C_1$) such that
$$
\gamma_0 e^{-V(x) } \leq m(t,x) \leq \gamma_1 e^{-V(x) } \qquad \forall (t,x)\in [0,T]\times \R^d\,.
$$  
\end{theorem}

\vskip0.5em
As we mentioned before,  the results on the classical optimal transport with quadratic cost are embedded in a larger class of results provided for the system 
\begin{equation}\label{mfpp0}
\begin{system}
&-u_t + H( Du)  = f(m) + \vep \log m + V(x)  &&\text{in }  Q, \\
&m_t  - \dive (m\,H_p(Du))  = 0  &&\text{in } Q,
\\
& H_p(Du)\cdot \vec\nu=0 &&\text{on $(0,T)\times \partial \Omega$,}
\\
&m(0,\cdot) = m_0, \; m(T,\cdot)  = m_1 &&\text{in } \Omega\,,
\end{system}
\end{equation}
where $\Omega$ is a smooth bounded convex subset of $\R^d$,  and $Q=(0,T)\times \Omega$.  Here $H$ is a strictly convex radial function, 
 with power-like growth ($H\simeq |Du|^q$ with $q>1$) and $f$ is a nondecreasing function, satisfying quite general conditions, which are precisely  given in the next Section. 

When $H$ is uniformly convex and with quadratic growth,  the  results that we  obtain for problem \rife{mfpp0} include  both the existence of smooth solutions (when the marginals $m_0, m_1$ are positive) and further global and local estimates independent of the positivity of $m$. In the following statement we gather  results which appear, in more details, in Theorem \ref{exi-smooth} and Proposition \ref{supm}. 

\begin{theorem}\label{main}
Let $\Omega$ be a  bounded smooth convex domain  in $\R^d$ and let $V\in W^{2,\infty}(\Omega)$. Let $H= h(|p|) $ for some (strictly) convex increasing function $h\in C^3([0,\infty))$ with $h'(0)=0$. Assume in addition that $H$ satisfies conditions \rife{Hpp}-\rife{hppp}, and 
%$$
%\exists\,\, \alpha, \beta>0\,:\, \quad \alpha I_d \leq H_{pp} \leq \beta I_d \qquad \forall p\in \R^d\,.
%$$ 
let $f\in C^2(0,\infty)$ be a nondecreasing function.  
%satisfying condition \rife{fprime}. 
Assume that $m_0,m_1\in W^{1,\infty}(\Omega)$ and $m_0,m_1>0$ in $\oo$. 

Then  there exists  a unique couple $(u,m)$ such that  $u\in C^{2,\alpha}(Q)\cap C^{1,\alpha}(\oQ)$, $m\in C^{1,\alpha}(Q)\cap C^{0,\alpha}(\oQ)$, $\into u(T)m_1=0$, $\into m(t)=1$ for all $t$, and  \rife{mfpp0} is satisfied in a classical sense.  We also have that $m>0$ in $\oQ$ and is a minimizer of the corresponding optimal transport problem, while  $u$ is  a solution to the elliptic problem
\be\label{quasi0}
\begin{system}
& -{\rm tr}\left({\mathcal A}({\mathcal D}u) \, {\mathcal D}^2u \right) + DV(x) \cdot H_p(Du)=0 &&\text{in } Q, \\%(0,1)\times \R^d, \\
& \,\, \, H_p(Du)\cdot \vec\nu=0   &&\text{in } (0,T)\times \partial \Omega, \\%(0,1)\times \R^d
& -u_t + H(Du)- V(x) =   f(m_1)+ \vep \log (m_1)   &&\hbox{at $t=T$, $x\in \Omega$,}
\\
&   -u_t + H(Du)- V(x) = f (m_0) + \vep \log(m_0) &&\hbox{at $t=0$, $x\in \Omega$,}
\end{system}
\ee
where $\mA$ is given by \rife{aij}. 
\vskip0.3em
Moreover, if $f$  satisfies
$$
\exists \,\, c_0, r_0 >0\,:\,\quad f'(r) \geq \frac{c_0} r  \qquad\forall r\geq r_0\,,
$$
then there exist constants $K_0, K_1$, independent of $\vep$, such that, for $1\leq p\leq \infty$: 
\be\label{globm1}
\| m(t)\|_{\elle p} \leq K_0\, (\|m_0\|_{\elle p}+ \|m_1\|_{\elle p}+1) \qquad \forall t\in [0,T]
\ee 
and  
\be\label{locm1}
\| m(t)\|_{\elle p} \leq K_1\, \left(t^{-q}+ (T-t)^{-q} \right) \qquad \forall t\in (0,T)
\ee
for some $q>0$. 
\end{theorem}

%We point out that Theorem \ref{main} gathers results which appear, in more details, in Theorem \ref{exi-smooth} and Proposition \ref{supm}. 
Let us comment on the above statement in relation to the previous literature which analyzed problem \rife{mfpp0} (often called mean-field game systems). 

In the first  part of Theorem \ref{main}, we prove that the system has smooth solutions whenever an $\vep$-entropy term is added; this generalizes Theorem \ref{OT1} (in which $H(p)=\frac{|p|^2}2$ and $f=0$). Previous results showing the existence of smooth solutions were obtained for the periodic case ($\Omega$ is the torus)  in \cite{Mu}, when a final condition is prescribed on $u$ rather than for the transport problem; this latter one was previously discussed in \cite{L-college}, as we already said. Recent results for the one-dimensional case can also be found in \cite{Gomes2}, \cite{munoz1}.  

{  Let us stress that the main difference of the transport problem, compared to the case when  a final condition is prescribed on $u$,   lies in the estimate of the sup-norms of $m$ and $u$. Indeed, while prescribing a final cost  fixes the $L^\infty$-norm of $u$, this is no longer the case for the transport problem, which is invariant by addition of  constants to $u$. This leads to a different strategy in order to close the gradient estimate, because this latter one depends on the $L^\infty$-bound of $u$, see Theorem \ref{bound}. In the present paper we circumvent this problem through the obtention of preliminary estimates on $m$ and through the choice of  a suitable normalization for $u$. Both those ingredients contain novel estimates,   locally and globally in time. In particular, this is where estimates \rife{globm1} and \rife{locm1} play a role. }Such global and local (in time) $L^p$ estimates   only depend on the behavior of $f(m)$ for large $m$, so they are not relying on the positivity of the marginals. Letting $\vep\to 0$, those estimates are exported for so-called weak solutions of mean field game systems, as introduced in \cite{Carda1}, \cite{CG} and studied in \cite{CGPT}, \cite{CMS}, \cite{OPS}. A statement in this direction can  be found in Theorem \ref{weaksol}.  Notice that those weak solutions coincide with relaxed minima of the functionals. In particular, the estimates \rife{globm1}, \rife{locm1} extend similar estimates proved by H. Lavenant and F. Santambrogio \cite{LaSa} for the case $H(p)= \frac{|p|^2}2$ using  flow-interchange techniques. Our approach is different and relies on displacement convexity inequalities which were first exploited by D. Gomes and T. Seneci for mean-field game systems on the torus \cite{GoSe}. To this respect, our approach can be seen as the Eulerian route  to the $L^p$ estimates of $m$, compared to the Lagrangian and gradient flow approach developed in \cite{LaSa}. {  Let us observe that the global estimates \rife{globm1} are classically known for the plain optimal transport problem \rife{bebr}, since Mc Cann's work \cite{McCann}, as a consequence of the geodesic convexity of $m\mapsto \into m^p$.  The local estimates \rife{locm1} are instead a new effect induced by entropic or congestion-like terms  in the functional.}

Finally,  we will extend some of the previous results to  more general Hamiltonians.  In Theorem \ref{exi-smooth} we include the case that the Hamiltonian has $q$-growth in the gradient, although the uniform convexity is still needed to get smooth solutions. Conversely, the case that $H= |Du|^q$ corresponds to possible degeneracy (or singularity) of $H_{pp}$. This   case is considered in Theorem \ref{exi-div},  where we obtain solutions in a weaker sense, namely the elliptic equation should be understood in divergence form. {  In both of these latter cases, where we quit the quadratic growth of the Hamiltonian,  we can no longer rely on  the displacement convexity estimates; our strategy is changed accordingly, and  we use the convexity of the potential $V(x)$ to recover preliminary $L^\infty$ bounds on $m$.}
%, nevertheless the optimal transport problem still possesses a solution $m\in C^{0,\alpha}(\oQ)$.

\vskip0.5em
Let us conclude by presenting the material of the next sections.
In Section \ref{not} we introduce precise assumptions and notations. In Section \ref{elli} we discuss the nonlinear elliptic problem \rife{quasi0} and we derive the main gradient estimates which are needed to construct smooth solutions. Here we essentially follow the same strategy introduced by P.L.-Lions and refined in \cite{Mu}, which leads to gradient bounds depending on the $L^\infty$-norm of $u$ (Theorem \ref{bound}). 
%In addition, we also show that $L^\infty$ bounds for $m$ can be deduced 

In Section \ref{pena}, we build solutions for a penalized approximation of the  optimal transportation problem,  in which the constraint on $m$ at time $t=0,t=T$ is relaxed. This is similar to the natural relaxation occurring in control theory, when exact controllability problems are relaxed into approximate controllability problems penalizing the final target. This construction is important in our strategy, since it allows us to show how a suitable normalization of $u$ can be controlled by $L^\infty$-bounds of $m$ (see Lemma \ref{norm}). In Section \ref{disconv} we extend to our setting the  displacement convexity inequalities obtained in \cite{GoSe} (see also  \cite{Gomes2}) and we develop new $L^p$ estimates for the density, both globally and locally  in time (Proposition \ref{supm}). This is where we recover and generalize the results of \cite{LaSa}.  In Section \ref{exi}, we deduce our main results, say Theorem \ref{exi-smooth} and \ref{exi-div}, with specific attention to the model case of  classical optimal transportation. We also prove a few extra results which show possible further developments;  for example, we consider the limit as the logarithmic term goes to zero, showing convergence towards weak solutions of mean-field game systems (Theorem \ref{weaksol}). 

%show in Theorem \ref{OTspace} the role played by the potential $V(x)$ . 
Finally, we have detailed in  the Appendix the construction of solutions to the elliptic equation  with nonlinear Neumann boundary conditions. Unfortunately, due to the specific origin of our elliptic problem, which is set in the cylindrical domain $(0,T)\times \Omega$, we cannot rely on classical results (mostly obtained in smooth domains or in nonsmooth domains with simpler operators). This is why we  provide a self-contained proof, where we use a reflection argument to handle the Neumann condition on $\partial \Omega$. It is only in this step that we use the radial structure of the Hamiltonian. This extra condition would not be used in the periodic case, or if $x$ was taken in  a compact manifold without boundary.

\section{Notations and assumptions}\label{not}

Throughout the paper, we denote by $\Omega$ an open bounded convex subset of $\R^d$, $d\geq 1$. We denote  by $\partial \Omega$ the boundary of $\Omega$ and by $\vec \nu$ the unit outward normal vector on $\partial \Omega$. We will assume that $\Omega$ is  of class $C^3$   
(i.e. $\partial \Omega$ is locally the graph of a $C^3$ function) and, for $x\in \Omega$,  we denote by $d(x):= $dist$(x, \partial\Omega)$ the distance function to the boundary. In particular, by smoothness of  $\Omega$, we have that $d(x)$ is a $C^2$ function in the set $\{x\in \Omega: 0<d(x)<\de_0\}$ for some $\de_0>0$.  In addition, we also have $D^2 d\leq 0$ since $\Omega$ is convex. 

Given $\Omega$, we denote by $\cP(\Omega)$ the space of probability measures in $\Omega$, and by $L^p(\Omega)$ the standard Lebesgue space, for $p\in [1,\infty]$. The space $W^{k,\infty}(\Omega)$ is the space of functions with bounded $k-th$ derivatives (in weak sense); in particular, $W^{1,\infty}(\Omega)$ coincides with the space of Lipschitz continuous functions, and $W^{2,\infty}(\Omega)$ with the space of $C^1$ functions whose first derivatives are Lipschitz continuous.

Given $T>0$ and $t\in (0,T)$, we consider the following first order evolution system in $(0,T)\times \Omega$:
\begin{equation}\label{mfpp}
\begin{system}
&-  u_t  + H( Du)  = f(m) + \vep \log m + V(x)  &&\text{in }  Q:= (0,T)\times \Omega, \\
& m_t   - \dive (m\,H_p(Du))  = 0  &&\text{in } Q,
\\
& H_p(Du)\cdot \vec\nu=0 &&\text{on $(0,T)\times \partial \Omega$,}
\\
&m(0,\cdot) = m_0, \; m(T,\cdot)  = m_1 &&\text{in } \Omega\,,
\end{system}
\end{equation}
where $u,m$ are functions of $(t,x)$ and $u_t, m_t$ denote the partial   derivatives in time, $Du$ the gradient vector of $x$-derivatives, and $\dive(\cdot) $ is the divergence operator in the $x$-variables ($\dive (F)= \sum_{i=1}^d \partial_{x_i} F_i$).

In \rife{mfpp}, we assume that  $V(x)$ is a Lipschitz continuous function (but it will be often required to be in $W^{2,\infty}(\Omega)$) and that
\be\label{fm}
\hbox{$f: (0,\infty)\to \R $ is a  nondecreasing $C^1$ function.}
\ee
Of course, we could have embedded the term $\vep \log m$ into $f$, but we decided to make it explicit the contribution coming from the entropy functional.

The function $H:\R^{d} \to \R$ (so-called Hamiltonian function) will be assumed smooth (at least of class $C^2$)
%\footnote{quanto regolare servano $H,f$ vorrei un po' ricontrollarlo. A posteriori dovrebbe essere $H$ di classe $C^3$  e $f$ di classe $C^2$ (per le stime necessarie), ma il risultato di Lieberman che uso richiederebbe i coefficienti della matrice (e  quindi $H_{pp}, f'$) di classe $C^{2,\alpha}$. Al momento cmq mi sembra un problema minore...} 
and  convex. 
In order to build smooth solutions (say, $u$ of class $C^2$) we will strengthen both the regularity and the convexity of $H$. In the simplest setting, this will require $H$ to be a $C^3$ function satisfying the following conditions (hereafter $H_p, H_{pp}, H_{ppp}$ denote  the derivatives of $H(p)$):
%\be\label{H}
%H:\R^{N} \to \R\quad \hbox{smooth and uniformly convex}
%\ee
%%\be\label{f}
%%f: (0,\infty) \to \R \quad \hbox{smooth and increasing}
%%\ee
%we set $y:=(t,x)$, for $t\in (0,T), x\in \Omega$, and ${\mathcal D}:= (\partial_t , D_x )$. We define
%\be\label{matrix}
%{\mathcal A}:= \begin{pmatrix} 1 & -H_p(Du) \\ -H_p(Du) & H_p(Du)\otimes H_p(Du)+ (mf'(m)+\vep) 
%H_{pp}(Du)\end{pmatrix}
%\ee
%where $m= \vfi(-u_t + H(Du)-V(x))$ and where we denote with $D.$ the $x$-derivative.   
%\be\label{coerc}
%H_p(p)\cdot p - H(p) \geq \gamma_0\, |p|^2 -\gamma_1 \qquad \forall p\in \R^d
%\ee
%for some $\gamma_0,\gamma_1>0$ and 
\be\label{Hpp}
\exists \, \alpha_H, \beta_H >0\, : \qquad \alpha_H \, I_d \leq H_{pp} (p) \leq \beta_H \, I_d \qquad \forall p\in \R^d
\ee
where $I_d$ is the identity matrix in $\R^d$, and 
%. We will also need in some results to require that
\be\label{hppp}
\hbox{$H_{ppp}(p)$ is uniformly bounded for  $p\in \R^d$.}
\ee
{  Conditions \rife{Hpp}-\rife{hppp} describe functions $H$ which are uniformly convex and with quadratic growth, whose model case is obviously   $H= \frac{|p|^2}2$.
This set of assumptions will be required in the displacement convexity estimates, Section \ref{disconv}, and in the statements which rely on  them.}

However, in most of our results the two above conditions can be generalized into the following form, modeled on the case of Hamiltonian with superlinear growth  of power $q$:
%\be\label{coercq}
%H_p(p)\cdot p - H(p) \geq \gamma_0\, |p|^q -\gamma_1 \qquad \forall p\in \R^d
%\ee
%for some $\gamma_0,\gamma_1>0$ and 
\be\label{Hppq}
\exists \, \alpha_H, \beta_H >0\,\hbox{and $\varpi\geq 0$} : \qquad \alpha_H \, (|p|+ \varpi)^{q-2}\, I_d \leq H_{pp} (p) \leq \beta_H\, (|p|+ \varpi)^{q-2}  \, I_d \qquad \forall p\in \R^d\setminus\{0\}
\ee
\be\label{hpppq}
|H_{ppp}(p)| \leq \gamma  (1+ |p|)^{\frac 32(q-2)}\qquad \forall p\in \R^d\,. 
\ee
%where $q>1$.
Let us notice that if $H\in C^2$ is strictly convex and  satisfies $\alpha \, |p|^{q-2} \leq H_{pp}(p) \leq \beta \, |p|^{q-2} $ only for $|p|$ large, then it satisfies \rife{Hppq} (with  $\varpi=1$) for every $p\in \R^d$.  

{  Assumptions \rife{Hppq}-\rife{hpppq}  describe functions $H$ with $q$-growth in the gradient  (thus generalizing \rife{Hpp}-\rife{hppp} which correspond to $q=2$), whose simplest model is given by $H= (|p|^2+\varpi)^{\frac{q}2}$.}
Of course, the case when $\varpi=0$  is special, and corresponds to singular, or degenerate, $H_{pp}$. In that case, we will be able to show only existence of weak, rather than classical, solutions.
\vskip0.4em
We now rephrase system \rife{mfpp} into a single elliptic equation in $(t,x)$ variables. To this purpose, we denote
$$
f^\vep(r):= f(r)+ \vep \log r\,,\qquad\quad \vfi(r):= (f^\vep )^{-1}(r)= (f+\vep \log )^{-1}(r)\,.
$$
Hence,  the first equation in \rife{mfpp} implies
\be\label{m}
m= \vfi(-u_t + H( Du) - V(x))\,.
\ee
Computing formally, we have,  using  the first equation, 
\begin{align*}
& (f'(m)+ \vep /m) m_t = [-  u_{tt} + H_p( Du)\cdot Du_t]
\\
\m
& (f'(m)+ \vep /m) \dive (m\,H_p(Du)) = [-Du_t\cdot H_p( Du)+ D^2 u H_p(Du)\cdot H_p(Du)- DV(x)\cdot H_p(Du)] 
\\
& \qquad \qquad \qquad \qquad +(f'(m)m+ \vep ){\rm tr}\left(H_{pp}(Du)D^2u \right)\,.
\end{align*}
Therefore, subtracting the previous two terms and using the second equation in \rife{mfpp}, we get 
%\rife{mfpp} is formally equivalent to
$$
\begin{cases}
 m= (f+\vep \log)^{-1}(-u_t + H( Du)-V(x) ) &
 %&\text{in } Q, \\%(0,1)\times \R^d, \\
 \\
 -u_{tt} + 2H_p( Du)\cdot Du_t - D^2 u H_p(Du)\cdot H_p(Du) - ( mf'(m)+ \vep){\rm tr}\left(H_{pp}(Du)D^2u \right)+ DV(x)\cdot H_p(Du)= 0  &
 %&\text{in } Q,
%\\%(0,1)\times \R^d,\\
%&m(0,\cdot) = m_0, \; m(1,\cdot)  = m_1 &&\text{in } \R^d.
\end{cases}
$$
%where we used that $\vfi'(-u_t + H( Du))= \frac1{f'(m)}$ by definition of $\vfi$. 
Eventually,  the second equation of the system can be shortly written as  a quasilinear equation in the space-time variables:
$$
-{\rm tr}\left({\mathcal A} \, {\mathcal D}^2u \right)+ DV(x)\cdot H_p(Du)= 0
$$
where
$$
{\mathcal D}:= (\partial_t , D_x )\,,\qquad {\mathcal A}:= \begin{pmatrix} 1 & -H_p(Du) \\ -H_p(Du) & H_p(Du)\otimes H_p(Du)+ (mf'(m)+\vep) 
H_{pp}(Du)\end{pmatrix}
$$
and we recall that $m$ is a function of ${\mathcal D}u$ (see \rife{m}).
Let us observe that ${\mathcal A} $ is  the contribution of two terms, namely
\be\label{aij}
{\mathcal A}:= {\mathcal A_1}+ {\mathcal A_2}= \begin{pmatrix} 1 & -H_p(Du) \\ -H_p(Du) & H_p(Du)\otimes H_p(Du)\end{pmatrix}
+ \begin{pmatrix} 0 & 0 \\ 0 & (mf'(m)+\vep) 
H_{pp}(Du)\end{pmatrix}
\ee
and it can be   checked that  ${\mathcal A}$ is   elliptic provided $H_{pp}>0$:
\be\label{ella}
\hbox{if $\eta=(s,\xi)$} \qquad\Rightarrow \quad {\mathcal A}\eta\cdot \eta = (s-H_p(Du)\cdot \xi)^2+ (mf'(m)+\vep) H_{pp}(Du)\xi\cdot \xi > 0 \quad \hbox{if $\eta\neq 0$}\,.
\ee
Let us notice that, when $\vep=0$, the operator becomes degenerate elliptic; in particular, the matrix ${\mathcal A}$ would  be elliptic in the set $\{m>0\}$ and would  possibly degenerate only in the set where $m$ vanishes, assuming $f$ increasing. Otherwise, even if $\vep>0$, $ {\mathcal A}$ can degenerate if $H_{pp}$ vanishes at some point.

%Indeed, assuming $H$ to be strictly convex, and since $f$ is  increasing, the matrix ${\mathcal A}$ turns out to be elliptic in the set $\{m>0\}$ and may possibly degenerate only in the set where $m$ vanishes. This actually occurs {\it unless $f$ is singular at $m=0$} and blows-up as a logarithm. 
%In the model example $f(m)=\log m$, $H=\frac12|p|^2$ (corresponding to the entropic regularization of the Wasserstein metric), we have $mf'(m)=1$ and ${\mathcal A}$ is uniformly elliptic.

We   point out that the space-time elliptic problem also enjoys a  divergence form structure, which is clearly induced by the continuity equation in \rife{mfpp}. This is very important for considering possibly weak formulations. Anyway, in the following we will first approach the elliptic problem in its non-divergence form, through the use of maximum principle.

Finally, we observe that the planning conditions which prescribe the marginals of $m$ turn into a nonlinear Neumann condition for $u$:
$$
 -u_t(0,\cdot) + H( Du(0,\cdot))= f^\vep(m_0)+ V \,\,;\qquad  -u_t(T,\cdot) + H( Du(T,\cdot))=  f^\vep(m_1)+V\,.
$$
To conclude, problem \rife{mfpp} can be rephrased as the following quasilinear elliptic problem:
\be\label{quasell}
\begin{system}
& -{\rm tr}\left({\mathcal A}(x,{\mathcal D}u) \, {\mathcal D}^2u \right)+   DV(x)\cdot H_p(Du) =0 &&\text{in } Q, \\%(0,1)\times \R^d, \\
& -u_t + H(Du)=   f^\vep (m_1) + V(x) &&\hbox{at $t=T$, $x\in \Omega$,}
\\
&   -u_t + H(Du) = f^\vep(m_0)+ V(x)   &&\hbox{at $t=0$, $x\in \Omega$,}
\\
& \,\, H_p(Du)\cdot \vec\nu=0 &&\text{on $(0,T)\times \partial \Omega$}
%\\%(0,1)\times \R^d,\\
%&m(0,\cdot) = m_0, \; m(1,\cdot)  = m_1 &&\text{in } \R^d.
\end{system}
\ee
where $f^\vep(r)= f(r)+ \vep \log r$. Unfortunately, in \rife{quasell}  the boundary condition takes the form of a nonlinear Neumann condition defined piecewisely on the 
time-space boundary. This makes it harder to  construct smooth solutions; it is known that, in general quasilinear problems with nonlinear boundary conditions,  Lipschitz solutions may even lack of $C^1$-regularity   if the domain is not sufficiently smooth. This is the only reason why we will simplify our setting requiring the Hamiltonian $H$ to be a radial function, namely that
\be\label{hrad}
H(p)= h(|p|)\qquad \hbox{for some $h\in C^1([0,\infty))$: $h'(0)=0$.}
\ee
This structure assumption on $H$  reduces the boundary condition on $\partial \Omega$ to a standard Neumann condition. The only point where this is needed appears in the Appendix, where we prove the $C^{1,\alpha}$ regularity (up to the boundary) for the   elliptic equation. Of course, all conditions given above on $H$ can be rephrased in terms of the real function $h(r)$, but we decided to keep the whole exposition for a general function $H$, since all results would stand in the general case up to  a modification of the $C^{1,\alpha}$ result in the Appendix. Moreover, all results would hold for general $H$, not necessarily radial, in case of compact state space, e.g. if the state variable $x$ belongs to the flat $d$-dimensional torus.

\section{Quasilinear elliptic equations}\label{elli}

In this Section we directly study the quasilinear elliptic equation
\be\label{quasi-ell}
-{\rm tr}\left({\mathcal A}(y,{\mathcal D}u) \, {\mathcal D}^2u \right)= f\,,
\ee
where $u= u(y)$, $y\in \mO\subset \R^n$.

We start with  a computational lemma which underlines the typical structure  used for gradient bounds of (possibly degenerate) quasilinear equations.

\begin{lemma}\label{struc} Let $u\in C^3(\mO)$ be a  solution to \rife{quasi-ell}, with $f\in C^1(\mO)$ and $\mathcal A(y,\eta) \in C^1(\mO\times \R^n)$. Then we have:
\vskip1em
\noindent (i)  for any $M\in \R$,  the function $(u+M)^2$ satisfies
\be\label{u2}
-{\rm tr}\left({\mathcal A}(y,{\mathcal D}u) \, {\mathcal D}^2 (u+M)^2 \right) = 2 f\, (u+M) - 2 \mA(y,\mD u) \mD u\cdot \mD u
\ee
(ii) for any $C^2$ function $\psi: \R^n \to \R$, we have that the function 
$$
w:= \psi(\mD u)
$$
satisfies
\be\label{conv}
-{\rm tr}\left({\mathcal A} \, {\mathcal D}^2 w \right) = \mD f\cdot \mD\psi +  {\rm tr}\left([\mD_\eta {\mathcal A} \cdot \mD w+\mD_y {\mathcal A} \cdot \mD \psi ] \,\mD^2 u \right) 
-   \sum\limits_{k,l=1}^n (\mD^2 \psi)_{k\ell } \left( \mA  \,\mD u_k \cdot \mD u_\ell\right)
\ee 
where $\mA$ and $\psi$ are computed on $\mD u$, and where $\mD_\eta {\mathcal A} \cdot \mD w$ and $\mD_y {\mathcal A} \cdot \mD \psi$ are  the matrices with $ij$-th component given, respectively,  by $\sum_{\ell=1}^n \partial_{\eta_\ell} (a_{ij} ) w_\ell$ and $\sum_{k=1}^n    \partial_{y_k} (a_{ij} )\psi_k$.
\end{lemma}

\proof To obtain \rife{u2}, it is enough to observe that
$$
{\mathcal D}^2 (u+M)^2= 2 (u+M) {\mathcal D}^2 u + 2 {\mathcal D}u \otimes \mD u\,.
$$
Applying the matrix $\mA(y,\mD u)$ and taking traces,    \rife{u2} follows from \rife{quasi-ell}.

In the case of (ii),  we compute  (denoting $z_k= \partial_k z$ the partial derivatives):
$$
w_i= \psi_k (\mD u)  u_{ki}\,,\qquad w_{ij}= \psi_k (\mD u)  u_{kij} + \psi_{k\ell} u_{ki}u_{\ell j}
$$
so that 
\begin{align*}
{\rm tr}\left({\mathcal A}(y,{\mathcal D}u) \, {\mathcal D}^2 w \right) & = \psi_k (\mD u) \partial_k \left[a_{ij}u_{ij}\right] - \psi_k (\mD u) \partial_{\eta_\ell} (a_{ij} ) u_{\ell k}u_{ij}
-  \partial_{y_k} (a_{ij})  \psi_k (\mD u) u_{ij}
+ \psi_{k\ell} a_{ij}   u_{ki}u_{\ell j}
\\
& 
= -\psi_k (\mD u) \partial_k f-  \left(\partial_{\eta_\ell} (a_{ij} ) w_\ell \right)   u_{ij} -  \partial_{y_k} (a_{ij}) \psi_k (\mD u)   u_{ij}   + (\mD^2 \psi)_{k\ell } \left( \mA  \mD u_k \cdot \mD u_\ell\right)
\end{align*}
\qed
 
\vskip1em
Now we consider our specific case, where the matrix $\mA$ is induced by the optimal transport problem and given by \rife{aij}. 
In particular, we  consider the solution $u$ to the elliptic problem
\be\label{eps-de}
\begin{system}
& -{\rm tr}\left({\mathcal A}(x,{\mathcal D}u) \, {\mathcal D}^2u \right)+ \rho u + DV(x)\cdot H_p(Du) =0 &&\text{in } Q, \\%(0,1)\times \R^d, \\
& -u_t + H(Du)= \de u+ f^\vep (m_1) + V(x) &&\hbox{at $t=T$, $x\in \Omega$,}
\\
&   -u_t + H(Du)+ \de u= f^\vep(m_0)+ V(x)   &&\hbox{at $t=0$, $x\in \Omega$,}
\\
& \,\, H_p(Du)\cdot \vec\nu=0 &&\text{on $(0,T)\times \partial \Omega$,}
%\\%(0,1)\times \R^d,\\
%&m(0,\cdot) = m_0, \; m(1,\cdot)  = m_1 &&\text{in } \R^d.
\end{system}
\ee
where $f^\vep(s)= f(s)+ \vep \log s$.  Hereafter we denote
$$
\Sigma_0:= \{0\} \times \Omega\,,\,\, \Sigma_T:= \{T\} \times \Omega\,.
$$
Throughout all this Section,  we assume that $H\in C^3(\R^d)$ and satisfies  conditions \rife{Hppq}-\rife{hpppq} for some $q>1$. Let us explicitly note that assumption \rife{Hppq}, together with the regularity of $H$, also implies the following: 
\be\label{coercq}
H_p(p)\cdot p - H(p) \geq \gamma_0\, |p|^q -\gamma_1 \qquad \forall p\in \R^d
\ee
for some $\gamma_0,\gamma_1>0$
%\footnote{as can be proved by writing 
%$$
%H_p(p)\cdot p - H(p)= \int_0^1 \frac{d}{dt} [H_p(tp)\cdot (tp) - H(tp)] dt =  \int_0^1  [H_{pp} (tp)p\cdot p] \,t\,dt  
%% \geq \alpha \int_0^1  t|p|^2(|p| t+ \varpi)^{q-2} dt  
%$$}
.
\vskip1em

We stress that the additional terms $\rho u$ and $\de u$ in \rife{eps-de}, added in the interior and, respectively, on the boundary, reinforce the coercivity of the elliptic operator and guarantee the uniqueness of solutions and the validity of the maximum principle.  We start with a bound which is uniform on $\rho$.

\begin{lemma}\label{sup} Let $u\in C^2(Q)\cap C^1(\overline Q)$ be a solution of \rife{eps-de}. Then we have 
$$
\de \|u\|_\infty \leq   \left( \|f^\vep (m_0)+V\|_\infty + \|f^\vep (m_1)+V\|_\infty  \right).
$$
%for some $C$ depending on $\|V\|_{W^{1,\infty}} , T, H, \Omega$.
\end{lemma}

\proof  Apply the maximum principle (respectively, the minimum principle) to $u-C$ (respectively, $u+C$), with $C=\frac1\de(\|f^\vep(m_0)+V\|_\infty + \|f^\vep(m_1)+V\|_\infty)$.
\qed

\vskip1em
The next step is the main tool in this approach and consists in a gradient bound for solutions to \rife{eps-de}. The idea  was introduced  by P.-L. Lions in \cite[Lessons 27/11 and 4/12 2009]{L-college},  following  the lines of the classical Bernstein  method.   It was recently refined in \cite{Mu} to handle the general case with $x$-dependent nonlinearities. In the aforementioned results,  $\Omega$ was assumed to be the flat torus. We extend here this  kind of estimates for the case of bounded domains with Neumann conditions.
\vskip1em
%We start with an upper bound on $m$.

We start with a consequence of Lemma \ref{struc}.

\begin{lemma}\label{ut} Let $u\in C^2(Q)\cap C^1(\overline Q)$ be a solution of \rife{eps-de}. Then we have
\be\label{lemut}
u_t(t,x) \leq \sup_{\Sigma_0\cup \Sigma_T} (u_t)_+ \qquad \hbox{and}\qquad |u_t(t,x)| \leq \sup_{\Sigma_0\cup \Sigma_T} |u_t| \qquad \forall (t,x)\in \oQ\,. 
\ee 
%both $u_t$ and $-u_t$ achieve their maximum on $\Sigma_0\cup \Sigma_T$.
\end{lemma}

\proof Since $\mathcal A$ is independent of $t$, from \rife{conv} we deduce that $w=u_t$ (and $w=-u_t$) solve
$$
-{\rm tr}\left({\mathcal A} \, {\mathcal D}^2 w \right) + H_{pp}(Du)  DV(x) \cdot Dw -  {\rm tr}\left( [\mD_\eta {\mathcal A} \cdot \mD w]\cD^2 u \right)+ \rho w=0
$$
so $w$ does not have positive maximum inside. On $\partial \Omega$ we have (taking time derivative of the Neumann condition)
\be\label{neutime}
Dw \cdot H_{pp}(Du)\vec \nu=0\,.
\ee
If  $w$ has a maximum on $\partial \Omega$,  the tangential derivative would vanish so that  
$$
Dw \cdot (H_{pp}(Du)\vec \nu) = Dw \cdot \vec \nu (H_{pp}(Du)\vec \nu\cdot \vec \nu) \,. 
$$
which implies $Dw \cdot \vec \nu=0$ by \rife{neutime}.  With a typical perturbation argument (replacing $w$ with $w+ \theta d(x)$ and letting $\theta\to0$),  we conclude that  the maximum of $w$ cannot be attained on $\partial\Omega$. Hence
$$
w \leq   \sup_{\Sigma_0\cup \Sigma_T} w_+ 
$$
for $w= \pm u_t$, which yields \rife{lemut}.
\qed

%\begin{lemma}\label{mound} Let $u\in C^2(Q)\cap C^1(\overline Q)$ be a solution of \rife{eps-de}. Then, there exists a constant $C$, independent of $\rho$ and $\de$, such that
%$$
% -u_t  + H(Du) \leq C  \,.
%$$
%\end{lemma}
%
%\proof  We use \rife{conv} with the convex function $\psi((s,\xi))= -s + H(\xi)$. This means that $w:= -u_t + H(Du)$ satisfies
%$$
%-{\rm tr}\left({\mathcal A} \, {\mathcal D}^2 w \right) + \rho (-u_t +H_p(Du)Du)\leq    -  {\rm tr}\left([\mD {\mathcal A} \cdot \mD w] \,\mD^2 u \right) \,. 
%$$
%Since $-u_t +H_p(Du)Du\geq w$, we deduce that
%$$
%-{\rm tr}\left({\mathcal A} \, {\mathcal D}^2 w \right) + \vep w\leq    -  {\rm tr}\left([\mD {\mathcal A} \cdot \mD w] \,\mD^2 u \right)\,.
%$$
%This implies that $w$ cannot have positive maximum inside the domain. By Lemma \ref{sup}, at $t=0,t=1$ we have
%$$
%w \leq C  \left( \|f(m_0)\|_\infty + \|f(m_1)\|_\infty \right)\,.
%$$
%By maximum principle, we conclude that $w$ is bounded above by a constant only depending on $\|f(m_0)\|_\infty$, $\|f(m_1)\|_\infty$. 
%\qed
%
%\vskip1em
%Notice that the above lemma guarantees that $m:= f^{-1}(-\partial_t u + H(Du))$ is bounded above.  
\vskip1em
Now we estimate the derivatives of $u$ in terms of the sup-norm of $u$. 
%Due to the nonlinear Neumann condition on $\partial \Omega$, we will assume that $\Omega$ is  a  strictly convex domain, in the sense that there exists $\gamma>0$ such that, in a neighborhood of the boundary, we have 
%\be\label{striconv}
%D^2 d(x) \leq -\gamma I_d \qquad \forall x\in \Omega\, :\, d(x) \leq \de_0
%\ee
%where $d(x)=$ dist$(x, \partial \Omega)$. This condition is not needed if $H$ is radial, and we believe it unnecessary even for the general case, but this may  require further technicalities which are beyond our present scopes.

\begin{theorem}\label{bound} Let $u\in C^3(Q)\cap C^2(\overline Q)$ be a solution  of \rife{eps-de}. Assume that $H$ is a $C^3$ radial function which satisfies conditions \rife{Hppq}--\rife{hpppq},   that $f\in C^2(0,\infty)$ and there exist    $ \beta, s_0>0$ such that 
\be\label{fprime}
| f''(s)|\, s^2 \leq \beta \,   \sqrt{(1+ f'(s)s)^3} \qquad \hbox{ for all $s\geq s_0$. }
\ee
%Assume In addition to \rife{coerc}--\rife{Hpp}, assume that $H\in C^3$ and $H_{ppp}$ is uniformly bounded on $\R^d$.  
%
Then, there exists a constant $C$, independent of $\rho$ and $\de$, such that
$$
\| \mD u\|_\infty  \leq C (1+ \|u\|_\infty)\,.
$$
Here the constant  $C$ depends on $\vep, \beta,   \|f^\vep(m_0)\|_{W^{1,\infty}}, \|f^\vep(m_1)\|_{W^{1,\infty}} , \|V\|_{W^{2,\infty}}$ and  on the functions $H,f$.
%\noindent Moreover, if $H$ is radial, the estimate is independent of $\gamma$ and the assumption \rife{striconv} is not needed.
\end{theorem}

\begin{remark}
\rm We stress that, in the proof below, the radial form of $H$ is only used to simplify the treatment of the boundary condition $H_p(Du)\cdot\vec\nu=0$. However, the same result of Theorem \ref{bound} holds true for general $H$ (not radial)  if we assume that $\Omega$ is  a  strictly convex domain, in the sense that there exists $\gamma>0$ such that, in a $\de_0$-neighborhood of the boundary, we have 
\be\label{striconv}
D^2 d(x) \leq -\gamma I_d \qquad \forall x\in \Omega\, :\, d(x) \leq \de_0
\ee
where $d(x)=$ dist$(x, \partial \Omega)$. This reinforced convexity would be enough to handle the nonlinear condition at $\partial \Omega$.
\end{remark}

\begin{remark} \rm {  Assumption \rife{fprime} is a technical condition, which however is satisfied by most natural examples, such as powers and exponentials. Functions 
which do not satisfy \rife{fprime} can be highly oscillating such as, for instance, if $f'(s)= s(1+\sin(s^\gamma))$, with $\gamma>1$.

Let us stress, however,  that the growth of $f$  becomes irrelevant once we obtain an $L^\infty$ bound on $m$ {\it prior to the gradient estimate} of $u$. This is what we will exploit later,  with different arguments, in Lemma \ref{mound} (using the convexity of the potential $V$) or in Proposition \ref{supm} (for quadratic growth Hamiltonians). 

At this stage, we decided to keep Theorem \ref{bound} independent of any a priori estimate of $m$, at the expense of requiring the additional condition  \rife{fprime}.}
\end{remark}

\proof 

Step 1 {\it (tuning the $u$)}:   we replace $u$ with 
$$
v:= u+ M- C_0\frac{(T-t)}T\,,\qquad \hbox{with $M= 2\|u\|_\infty+1$, $\,\,\,C_0= 2M= 2(1+ 2\|u\|_\infty)$.}
$$
On account of \rife{eps-de}, the function $v$ is therefore a solution to the following problem
\be\label{pbv}
\begin{system}
& -{\rm tr}\left({\mathcal A}({\mathcal D}u) \, {\mathcal D}^2v \right)+ DV(x)\cdot H_p(Dv)+ \rho v =\rho(M-C_0\frac{(T-t)}T) &&\text{in } Q, \\
& -v_t + H(Dv)- V(x)= \de u+ f^\vep (m_1)  -C_0/T&&\hbox{at $t=T$, $x\in \Omega$,}
\\
&   -v_t + H(Dv)- V(x) + \de u= f^\vep(m_0)  -C_0/T &&\hbox{at $t=0$, $x\in \Omega$}
\\
& \,\,\, H_p(Dv)\cdot \vec\nu=0 && \hbox{on $(0,T)\times \partial \Omega$}
%\\%(0,1)\times \R^d,\\
%&m(0,\cdot) = m_0, \; m(1,\cdot)  = m_1 &&\text{in } \R^d.
\end{system}
\ee
where $f^\vep(s)= f(s)+ \vep \log s$. 
We observe that 
$$
\|v\|_\infty \leq C (1+\|u\|_\infty)
$$
where, here and below, we denote by $C$ any (possibly different) constant independent of $\rho$ and $\de$. We also notice that $v$ satisfies
\be\label{vbound}
v(T) \geq 1+ \|u\|_\infty \geq 1\,,\qquad v(0) \leq - (1+ \|u\|_\infty) \leq -1 \,.
\ee
Now we define the function
$$
z:= H(Dv) + \frac\lam 2 v^2
$$
where
\be\label{lamb}
\lambda = \frac \sigma{1+ \|u\|_\infty^2}   
\ee
for  some (small) constant $\sigma$ to be chosen later. The goal now is   to estimate the function $z$ through the maximum principle.
\vskip1em
Step 2 ({\it the boundary behavior}). Suppose that $z$ attains a maximum at the boundary $t=T$,  in which case we have $Dz=0$ and $z_t \geq 0$ on the maximum point. By using the boundary condition for $v$, we compute
\begin{align*}
z_t& = H_p(Dv) \cdot  Dv_t+ \lambda v\, v_t 
\\
& = H_p(Dv) \cdot DH(Dv) - H_p(Dv)\cdot DV(x) - H_p(Dv)\cdot D f^\vep(m_1) - \de H_p(Dv)\cdot Dv  \\
& \qquad + \lam v\, ( C_0/T +   H(Dv) - \de u - f^\vep(m_1)- V)
\\
& = H_p(Dv) \cdot Dz - \de H_p(Dv)\cdot Dv  - H_p(Dv)\cdot DV(x) - H_p(Dv)\cdot D f^\vep(m_1) \\
& \qquad - \lam v \left( H_p(Dv)Dv-H(Dv)     + f^\vep(m_1)+ V - C_0/T +\de u\right) 
\end{align*}
hence
$$
z_t- H_p(Dv) \cdot Dz + \de  H_p(Dv)\cdot Dv  \leq   \left(\|DV\|_\infty+ \|Df^\vep (m_1)\|_\infty\right) |H_p(Dv)| - \lam v \left( H_p(Dv)Dv-H(Dv)  - K\right)
$$
where $K   = K ( \|f^\vep (m_1)\|_\infty ,  \| V\|_\infty ,    \|\de u\|_\infty , T)$.
\vskip0.4em
Due to \rife{coercq},  we can suppose that $H_p(Dv)Dv-H(Dv)  >2K$, otherwise we have $ \gamma_0 |Dv|^q \leq \gamma_1 + 2K$, hence 
$\max z \leq C\, K + \lam \|u\|_\infty^2$ and we are done.  
Therefore,   using $v(T)\geq 1+ \|u\|_\infty$, \rife{lamb} and \rife{coercq}, and since $|H_p(Dv)| \lesssim |Dv|^{q-1}$ by \rife{Hppq}, we estimate:
\begin{align*}
z_t- H_p(Dv)\cdot Dz  & \leq   \left(\|DV\|_\infty+ \|Df^\vep (m_1)\|_\infty\right) |H_p(Dv)| - C\, \lam \, (1+\|u\|_\infty) [H_p(Dv)Dv-H(Dv)] 
\\
& \leq C \left(\|DV\|_\infty+ \|Df^\vep (m_1)\|_\infty\right) |Dv|^{q-1} - \sigma \frac C{(1+\|u\|_\infty)} |Dv|^q
\end{align*}
and we conclude that
$$
z_t- H_p(Dv)\cdot Dz  < 0  \quad \hbox{if $|Dv| > \frac C\sigma \left(\|DV\|_\infty+ \|Df^\vep(m_1)\|_\infty\right)(1+\|u\|_\infty)$.}
$$
This implies that $z$ cannot have a maximum at $t=T$ unless $|Dv| \leq \frac {K_0}\sigma(1+\|u\|_\infty)$
for some constant $K_0>0$ independent of $\rho,\de$. Notice that the above argument only needs $H_p(Dv)\cdot Dz=0$ on the maximum points, which holds even on $\partial\Omega$ because $H_p(Dv)$ is tangential. 

Similarly we reason for $t=0$ using that $v(0) \leq - 1- \|u\|_\infty \leq -1$.  Thus, we conclude that $z$ cannot have maximum at $t=0, t=T$ unless $|Dv|$ is uniformly bounded.

%We have
%\begin{align*}
%z_t & = H_p(Dv)\cdot D^2v Dv + \de |Dv|^2 - Dv\cdot D f(m_0)  + \lam v\, ( H(Dv) + \de v+ C_0 - \de (M-C_0) - f(m_0))
%\\
%& = H_p(Dv)\cdot Dz + 2\de z -Dv\cdot DV(x) - Dv\cdot D f(m_0) - \lam v \left( H_p(Dv)Dv-H(Dv)     + f(m_0) - C_0+\de (M-C_0)\right) 
%\end{align*}
%which implies, for $K \simeq \|f(m_0)\|_\infty + C \|u\|_\infty$, that  
%$$
%z_t- H_p(Dv)\cdot Dz - 2\de z \geq  -\left( \|DV\|_\infty+ \|Df(m_0)\|_\infty\right) |Dv| - \lam v \left( H_p(Dv)Dv-H(Dv) - K\right)
%$$
%since $v(0)\leq 0$. In fact we have $ v(0) \leq - 1- \|u\|_\infty \leq -1$ and we conclude exactly as before that
%$$
%z_t- H_p(Dv)\cdot Dz - 2\de z \geq 0\quad \hbox{provided $ |Dv|  > \frac C\sigma \left( \|DV\|_\infty+ \|Df(m_0)\|_\infty\right)(1+\|u\|_\infty)$.}
%$$
%Again, this means that $z$ cannot have a maximum at $t=0$ unless $|Dv| \leq C$
%for $C= C(\sigma, \|DV\|_\infty, \|f(m_0)\|_{W^{1,\infty}}, \|u\|_\infty )$.
\vskip1em
Finally, we look at $(0,T)\times \partial \Omega$, where we use the Neumann condition for $v$. 
%$$
%Dz \cdot \vec \nu=   D^2v Dv\cdot  \vec \nu 
%= D ( Dv\cdot \vec\nu)\cdot Dv- D\vec \nu DvDv\,.
%$$
%Since $Dv\cdot \vec\nu=0$ on $\partial \Omega$, its gradient vector is orthogonal to the boundary, so  $D ( Dv\cdot \vec\nu)\cdot Dv=0$ because $v$ has zero normal derivative. Moreover, if $d(x)$ denotes the distance from the boundary, we have  $\vec \nu= -Dd(x)$, hence $D\vec \nu= D^2 d$ is negative since $\Omega$ is convex. We conclude that 
%$$
%Dz \cdot \vec \nu\leq 0\qquad \hbox{on $\partial \Omega \times (0,T)$}
%$$
%\vskip1em
Here we have
\begin{align*}
Dz\cdot Dd(x) & =  H_p(Dv)\cdot [D^2v  Dd(x)] + \lam v\, Dv\cdot  Dd(x) \\
& = H_p(Dv)\cdot D (Dv\cdot Dd(x))- D^2 d(x) H_p(Dv)\cdot Dv
\end{align*}
where we used that $H_p(Dv)\cdot \vec\nu =0$ implies $Dv\cdot Dd(x)=0$ on $\partial\Omega$ (because $H$ is radial). This also implies    that   $D(Dv\cdot Dd(x))$    is parallel to the normal direction,  so it is orthogonal to $H_p(Dv)$ using again the Neumann condition. Hence  $H_p(Dv)\cdot D (Dv\cdot Dd(x))=0$. Finally, using that $H_p(Dv)$ is parallel to $Dv$ and that $\Omega$ is convex,  we conclude that   
$$
Dz\cdot Dd(x)= - D^2 d(x) H_p(Dv)\cdot Dv\geq 0\,.
$$
This implies (up to the usual perturbation argument which consists in replacing $u$ with $u+ \theta d(x)$ and letting $\theta\to 0$) that no maximum point of $z$ can occur on $\partial\Omega$.
\vskip1em

Step 3 ({\it the interior estimate}).  We use \rife{u2} and we get
$$
-{\rm tr}\left({\mathcal A}(x,{\mathcal D}u) \, {\mathcal D}^2 \left[\frac{v^2}2\right] \right)+ DV(x) \cdot H_p(Dv)\, v + \rho  \,v^2  =   \rho \left(M-C_0\frac{(T-t)}T\right)\, v -   \mA(\mD u) \mD v\cdot \mD v\,.
$$
Then we use 
\rife{conv} with $\psi( v_t, Dv)= H(Dv)$ and we get
\be\label{hdueq}
\begin{split}
-{\rm tr}\left({\mathcal A} \, {\mathcal D}^2 \left[H(Dv) \right]\right) & + H_{pp}(Dv) DV \cdot D^2vH_p(Dv) + D^2V(x) H_p(Dv)H_p( Dv)  + \rho H_p(Dv)\cdot Dv \\
& =     {\rm tr}\left(\left\{ [\mD_\eta {\mathcal A} \cdot \mD [H(Dv)] ] +[ D_x {\mathcal A} \cdot H_p(Dv)]\right\} \,\mD^2 v \right) 
-   \sum\limits_{k,\ell=1}^d  H_{p_kp_\ell}(Dv) \left( \mA  \,\mD v_k \cdot \mD v_\ell\right)
\end{split}
 \ee 
where $\mA$ is computed on $(x,{\mathcal D}u)$.  Summing up the previous equalities and noticing that $D^2vH_p(Dv)= Dz-\lambda \, v Dv$, we deduce
 \begin{align*}
& -{\rm tr}\left({\mathcal A} \, {\mathcal D}^2 z \right)+ H_{pp}DV \cdot Dz+ \rho z + \lambda v\,  \left(DV\cdot H_p-H_{pp}DV\cdot Dv\right)   +   D^2V(x) H_p \cdot H_p 
\\
& + \rho\left( H_p\cdot Dv- H + \lam \frac{v^2}2\right)  =  \lam \,\rho \left(M-C_0\frac{(T-t)}T\right)\, v -   \lam\,\mA  \mD v\cdot \mD v
\\  &+  {\rm tr}\left(\left\{ [\mD_\eta {\mathcal A} \cdot \mD z]+[ D_x {\mathcal A} \cdot H_p ]\right\} \,\mD^2 v \right)  
-  \lam \sum\limits_{i,j,\ell=1}^n v \partial_{\eta_\ell} (a_{ij}) \partial_\ell v \,v_{ij}
 -  \sum\limits_{k,\ell=1}^d  H_{p_kp_\ell}\left( \mA  \,\mD v_k \cdot \mD v_\ell\right)\,
 \end{align*}
 where $H_p, H_{pp}$ are computed on $Dv$. 
 The latter term in the left-hand side of the equality can be dropped by positivity.  Hence we deduce
 \be\label{list1}
 \begin{split}
& -{\rm tr}\left({\mathcal A} \, {\mathcal D}^2 z \right)+ H_{pp}DV \cdot Dz+ \rho z 
 + \lambda v\,  \left(DV\cdot H_p-H_{pp}DV\cdot Dv\right)   +   D^2V(x) H_p \cdot H_p  
 \\
&  + \lam\,\mA  \mD v\cdot \mD v + \sum\limits_{k,\ell=1}^d  H_{p_kp_\ell}\left( \mA  \,\mD v_k \cdot \mD v_\ell\right) \leq   \lam \,|v|\, \rho\, M    
\\  &+  {\rm tr}\left(\left\{ [\mD_\eta {\mathcal A} \cdot \mD z]+[ D_x {\mathcal A} \cdot H_p ]\right\} \,\mD^2 v \right)  
-  \lam \sum\limits_{i,j,\ell=1}^n v \partial_{\eta_\ell} (a_{ij}) \partial_\ell v \,v_{ij}\,.
 \end{split}
 \ee
We use the precise definition of ${\mathcal A}$ and  the coercivity of $H_{pp}$ from \rife{Hppq}. Thus there exists $p_0,\gamma_H>0$ such that,  if $|Dv|>p_0$, we have
\begin{align*}
\mA  \mD v\cdot \mD v  & \geq  |-v_t+ H_p(Dv)  \cdot Dv |^2+   \gamma_H(\vep + \chi(m)) |Dv|^q \\
\sum\limits_{k,\ell=1}^d  H_{p_kp_\ell}\left( \mA  \,\mD v_k \cdot \mD v_\ell\right) &  \geq \gamma_H\,|Dv|^{q-2} \left( |Dv_t- D^2v H_p(Dv)|^2  + \gamma_H  \,|Dv|^{q-2} (\vep + \chi(m))  \, |D^2 v|^2\right)
\end{align*} 
where we denoted $\chi(m)= f'(m)m$. Using these inequalities in \rife{list1} we get
\be\label{eqz}\begin{split} &
-{\rm tr}\left({\mathcal A} \, {\mathcal D}^2 z \right)+ H_{pp}DV \cdot Dz+ \rho z 
%\\
%& \qquad 
+ \lambda v\,  \left(DV\cdot H_p-H_{pp}DV\cdot Dv\right) +   D^2V H_p \cdot H_p 
%+ \rho\left( H_p\cdot Dv- H\right) 
\\
& 
 +\lam \,\, \big( |-v_t+ H_p(Dv)  \cdot Dv |^2+   \gamma_H(\vep + \chi(m)) |Dv|^q \big)
\\ & \qquad  
 + \gamma_H\,|Dv|^{q-2}\, |Dv_t- D^2v H_p(Dv)|^2  + \gamma_H^2 \,|Dv|^{2(q-2)} (\vep + \chi(m))  \, |D^2 v|^2  
 \\ & \leq    \lam\,|v|\, \rho M  
   +  {\rm tr}\left([\mD_\eta {\mathcal A} \cdot \mD z] \,\mD^2 v \right) + {\rm tr}\left([ D_x {\mathcal A} \cdot H_p ]\,\mD^2 v \right) 
%   \\
%   & \qquad\qquad  - \frac{m\chi'(m)}{\chi(m)+\vep} \left(DV \cdot H_p(Dv)\right) {\rm tr}\left(H_{pp} \,D^2 v \right)  
-  \lam \sum\limits_{i,j,\ell=1}^n v \partial_{\eta_\ell} (a_{ij}) \partial_\ell v \,v_{ij}\,.
\end{split}
 \ee
 We now estimate last two terms. By definition of $\mA$, using \rife{m} and that $(f(m)+ \vep \log m)'= (\chi(m)+ \vep) m^{-1}$, we have
\be\label{nasty}
\begin{split}
  \sum\limits_{i,j,\ell=1}^n   \partial_{\eta_\ell} (a_{ij}) \partial_\ell v\, v_{ij} & = 2  \, H_{pp}Dv[D^2v\, H_p(Dv)-Dv_t]  +  \sum\limits_{i,j=1}^d \sum\limits_{\ell=1}^n  \partial_{\eta_\ell} ((\vep+ \chi(m)) H_{p_ip_j}) \partial_\ell v\, v_{ij}
 \\
 & = 2  \, H_{pp}Dv[D^2v\, H_p(Dv)-Dv_t]   + 
 \sum\limits_{i,j=1}^d \frac{m\chi'(m)}{\chi(m)+\vep }  \left(  - v_t + H_p(Dv)\cdot Dv\right) H_{p_ip_j} \, v_{ij} 
 \\
 & + \sum\limits_{i,j=1}^d \sum\limits_{\ell=1}^d  (\vep + \chi(m))    [H_{p_ip_jp_\ell}] \partial_\ell v\, v_{ij}
 \,
 \end{split}
 \ee
 where $H_{pp}, H_{ppp}$ are computed on $Dv$.
 %where $h(m)= f'(m)m$.
%Observe that 
%$$
%\sum\limits_{\ell=1}^n \partial_\ell (f'(m)m ) \partial_\ell v = \partial_r( f'(\vfi(r))\vfi(r))\mathop{|}_{r=f(m)} \left(C_0 - v_t + H_p(Du)\cdot Dv\right)
%$$
%hence
%\begin{align*} & 
%\sum\limits_{i,j=1}^N \sum\limits_{\ell=1}^n  \partial_\ell (f'(m)m H_{p_ip_j}) \partial_\ell v\, v_{ij}=
%\sum\limits_{i,j=1}^N \partial_r( f'(\vfi(\cdot )\vfi(\cdot) ) \left(- v_t + H_p(Du)\cdot Dv\right) H_{p_ip_j} \, v_{ij} 
%\\
%& \qquad + C_0 \sum\limits_{i,j=1}^N \partial_r( f'(\vfi(\cdot )\vfi(\cdot) )) H_{p_ip_j} \, v_{ij}
%+ \sum\limits_{i,j=1}^N \sum\limits_{\ell=1}^n  f'(m)m  \partial_\ell [H_{p_ip_j}] \partial_\ell v\, v_{ij}
%\end{align*}
We observe that assumption \rife{fprime} implies (and, for fixed $\vep$,  is actually equivalent to) $ m \chi'(m) \leq c (\vep + \chi(m))^{\frac32}$ for some constant $c>0$, and for $m>s_0>0$. Let us suppose by now that $m>s_0$ holds true, so that we can use this condition.  Then, using also conditions \rife{Hppq}-\rife{hpppq} to handle  $H_{pp}$ and $H_{ppp}$,
%and that $f'(m)m\in C^1$ and $m$ is bounded above by Lemma \ref{mound},   
we can  estimate the quantity in   \rife{nasty} multiplied by $\lambda v$ as follows:
 \begin{align*}
&   \lam \,\, \left| \sum\limits_{i,j,\ell}  v \partial_{\eta_\ell} (a_{ij}) \partial_\ell v \, v_{ij} \right|  \leq  \lam\, \frac{\gamma_H} 2 (\vep + \chi(m)) | D v|^q + \frac\lam 2 |- v_t + H_p(Dv) \cdot Dv|^2 
  \\
  & \qquad \quad + C\, \lam\, v^2\,  \frac {|Dv|^{q-2}}{\vep + \chi(m)}  |D^2v\, H_p(Dv)-Dv_t|^2+     \lambda v^2\,  C  \, (\vep + \chi(m))|Dv|^{2(q-2)}  |D^2v|^2  
  %+ \lam \, v^2C_0^2\frac{\left(\frac{mh'(m)}{\chi(m)+\vep }\right)^2}{(\vep + \chi(m))} 
  \\
& \leq  \lam\, \frac{\gamma_H} 2 (\vep + \chi(m)) | D v|^q + \frac\lam 2 |- v_t + H_p(Dv) \cdot Dv|^2 \\
& \quad + \sigma  \frac C{\vep } |Dv|^{q-2}\, |D^2v\, H_p(Dv)-Dv_t|^2+    \sigma C\, (\vep + \chi(m))  |Dv|^{2(q-2)} \,|D^2v|^2  \,
%+ \sigma  \, C\, (1+ \|u\|_\infty)^2 \,,
\end{align*}
where we used  the choice of $\lambda$ (see \rife{lamb}). 
Let us stress that, if $\sigma$ is sufficiently small, last two terms  are absorbed in the left-hand side of \rife{eqz}.  Similarly we estimate, using again \rife{fprime},
\begin{align*} {\rm tr}\left([ D_x {\mathcal A} \cdot H_p ]\,\mD^2 v \right) & =  
 -\frac{m\chi'(m)}{\chi(m)+\vep} \left(DV \cdot H_p \right) {\rm tr}\left(H_{pp} \,D^2 v \right) \\
 & \leq  \frac{\gamma_H^2}2(\vep + \chi(m)) |Dv|^{2(q-2)}\, |D^2v|^2 + C\, \|DV\|_\infty^2\, |Dv|^{2(q-1)}\,,
\end{align*}
and
\begin{align*} &  
\lambda v\,  \left(DV\cdot H_p-H_{pp}DV\cdot Dv\right) +   D^2V H_p\cdot H_p   \geq   -  \sigma  - C\,  (\|DV\|_\infty^2+ \|D^2V\|_\infty) \, |Dv|^{2(q-1)}\,.
\end{align*}
%where we used the choice of $\lambda$ (see \rife{lamb}), the value of $C_0$ and the ellipticity   of $\mA$ (see \rife{ella}), where we called  $\alpha$  the ellipticity constant of $\mA$, namely  $\mA (Du) \eta\cdot \eta \geq \alpha |\eta|^2$. 
Putting all together, we choose  $\sigma$ sufficiently small and we deduce from \rife{eqz} 
 \begin{align*} & 
-{\rm tr}\left({\mathcal A} \, {\mathcal D}^2 z \right)+ H_{pp}DV \cdot Dz+ \rho z 
%\\
%& \qquad 
%+ \rho\left( H_p\cdot Dv- H\right) 
%\\
%& 
 +\frac \lam2 \left\{ |-v_t + H_p(Dv) \cdot Dv |^2+ \frac {\gamma_H}2(\vep + \chi(m)) |Dv|^q\right\}
\\ & \quad  
 + \frac{\gamma_H}2|Dv|^{q-2}\, |Dv_t- D^2v H_p(Dv)|^2  + \frac{\gamma_H^2}4 (\vep + \chi(m)) |Dv|^{2(q-2)} \, |D^2 v|^2  
 \\ & \leq   \lam\,|v|\, \rho M  
   +  {\rm tr}\left([\mD_\eta {\mathcal A} \cdot \mD z] \,\mD^2 v \right) 
+ C  +  C\, (\|DV\|_\infty^2+ \|D^2V\|_\infty) \, |Dv|^{2(q-1)}    \,.
\end{align*}
Dropping some positive terms, and noticing that $\lam\,|v|\, \rho M \leq \rho\, C $ (by choice of $\lam, M$), we get
\be\label{preut}
\begin{split} & 
-{\rm tr}\left({\mathcal A} \, {\mathcal D}^2 z \right)+ H_{pp}DV \cdot Dz + \rho z 
%\\
%& \qquad 
%+ \rho\left( H_p\cdot Dv- H\right) 
%\\
%& 
 +\frac \lam2 \left\{ |-v_t + H_p(Dv) \cdot Dv |^2+ \frac {\gamma_H}2(\vep + \chi(m)) |Dv|^q\right\}
%\\ & \quad  
 %+ \frac{\gamma_H}2\, |Dv_t- D^2v H_p(Dv)|^2  + \frac{\gamma_H^2}2 (\vep + \chi(m))  \, |D^2 v|^2  
 \\ & \leq     C\,  \rho   +  {\rm tr}\left([\mD_\eta {\mathcal A} \cdot \mD z] \,\mD^2 v \right) 
+ C\,   +  C\, (\|DV\|_\infty^2+ \|D^2V\|_\infty) \, |Dv|^{2(q-1)}    \,.
\end{split}
\ee
It would be finished if not for the term containing the drift $V(x)$. To handle this part, we need to further exploit the quantity $|-v_t + H_p(Dv) \cdot Dv |^2$. 
In fact, by Lemma \ref{ut} we know that $(u_t)_+$ attains its maximum at   $\Sigma_0 \cup \Sigma_T$. Suppose it holds at $t=0$. Then  we have, for some point $x_0\in \Omega$, 
\begin{align*}
u_t(t,x)  & \leq  u_t(0,x_0)\leq  H(Du(0, x_0))- f(m_0(x_0))-V(x_0)+ \de \|u\|_\infty 
\\ & \leq \max_{\oQ} z + C \left( \|f^\vep(m_0)\|_\infty+ \|V\|_\infty   \right)
\end{align*}
where we used Lemma \ref{sup}. Similarly we reason if the maximum occurs at $t=T$.  Hence,  at any interior maximum point of $z$, we have
$$
H(Dv)-v_t= \max_{\oQ} z- u_t-\frac\lam2 v^2-\frac {C_0}T \geq -K
%|v_t- H(Dv)| \leq C \left( \|f(m_i)\|_\infty+ \|V\|_\infty+ C_0+ \sigma \right)
$$
for some $K$ depending on $\|f^\vep(m_0)\|_\infty, \|f^\vep(m_1)\|_\infty, \|V\|_\infty$.  Notice that this also yields  $H(Du)-u_t- V(x) \geq -K - \|V\|_\infty$ at any interior maximum point, so that $m > \kappa>0$ for some $\kappa>0$ only depending on $\|f^\vep(m_0)\|_\infty, \|f^\vep(m_1)\|_\infty, \|V\|_\infty$. This justifies that, at the interior maximum point,  we used condition \rife{fprime} (where we can assume $s_0<\kappa$ with no loss of generality).   On account of \rife{coercq} we get 
$$
H_p(Dv)\cdot Dv-v_t  \geq \gamma_0\, |Dv|^q - \gamma_1 -  K \,.
%c_0 \, |Dv|^4 - C \left( \|f(m_i)\|_\infty^2+ \|V\|_\infty^2+ C_0^2+ 1 \right)\,.
$$
Therefore, we deduce from \rife{preut}
%
%which implies, using the ellipticity of $\mA$:
% \begin{align*} &
%-{\rm tr}\left({\mathcal A} \, {\mathcal D}^2 z \right)+ \vep z +\lam\, \frac\alpha 2  |\mD v|^2  + |Dv_t- D^2v H_p(Dv)|^2 + mf'(m) \, {\rm tr}  \left(H_{pp}(Dv)D^2vD^2v  \right)
% \leq   \lam\,\vep (M-C_0(1-t))\, v  
%\\ &  -  {\rm tr}\left([\mD {\mathcal A} \cdot \mD z] \,\mD^2 v \right)  +     \sigma\, C\left(  |D^2v\, H_p-Dv_t|^2+|D^2v|^2\right) \,.
%\end{align*}
%Now we choose $\sigma$ sufficiently small (only depending on the ellipticity of $f'(s)s\, H_{pp}$) and last term is absorbed in the left-hand side.  We conclude that
\begin{align*}
& -{\rm tr}\left({\mathcal A} \, {\mathcal D}^2 z \right)+ H_{pp}DV \cdot Dz -    {\rm tr}\left([\mD_\eta {\mathcal A} \cdot \mD z] \,\mD^2 v \right) + \rho z 
+ c_0 \lambda |Dv|^{2q} \\
& 
\leq   \tilde K +     C\, (\|DV\|_\infty^2+ \|D^2V\|_\infty) \, |Dv|^{2(q-1)} \,,
\end{align*}
for some $c_0>0$ and some $\tilde K$ depending on   $\|f^\vep(m_0)\|_\infty, \|f^\vep(m_1)\|_\infty, \|V\|_{\infty}$.  
Recalling the choice of $\lambda$, we conclude that, on any interior maximum point of $z$, we have
$$
|Dv| \leq C (1+ \|u\|_\infty)
$$
for some $C$ depending on $\vep, \|V\|_{W^{2,\infty}}, \|f^\vep(m_i)\|_{W^{1,\infty}}$ and the functions $H, f$.
\vskip1em

Step 4 ({\it conclusion}). On account of the previous steps, we estimate the maximum of $z$ and in turn this yields
$$
\|Dv\|_\infty  \leq C \,,\qquad \hbox{where $C=C(\vep, \|V\|_{W^{2,\infty}}, \|f^\vep(m_0)\|_{W^{1,\infty}}, \|f^\vep(m_1)\|_{W^{1,\infty}}, \|u\|_\infty)$.}
$$
This estimate, together with Lemma \ref{sup}, imply a similar bound for $u_t$ at $t=0,t=T$.  Using Lemma \ref{ut}, we conclude with a full estimate for $|u_t|$. Hence, the desired estimate follows for $\mD u$.
\qed

\vskip1em
We  now give a further  application  of Lemma \ref{struc}  to  estimate the density $m$ in the optimal transport problem.
Indeed, the next  lemma gives sufficient conditions under which  $m$, defined from \rife{m} as $ \vfi (-u_t + H(Du)-V(x) )$, is bounded above.

 \begin{lemma}\label{mound} Let $u\in C^2(Q)\cap C^1(\overline Q)$ be a solution of \rife{eps-de}. Assume that $V\in W^{2,\infty}(\Omega)$ is convex, and that  $f'(m)m$ is nondecreasing.  Then, there exists a constant $C$ such that
$$
 -u_t  + H(Du) \leq C  \,,
$$
where $C= C(\|m_0\|_\infty, \|m_1\|_\infty, \|V\|_\infty, \de \| u\|_\infty, H,f)$.
\end{lemma}

\proof  We follow the computations  of  \rife{conv} with the convex function $\psi(s,\xi)= -s + H(\xi)$, but we specialize to the matrix $\mathcal A$ given by \rife{aij}.  This means that $w:= -u_t + H(Du)$ satisfies
\begin{align*}
& -{\rm tr}\left({\mathcal A} \, {\mathcal D}^2 w \right)   + \rho (-u_t +H_p(Du)Du) +  H_{pp}DV\cdot Dw \leq  - D^2 V H_p(Du)\cdot H_p(Du) 
\\
& \qquad   +  {\rm tr}\left([\mD_\eta {\mathcal A_1} \cdot \mD w] \,\mD^2 u \right)  +  (\vep + mf'(m)) H_{p_ip_jp_\ell} w_\ell \,u_{ij} 
- (mf'(m))' (m_t- H_p(Du)\cdot Dm) \hbox{tr}(H_{pp}D^2 u)
\end{align*}
where $H_{pp}, H_{ppp}$ are always computed on $Du$.
Since $m_t- H_p(Du)\cdot Dm= m\,\hbox{tr}(H_{pp}D^2 u)$, using the convexity of $V$ and the assumption that $mf'(m)$ is nondecreasing we have
$$
D^2 V H_p(Du)\cdot H_p(Du) + (mf'(m))' (m_t- H_p(Du)\cdot Dm) \hbox{tr}(H_{pp}D^2 u)\geq 0\,.
$$
Hence we deduce that $w$ satisfies
$$
-{\rm tr}\left({\mathcal A} \, {\mathcal D}^2 w \right)   + \rho (-u_t +H_p(Du)Du) +  H_{pp}DV\cdot Dw \leq     {\rm tr}\left([\mD_\eta  {\mathcal A_1} \cdot \mD w] \,\mD^2 u \right)  +  (\vep + mf'(m)) H_{p_ip_jp_\ell} w_\ell \,u_{ij}\,. 
$$
Since \rife{coercq} implies $-u_t +H_p(Du)Du\geq w-\gamma_1$, we deduce that
$$
-{\rm tr}\left({\mathcal A} \, {\mathcal D}^2 w \right) + \rho( w-\gamma_1) +  H_{pp}DV\cdot Dw \leq      {\rm tr}\left([\mD_\eta  {\mathcal A_1} \cdot \mD w] \,\mD^2 u \right)  +  (\vep + mf'(m)) H_{p_ip_jp_\ell} w_\ell \,u_{ij}\,. 
$$
This implies that $w\leq \gamma_1$ on any   maximum point inside the domain. Using the bound   at $t=0,t=T$ we have
$$
w \leq C  \left( 1+ \|m_0\|_\infty + \|m_1\|_\infty + \|V\|_\infty\right)+ \de \|u\|_\infty\,.
$$
On the boundary $\partial \Omega$, we have (with same arguments as in the previous propositions) that $Dw \cdot  \vec\nu\leq 0$ and no maximum point can occur. Finally, we conclude that $w$ is bounded above by a constant only depending on $\gamma_1, \| m_0\|_\infty$, $\|m_1\|_\infty, \|V\|_\infty$, $\de\|u\|_\infty$. 
\qed

\vskip2em

Finally, the estimates of the present Section, and specifically the gradient bound of  Theorem \ref{bound}, lead to the  existence of solutions to \rife{eps-de}.  The proof of the following result, based on a continuity method,  is quite  technical and we postpone it to the Appendix.

\begin{theorem}\label{lieb}  Let  $\rho, \de>0$. Assume that $f\in C^2(0,\infty)$ satisfies \rife{fprime}, and that $H\in C^3(\R^d)$ satisfies \rife{hrad} and the growth conditions \rife{Hppq}-\rife{hpppq} (for some $q>1$ and $\varpi>0$). Let $V\in W^{2,\infty}(\Omega)$ and $m_0, m_1 \in W^{1,\infty}(\Omega)$ with $m_0, m_1 >0$ in $\oo$. Then  there exists  a  solution $u\in C^{2,\alpha}(Q)\cap C^{1,\alpha}(\oQ)$ to problem \rife{eps-de} and we have
\be\label{lastpray}
\| u\|_{C^{1,\alpha}(\oQ)} \leq K
\ee
for some $K= K(\vep, \|V\|_{W^{2,\infty}}, \|f^\vep(m_0)\|_{W^{1,\infty}}, \|f^\vep(m_1)\|_{W^{1,\infty}}, \|u\|_\infty)$ (independent of $\rho, \de$). 
\end{theorem}

\section{A penalized optimal transport problem}\label{pena}

At this stage we wish to  let the parameters $\rho,\de$ in \rife{eps-de} tend  to zero: they were only needed in the existence result to control the sup norm of $u$. First of all, by letting $\rho \to 0$, we find  a smooth solution of a control/transport problem with penalized initial-terminal conditions, which may have an interest in its own.

\begin{theorem}\label{delta} Let $f,H, V, m_0, m_1 $ satisfy the same conditions of Theorem \ref{lieb}.   Then, for any $\de>0$ there exists    $u\in C^{2,\alpha}(Q)\cap C^{1,\alpha}(\oQ)$ which is a solution  to the problem
\be\label{pb-delta}
\begin{system}
& -{\rm tr}\left({\mathcal A}({\mathcal D}u) \, {\mathcal D}^2u \right) + DV(x)\cdot H_p(Du)=0 &&\text{in } Q, \\%(0,1)\times \R^d, \\
& -u_t + H(Du)= \de u+ f^\vep(m_1)+ V(x)  &&\hbox{at $t=T$, $x\in \Omega$,}
\\
&   -u_t + H(Du)+ \de u= f^\vep(m_0) + V(x) &&\hbox{at $t=0$, $x\in \Omega$,}
\\
& \,\, H_p(Du)\cdot \vec\nu=0 &&\text{on $(0,T)\times \partial \Omega$}
\end{system}
\ee
where $\mA$ is given by \rife{aij}. 

In particular, if we set $m:= (f^\vep)^{-1}(-u_t + H(Du)-V(x))$, then $u\in C^{2,\alpha}(Q)\cap C^{1,\alpha}(\oQ)$ is a  solution of \rife{pb-delta} if and only if the couple $(u,m)$ is a    solution of the penalized   problem
\be\label{MFG-delta}
\begin{system}
&-u_t + H( Du)  = f(m) + \vep \log m + V(x)&&\text{in } Q, \\%(0,1)\times \R^d, \\
&m_t  - \dive (m\,H_p(Du))  = 0  &&\text{in } Q,
\\
& \,\, H_p(Du)\cdot \vec\nu=0 &&\text{on $(0,T)\times \partial \Omega$,}
\\%(0,1)\times \R^d,\\
&u(0) = \frac{f^\vep (m_0)- f^\vep (m(0))}\de&&\text{in } \Omega,
\\
& u(T)  = \frac{f^\vep (m(T))-f^\vep (m_1)} \de&&\text{in } \Omega,
\end{system}
\ee
and  $m\in C^{1,\alpha}(Q) \cap C^{0,\alpha}(\oQ)$, $m>0$ in $\oQ$.
\end{theorem}

\proof Let $u_{\rho,\de}$ be the solution to \rife{eps-de} given by Theorem \ref{lieb}.  As a consequence of Lemma \ref{sup} and Theorem \ref{bound}, $u_{\rho,\de}$ is bounded in $W^{1,\infty}(Q)$ uniformly with respect to $\rho$. We deduce from Theorem \ref{lieb} that $u_{\rho,\de}$ is also bounded in $C^{1,\alpha}(\bar Q)$ for some $\alpha>0$. This also yields that $u_{\rho,\de}$ is bounded in $C^{2,\alpha}(Q)$ by classical Schauder's elliptic regularity. As $\rho\to 0$, we deduce the existence of a  smooth solution $u$ to \rife{pb-delta}. Defining $m$ as in \rife{m}, we have that $m\in C^{0,\alpha}(\bar Q)\cap C^{1,\alpha}(Q)$, and in addition $m>0$ because $-u_t + H(Du)$ is bounded and  $\vfi$ maps compact sets into compact sets of $(0,\infty)$. The equivalence between the elliptic equation \rife{pb-delta} and the system \rife{MFG-delta} is  straightforward for smooth solutions and the nonlinear Neumann condition for $u$ at $t=0,T$ turns into the initial-terminal conditions  
of system \rife{MFG-delta}.
\qed

\vskip1em

The next step will consist in the limit as $\de\to 0$. Here the main issue is  to normalize  $u$ in order to stabilize the sup norm.  In the next result we  show that a control on the $L^\infty$-norm of the density $m$ yields  a control on a natural normalization of  $u$.
%To this purpose,   we exploit the duality between sub solutions of Hamilton-Jacobi equations and geodesics  of the Wasserstein distance.

\begin{lemma}\label{norm} Assume that $m_0, m_1\in \elle\infty$, and that $(u,m)$ is a (smooth) solution of \rife{mfpp} or, alternatively,  a solution of \rife{MFG-delta} with  $\de>0$.
%
%
%t least one of the two following conditions hold:
%
%(i) $H$ satisfies \rife{coerc}-\rife{Hpp}.
%
%(ii) $V$ is convex and $f'(m)m$ is nondecreasing.
If we set
\be\label{normou}
\hat u := u - \into u (T)\,m_1\, dx
\ee
then  $\hat u $ satisfies the global estimate
\be\label{hatu1}
\| \hat u\|_\infty \leq C(\| m \|_\infty , \|  \vep \log (m_1) \|_\infty, \|  \vep \log (m_0) \|_\infty,  \| V \|_{W^{1,\infty}}, T, H, \Omega) 
% \| (m_0e^{\frac V\vep})^{-1}\|_\infty, \|\|m_0e^{\frac V\vep}\|_\infty) 
\ee
and the local estimate
\be\label{hatu2}
- \frac {C_0}{t^{\frac1{q-1}}} \leq \hat u(t,x) \leq \frac {C_1}{(T-t)^{\frac1{q-1}}} \qquad \forall (t,x)\in Q
\ee
where $C_0$, $C_1$ only depend on $\sup(f^\vep (m)+ V)$.
%In particularis bounded in $W^{1,\infty}(Q)$  independently of $\de$.
\end{lemma}

\proof  
%First observe that 
%$$
%\| m_\de(0)\|_\infty \leq  (f^\vep)^{-1} ( \de\|u_\de\|_\infty  + f^\vep(\|m_0\|_\infty)
%$$
%and then, due to the bound of $\de\|u_\de\|_\infty$, we have that $\| m_\de(0)\|_\infty$ is bounded. Similarly we have for $\| m_\de(T)\|_\infty$. Therefore, 
%if condition (i) is satisfied, we have that  $m_\de$ is uniformly bounded above by Proposition \ref{supm}.  Otherwise, if (ii) is satisfied, $m_\de$ is bounded above directly from Lemma \ref{mound}. 
We consider the  product between $u(t)$ and $m(t)$; using  \rife{mfpp} (or alternatively \rife{MFG-delta}) we have
$$
\into u(T) m(T) dx -\into u(0)m(0)dx  =  - \int_0^T \into m[H_p(Du)\cdot Du- H(Du)]dxdt - \int_0^T \into (f^\vep (m)+V)m\,dxdt\,,
$$
which implies, by convexity of $H$ and since $f^\vep(s)s$ is bounded below,  
$$
\into u(T) m(T) dx -\into u(0)m(0)dx  \leq C \,.
$$
Here and below all constants are independent of $\de$ and $\vep $ (for $\vep\leq 1$).  This implies
\begin{align*}
 \into u(T)(m (T) -m_1) dx -\into u (0)(m (0)-m_0) dx & 
  \leq C  + \into u (0) m_0 dx- \into u (T)  m_1  dx 
\\
& = C+ \into \hat u (0) m_0 dx
\end{align*}
where we used the definition of $\hat u $ and the fact that   $m_0$ has unit mass.  The left-hand side vanishes if $(u,m)$ is a solution of \rife{mfpp}. Alternatively, if $(u,m)$ is a solution of \rife{MFG-delta}, we have 
\be\label{margi}
\frac1\de \into f^\vep(m(T))-f^\vep(m_1))(m (T) -m_1) dx +\frac1\de \into f^\vep(m (0))-f^\vep(m_0))(m (0)-m_0) dx
\leq C+ \into \hat u (0) m_0 dx\,,
\ee
and the left-hand side is nonnegative by monotonicity of $f^\vep$. Therefore, for both problems \rife{mfpp} and \rife{MFG-delta} we get
\be\label{at0}
\into \hat u (0) m_0 dx \geq -C\,.
\ee 
For any given measure $\tilde m$, let us now consider the $p$-Wasserstein geodesic connecting $m_1$ and $\tilde m$ in $[s,T]$,  specifically we consider the solution $\mu$ to the continuity equation
$$
\begin{system}
& \mu_t  - \dive (\mu\,v)  = 0  &&\text{in } (s,T)\times \Omega,
\\
& v\cdot \vec\nu=0 && \text{in} (s,T)\times \partial \Omega 
\\
& \mu(T)  =m_1\,, \mu(s)=\tilde m &&\text{in } \Omega,
\end{system}
$$
which satisfies 
$$
 \int_s^T \into  |v|^{p}\,d\mu = \frac {c_p}{(T-s)^{p-1}}W_p(m_1, \tilde m)^p\,.
$$
Here we use this geodesic with $p=q/(q-1)$, where $q$ is given by the growth of the Hamiltonian  \rife{Hppq}.  Let us set
$M:= \sup_Q [f^\vep (m)+ V]$. Since $\hat u $ satisfies the inequality $-  \hat u_t  + H(D\hat u )\leq M$, multiplying   by $\mu$ and integrating we get
$$
\into \hat u (s) \tilde m\, dx \leq    \int_s^T \into v\cdot D\hat u  \, \mu\, dxdt - \int_s^T \into H(D\hat u ) \, \mu\, dx dt  + M (T-s)
$$
where we used that $\into \hat u \, m_1\, dx =0$. Using the coercivity of $H$ we deduce
$$
\into \hat u (s) \tilde m\, dx \leq  C  (M(T-s)+  \frac 1{(T-s)^{\frac1{q-1}}}W_{q'}(m_1, \tilde m)^{q'} )\,.
$$
Since this holds for any measure $\tilde m$, we deduce that there exists a constant $C_1$ (only depending on $M, T$), such that
\be\label{above}
\hat u (s,x) \leq \frac {C_1}{(T-s)^{\frac1{q-1}}}  \,. 
\ee
Reasoning in  a similar way, namely using the geodesic between $m_0$ and any measure $\tilde m$, and the bound \rife{at0}, we conclude that there exists a constant $C_0$  such that
\be\label{below}
\hat u (s,x) \geq - \frac {C_0}{s^{\frac1{q-1}}} \,. 
\ee
This concludes the proof of \rife{hatu2}.   Now we transform the previous two bounds into a global $L^\infty$-bound for $\hat u $, by using the strict positivity of the marginals. 
%Notice that  a local bound in $(0,1)$ could be deduced as above by using the geodesics rescaled in time $(0,s)$ or $(s,1)$, however in what follows  we directly obtain a global bound which stems from \rife{above}-\rife{below}.
To this purpose, we consider  $v:= \hat u  -At+\theta (t-T)^2+ \tilde \theta d(x)$, for small $\theta, \tilde \theta$. Then one can see from  \rife{pb-delta}
%$k_0, k_1$ can be chosen sufficiently large so that 
that $v$ cannot have positive maximum  inside $Q$,  nor on $\partial\Omega$, up to choosing $\tilde \theta = o(\theta)$. 
We claim that, if $m_1>0$ in $\oo$, then  $v$ cannot have maximum at $t=T$ either, provided $A$ is sufficiently large. Indeed, if $v$ has a maximum at $t=T$, then  we must have 
$$
0 \leq  \hat u_t (T,x) -A = H(D\hat u ) - \de u - V(x) - f^\vep (m_1) - A \leq  \de\|u\|_\infty+ \|V\|_\infty - f^\vep (\min m_1) - A 
$$
and we get a contradiction if $A$ is large enough (we use Lemma \ref{sup} if $\de>0$). Then  we deduce that $v$ attains its maximum at $t=0$; using  \rife{above} and letting $\theta, \tilde \theta\to 0$, we deduce that   
%$$
%\hat u  (t,x)\leq C
%%+ k_0 (t-T)^2+ k_1 d(x) - At \leq \max_\Omega  \hat u  (x,0) + k_0 T^2 + c\, k_1   \leq C 
%$$
%by \rife{above}.   This means that there exists $K_0$, independent of $\de$, such that
$$
\hat u  \leq K_0 \qquad \forall (t,x)\in \bar Q\,
$$
for some $K_0$, independent of $\de$.
 Similarly we reason with the estimate from below provided $m_0$ is strictly positive (so that $f^\vep (m_0)$ is bounded below).
%We conclude
% Using \rife{below} we can prove that there exists $K_1$, independent of $\de$, such that
%$$
%\hat u_\de(t,x) \geq -  K_1 \qquad \forall (t,x)\in \bar Q\,.
%$$
This concludes with the $L^\infty$ bound \rife{hatu1}.  
%Finally, once $\|\hat u_\de\|_\infty$ is bounded independently of $\de$, we can use Theorem \ref{bound} in order to get  a   bound for $\|\cD \hat u_\de\|_\infty$ as well. It is enough to notice that $\hat u_\de$ satisfies the same equation as $u_\de$ and the same boundary conditions except for adding (or subtracting) to  $f(m_i)$ the   term $\de\into u_\de(T)m_1\,dx$ which is  a bounded constant due to Lemma \ref{sup}.
\qed

\vskip1em

Thanks to Lemma \ref{norm}, the solution $u_\de$ of \rife{pb-delta} will be controlled in sup-norm after the normalization \rife{normou} provided $m_\de$ is bounded in $L^\infty(Q)$.  Sufficient conditions for this latter bound were already given in Lemma \ref{mound}. In the next Section, we see an alternative way of controlling the density through displacement convexity arguments. This is the reason why we postpone to Section \ref{exi} the final convergence $\de\to 0$ in \rife{MFG-delta}.

\section{Displacement convexity and density estimates}\label{disconv}

One of the main advantages in the construction of smooth solutions is the possibility to use the so-called Eulerian approach to displacement convexity estimates (see \cite[Chapter 5.4]{Vi}).  In the context of mean-field game systems, general displacement convexity inequalities were first addressed in \cite{GoSe} in the periodic case.  
In this Section we extend  their result (and later developments in \cite{Gomes2}) to the case of  Neumann boundary conditions and to general Lipschitz potentials $V$, and then we use it to get local and global bounds for the density.
We start with an adaptation of \cite[Thm 1.1]{GoSe}.  Let us stress that in this Section we don't use the radial structure of the Hamiltonian.

\begin{proposition}\label{displace}
Let $\Omega$ be a smooth convex bounded domain in $\R^d$. Let $u\in C^2(\oQ), m\in C^1(\oQ)$ be classical solutions to the system
\be\label{solosys}
\begin{cases}
 -u_t + H( Du)  = f (m) + V(x)&\text{in } Q, \\%(0,1)\times \R^d, \\
 m_t  - \dive (m\,H_p(Du))  = 0  &\text{in } Q, \\
 H_p(Du)\cdot \vec \nu=0 & \hbox{on $(0,T)\times \partial \Omega$}
 \end{cases}
 \ee
 where $f\in C^1(0,\infty)$, $V\in W^{1,\infty}(\Omega)$,   $H\in C^3$.  
 \vskip0.3em
\noindent  Let  $U:(0,\infty) \to \R$ be a $C^1$ function such that 
 $$
 P(r):= U'(r)r- U(r)  \geq 0\,.
 $$
Then  we have
\be\label{displa}
\begin{split}
\frac{d^2}{dt^2}\into U(m(t)) & \geq \into \left(P'(m)m- P(m)+ \frac1d P(m)\right) [\dive(H_p(Du))]^2 
\\
& + \into P'(m)f'(m) (H_{pp}(Du)Dm\cdot Dm) + \into P'(m) (H_{pp}(Du)Dm\cdot DV(x))
\end{split}
\ee 
\end{proposition}

\proof We borrow most of the  computations from \cite[Thm 1.1]{GoSe}, while we take care of  two new ingredients, i.e. the term with $V(x)$ and the boundary conditions. For the reader's convenience, we write all the required steps. First of all, using the continuity equation we have
\begin{align*}
\frac{d}{dt}\into U(m(t)) & = \into U'(m) m_t = \into U'(m) \dive(m H_p(Du))  \\
& 
= \into U'(m)m \dive(H_p(Du)) + \into DU(m) \cdot H_p(Du) \\
& = \into P(m) \dive(H_p(Du)) 
\end{align*}
where we used the Neumann condition (in the last step) and the definition of $P(r)= U'(r)r-U(r)$. Hence, using again the equation of $m$ and integrating by parts,
\be\label{half}
\begin{split}
& \frac{d^2}{dt^2}\into U(m(t))  = \into P'(m) [m (\dive(H_p(Du)) )+ Dm\cdot H_p(Du)] \dive(H_p(Du))  + \into P(m)    \dive(H_{pp}Du_t) 
\\
& = \into P'(m) m (\dive(H_p(Du)) )^2 - \into P(m) \dive [H_p(Du) \dive(H_p(Du)) ] + \into P(m)    \dive(H_{pp}Du_t) 
\\
& = \into (P'(m) m-P(m))  (\dive(H_p(Du)) )^2 - \into P(m)  H_p(Du) \cdot D[\dive(H_p(Du)) ] + \into P(m)    \dive(H_{pp}Du_t) 
\end{split}
\ee
where $H_{pp}$ is computed on $Du$. 
Now we use the equation of $u$ in the last term:
\be\label{half2}
\begin{split}
\into P(m) \dive(H_{pp}Du_t) & = \into P(m) \dive(H_{pp}D[u_t-H(Du)])  + \into P(m) \dive(H_{pp}DH(Du))
\\ & = \int_{\partial\Omega} P(m) D[u_t-H(Du)]\cdot H_{pp}\vec\nu - \into P'(m) H_{pp}Dm\cdot D[u_t- H(Du)]
\\ & \qquad + \into P(m) \dive(H_{pp}DH(Du))
\\ & = \int_{\partial\Omega} P(m) D[u_t-H(Du)]\cdot H_{pp}\vec\nu + \into P'(m) H_{pp}Dm\cdot D[f(m)+ V(x)]
\\ & \qquad + \into P(m) \dive(H_{pp}DH(Du))\,.
\end{split}
\ee
We treat the boundary term   writing $\vec\nu= - Dd(x)$, and using the convexity of $\Omega$ (which implies that $d(x)$ is concave).
We have
\begin{align*}
DH(Du)\cdot H_{pp}(Du)Dd & =   H_p(Du) \cdot D^2 u H_{pp}(Du)Dd \\
& =  H_p(Du) \cdot D (H_p(Du) \cdot Dd) - D^2 d H_{p}(Du)H_p(Du) \\
& = - D^2 d H_{p}(Du)H_p(Du) \geq 0
\end{align*}
where we used that $H_p(Du)$ is tangential and $H_p(Du) \cdot Dd$ has zero tangential gradient. Since we also have  
$$
Du_t\cdot H_{pp}(Du)Dd= \partial_t (H_p(Du)\cdot Dd)=0
$$
we conclude that 
$$
\int_{\partial\Omega} P(m) D[u_t-H(Du)]\cdot H_{pp}\vec\nu \geq 0\,.
$$
 Using this information in \rife{half2} and inserting that in \rife{half} we get
 \be\label{half3}
\begin{split}
& \frac{d^2}{dt^2}\into U(m(t))  \geq  \into (P'(m) m-P(m))  (\dive(H_p(Du)) )^2  + \into P'(m) H_{pp}Dm\cdot D[f(m)+ V(x)] \\
& + \into P(m) \dive(H_{pp}DH(Du))- \into P(m)  H_p(Du) \cdot D[\dive(H_p(Du)) ]  \,.
\end{split}
\ee
We develop jointly last two terms observing 
%(with Einstein's convention on repeated indices)
\begin{align*}
& \dive(H_{pp}DH(Du))-  H_p(Du) \cdot D[\dive(H_p(Du)) ] 
%\\
%& 
=   H_{p_ip_k}H_{p_\ell p_j} u_{ij}u_{\ell k} = {\rm Tr}([H_{pp}D^2u]^2)
\end{align*}
Using that ${\rm Tr}([AB]^2)\geq \frac 1d [{\rm Tr}(AB)]^2$ for $A,B$ symmetric with $A\geq 0$ (see e.g. \cite[Lemma A.1]{GoSe}) we conclude that
$$
\dive(H_{pp}DH(Du))-  H_p(Du) \cdot D[\dive(H_p(Du)) ] \geq  \frac1d(\dive(H_p(Du)) )^2
$$
and then from \rife{half3} we get
$$
 \frac{d^2}{dt^2}\into U(m(t))  \geq  \into \left(P'(m) m-\left(1- \frac1d\right) P(m)\right)  (\dive(H_p(Du)) )^2  + \into P'(m) H_{pp}Dm\cdot D[f(m)+ V(x)] 
$$
which is \rife{displa}.
\qed

\vskip1em
We now deduce a priori estimates from the displacement convexity inequality \rife{displa}. To this purpose, we require the Hamiltonian $H$ to be of quadratic type, namely satisfying condition \rife{Hpp}.
While the global (in time) estimates are the  typical application of displacement convexity (they are obtained e.g. in \cite{GoSe}, \cite{Gomes2}), the local estimates are more interesting, since they show a  regularizing effect due to the congestion term in the cost functional of optimal transport. Similar local estimates were obtained in \cite{LaSa} for quadratic Hamiltonian with a very tricky use of flow-interchange techniques together with the variational interpretation of the functional as geodesic in Wasserstein space.  Our next Proposition  provides an alternative, different proof of the results in \cite{LaSa}, in a slightly  broader setting. 

\begin{proposition}\label{supm}
Assume that $V\in W^{1,\infty}(\Omega)$  and $f\in C^1(0,\infty)$   satisfies
\be\label{coercf}
\exists\,\, c_0, \, r_0>0\,:\qquad f'(r) \geq c_0 r^{-1} \qquad\forall r\geq r_0\,.
\ee
Let $(u,m)$ satisfy \rife{solosys}, where $H$ is a $C^2$ function satisfying \rife{Hpp}. Then, for every $1\leq p \leq \infty$ we have the global  estimate
\be\label{supglobm}
\| m(t)\|_p \leq K_0\, (\|m_0\|_p+ \|m_1\|_p+1) \qquad \forall t\in [0,T]
\ee 
and the local estimate
\be\label{suplocm}
\| m(t)\|_p \leq K_1\, \left(t^{-q}+ (T-t)^{-q} \right) \qquad \forall t\in (0,T)
\ee
for some constants $K_0, K_1$ depending on $c_0, r_0, p, d, T, |\Omega|, \|DV\|_\infty, \alpha_H, \beta_H$, and for some $q>0$.
\end{proposition}

\proof  The proof is done in two steps. We first establish global and local estimates in $L^p$-norm, for $p<\infty$, and then we ugrade the estimates to the sup-norm.

%{\it Step 0.}\quad We use $U(r)= r\int_0^r \frac1{\sqrt{1+s^2}+1}ds$ in \rife{displa}. This yields $P(r)= \sqrt{1+r^2}-1$ and we get
%$$
%\alpha_H\into \frac{f'(m)m}{\sqrt{1+m^2}} |Dm|^2 + \frac1d\into m [\dive(H_p(Du))]^2 + \into P'(m)Dm\cdot H_{pp}DV \leq  \frac{d^2}{dt^2}\into U(m(t))
%$$
%which implies
%$$
%\frac{\alpha_H}2\into \frac{f'(m)m}{\sqrt{1+m^2}} |Dm|^2 \leq \frac{d^2}{dt^2}\into U(m(t)) + c_H \into \frac{P'(m)}{f'(m)}|V|^2
%$$
%Using that $f$ is increasing and satisfies \rife{coercf}, that $P'\leq 1$ and that $U(m)\leq C\, m\log m$, we deduce
%$$
% \into  \frac{f'(m)m}{\sqrt{1+m^2}} |Dm|^2 \leq C + \frac{d^2}{dt^2}\into U(m(t))\,.
%$$
{\it Step 1.}\quad 
We use \rife{displa} with $U(r)= (r-r_0)_+^p$, $p\geq 2$, where $r_0$ is given by \rife{coercf}. We notice that $U\in C^1\cap W^{2,\infty}$ and we have $P(r)\geq 0$ and  $P'(r)= U''(r)r= p(p-1)(r-r_0)_+^{p-2}r \mathds1_{\{r>r_0\}}$. Moreover, we have $P'(r)r-P(r)\geq 0$. We deduce that 
\begin{align*}
\frac{d^2}{dt^2}\into (m(t)-r_0)_+^p & \geq  p(p-1) \left\{ \into  (m-r_0)_+^{p-2}\mathds1_{\{m>r_0\}} \left[ f'(m)m (H_{pp}Dm\cdot Dm) + m (H_{pp}Dm\cdot DV(x))\right] \right\}
\\
& \geq p(p-1) \left\{ c_0\frac{\alpha_H}2 \into  (m-r_0)_+^{p-2}\mathds1_{\{m>r_0\}} |Dm|^2  -  c_H\, \|DV\|_\infty^2 \left(  \into  (m-r_0)_+^{p} +  r_0^p |\Omega|\right)\right\}
\end{align*}
where we used assumption \rife{coercf}  and that $\alpha_H Id\leq H_{pp} \leq \beta_H Id$.  Setting $\mu:= (m-r_0)_+$, we rephrase the above inequality as
\be\label{start}
\frac{d^2}{dt^2}\into \mu(t)^p  \geq p(p-1) \left\{ c_1 \into  \mu^{p-2}|D\mu|^2  -  c_2\, \|DV\|_\infty^2 \left(\into \mu^{p} +1\right) \right\}\,.
\ee
%We use \rife{displa} with $U(r)= r^p$, $p>1$, hence $P(r)= (p-1)r^p$;  we deduce that 
%\be\label{start}
%\begin{split}
%\frac{d^2}{dt^2}\into m(t)^p & \geq  p(p-1) \left\{ \into  m^{p-1}f'(m) (H_{pp}Dm\cdot Dm) + \into m^{p-1} (H_{pp}Dm\cdot DV(x))\right\}
%\\
%& \geq p(p-1) \left\{ c_0\frac{\gamma_H}2 \into  m^{p-2}|Dm|^2  -  c_H\, \|DV\|_\infty^2 \into m^{p} -  c\right\}
%\end{split}
%\ee
%where we used assumption \rife{coercf}  and that $\gamma_H Id\leq H_{pp} \leq \beta_H Id$.  We notice that there is no loss of generality in assuming that $r_0=1$ in \rife{coercf} and that, as soon as $p\geq 3$, we have 
%$$
%\into |Dm|^2 m^{p-2}  \mathds1_{\{m<1\}}\leq c \into m |Dm|^2 \mathds1_{\{m<1\}}\leq c
%$$
%as it follows using $U(r)= \sqrt{1+r^2}-1$ (hence $P= 1- (1+ r^2)^{-1/2}$) in \rife{displa} and the $L^1$ bound on $m$. 
By Sobolev and Poincar\'e-Wirtinger inequality we have (for $2^*= \frac {2d}{d-2}$ if $d>2$, or $2^*$ any sufficiently large number if $d= 2$)
\begin{align*}
\frac{p^2}4\into  \mu^{p-2}|D\mu|^2   & =  \into |D \mu^{\frac p2}|^2 \geq   C_S \left(\into |\mu^{\frac p2}- (|\Omega|^{-1}\into \mu^{\frac p2}) |^{2^*}\right)^{\frac 2{2^*}}
\\
& \geq c_3 \left(\into \mu^{p\frac{2^*}2}\right)^{\frac2{2^*}} - c_4 \left(|\Omega|^{-1}\into \mu^{\frac p2}\right)^2
\\ & \geq c_3 \left(\into \mu^{p\frac{2^*}2}\right)^{\frac2{2^*}} - c_4  |\Omega|^{-1}\into \mu^p
\end{align*}
while using the interpolation inequality for $1<p< p\frac{2^*}2$ and  $\|m(t)\|_{\elle1}=1$, we have
$$
\into \mu^p \leq     \left(\into \mu^{p\frac{2^*}2}\right)^{\frac{p-1}{p\frac{2^*}2-1} }\,.  
$$
Using the exact value of $2^*$ we conclude that
$$
\frac{p^2}4\into  \mu^{p-2}|D\mu|^2 \geq c_3\left( \into \mu^p\right)^{1+ \frac{2}{d(p-1)}} - c_4  |\Omega|^{-1}\into \mu^p\,.
$$
From \rife{start} we deduce (using e.g. $p(p-1)\geq p^2/2$ for $p\ge 2  $)
$$
\frac{d^2}{dt^2}\into \mu(t)^p \geq    c_5 \left( \into \mu(t)^p\right)^{1+ \frac{2}{d(p-1)}}-  c_6(1+ p^2 \|DV\|_\infty^2)\into \mu(t)^{p} - c_7\, p^2\|DV\|_\infty^2
$$
which implies
\be\label{ode}
-\vfi''+  \, c_5 \vfi^{1+\frac 2{d(p-1)}} - c_8 \vfi \leq c_7 \, p^2 \|DV\|_\infty^2
\ee
for $c_8= c_6(1+ p^2 \|DV\|_\infty^2)$ and $\vfi(t)= \into \mu(t)^p$. Now we apply the maximum principle to $\vfi$ and we get
$$
\vfi(t) \leq \max \left\{ \vfi(0), \vfi(T),  C_p\right\}
$$ 
for some constant $C_p$ depending on $c_0, r_0, p, d, \|DV\|_\infty, \alpha_H, \beta_H$. This yields the global estimate in $L^p$-norm
\be\label{globp}
\|m(t)\|_p \leq    \|m_0\|_p+ \|m_1\|_p+ C \left( c_p\|DV\|_\infty+1 \right)\,.
\ee
But we also deduce a  local estimate from \rife{ode}, due to the superlinear term.  In fact, the function 
$$
\bar \vfi(t):= L\left(t^{-\alpha}+ (T-t)^{-\alpha}\right) \,,\qquad \alpha:= d(p-1)
$$
is a supersolution of \rife{ode} for $L$ sufficiently large (possibly depending on $p$). By comparison between $\vfi$ and $\bar \vfi$ we deduce the local estimate:
\be\label{locp}
\|m(t)\|_p \leq     K  \left( t^{-\frac d{p'}}+ (T-t)^{-\frac d{p'}} \right)\, 
\ee
for some $K$ depending on $p, d, \|DV\|_\infty, T, r_0, \alpha_H, \beta_H$.
\vskip0.5em
{\it Step 2.} \quad We point out that, if not for the presence of the field $V(x)$, we could have  deduced the $L^\infty$ bound by letting $p\to \infty$ in the estimates of the $L^p$ norm from the precised form \rife{globp}. Due to the term containing $V(x)$ (the linear perturbation term in \rife{ode}), the Moser-type iteration would be more involved  and we exploit a different argument, which might have  an interest in its own. Namely, we obtain $L^\infty$ bounds with the classical  {\it level set  iteration type argument} which dates back to E. De Giorgi and G. Stampacchia, see e.g.  \cite{St}.  

To this purpose, we use \rife{displa} with $U(r)= \frac{(r-k)_+^2}2$, where now $k>0$ is a positive parameter. 
%We notice that $U$ is a $C^1$ function and 
%$$
%\begin{cases}
%P(r)= U'(r)r- U(r) = (r-k)_+r- \frac{(r-k)_+^2}2\geq 0\,,
%& \\
%P'(r)= r\, \mathds1_{\{r>k\}} \,,  
%& \\
%P'(r)r- P(r)= kr \mathds1_{\{r>k\}}+\frac{(r-k)_+^2}2\geq 0\,. & 
%\end{cases}
%$$
%Hence, dropping positive terms, we obtain from 
As before, using \rife{coercf}  (we can suppose that $k\geq r_0$) and \rife{Hpp}, we obtain from \rife{displa}
\be\label{stampa0}\begin{split}
 \frac{d^2}{dt^2}\frac12 \into  (m(t)-k)_+^2   & \geq  \into mf'(m) H_{pp}Dm\cdot Dm \, \mathds1_{\{m>k\}}  + \into  m\mathds1_{\{m>k\}} H_{pp}Dm\cdot DV(x) 
 \\
 & \geq c_0\frac{\alpha_H}2 \into |D(m-k)_+|^2 - c_H \into m^2|DV|^2 \mathds1_{\{m>k\}}  \,.
\end{split}
\ee
%where we used as before \rife{fprime} (we can suppose that $k\geq r_0$) and assumption \rife{Hpp}.  
%Using as before Sobolev and Poincar\'e-Wirtinger inequality
%we have
%\begin{align*}
% \into |D(m-k)_+|^2  & \geq c_1 \left(\into  (m-k)_+^{2^*}\right)^{\frac2{2^*}} - c_2  |\Omega|^{-1} \left(\into (m-k)_+\right)^2
% \\ & \geq c_1  \left\{\left(\into  (m-k)_+^{2}\right)^{\frac{2(2^*-1)}{2^*}} \left(\into  (m-k)_+\right)^{-\frac{2(2^*-2)}{2^*}} \right\}- c_2  |\Omega|^{-1} \left(\into (m-k)_+\right)^2
%\end{align*} 
%where we also interpolated the $L^2$-norm between $L^1$ and $L^{2^*}$. Using the precise value of $2^*$, we insert this inequality in \rife{stampa0} and we get
%\begin{align*}
%\frac{d^2}{dt^2}\frac12 \into  (m(t)-k)_+^2 & \geq c_3 \left(\into  (m-k)_+\right)^{-\frac 4d  } \left(\into  (m-k)_+^{2}\right)^{1+ \frac 2d}
%\\ & \quad - c_2  |\Omega|^{-1} \left(\into (m-k)_+\right)^2- c_H \into m^2|DV|^2 \mathds1_{\{m>k\}} 
%\end{align*}
%which means
%$$
%-\vfi'' + 2c_3 a_k \vfi^{1+\frac2d} \leq b_k:= 2c_2  |\Omega|^{-1} \left(\into (m-k)_+\right)^2+ 2c_H \into m^2|DV|^2 \mathds1_{\{m>k\}} 
%$$
%for $\vfi= \into  (m(t)-k)_+^2$ and $a_k= \| m_k\|_1^{-\frac4d}$. We choose $k\geq \max (\|m_0\|_\infty, \|m_1\|_\infty)$, so that $\vfi$ vanishes at $t=0, T$. Hence by maximum principle we get
%$$
%\vfi(t):= \into  (m(t)-k)_+^2 \leq 
%$$

We choose $k\geq \max (\|m_0\|_\infty, \|m_1\|_\infty)$; then  we observe that
$$
\begin{cases}
\frac{d}{dt} \left(\into  (m(t)-k)_+^2 \right)  \mathop{\big|}_{ t=T} \leq 0\,,& \\
\noalign{\medskip}
 \frac{d}{dt} \left(\into  (m(t)-k)_+^2 \right)  \mathop{\big|}_{ t=0} \geq 0\,, & \end{cases}\qquad \forall k\geq \max (\|m_0\|_\infty, \|m_1\|_\infty)
$$
In particular,  if we integrate in time we get from \rife{stampa0}
\be\label{stampa1}
c_0\frac{\alpha_H}2 \int_0^T\into |D(m-k)_+|^2 \leq c_H \int_0^T\into m^2|DV|^2 \mathds1_{\{m>k\}}  \,.
\ee
Similarly, if we first integrate \rife{stampa0} in $(t,T)$, we have
$$
   \frac{d}{dt}  \into  (m(t)-k)_+^2 \leq c_H \int_t^T\into m^2|DV|^2 \mathds1_{\{m>k\}}  \leq   c_H \|DV\|_\infty^2\int_Q m^2\mathds1_{\{m>k\}} 
$$
which yields, since $k\geq \|m_0\|_\infty$,
\be\label{stampa2}
\qquad \forall t>0\,, \qquad   \into  (m(t)-k)_+^2 =    \int_0^t \frac{d}{ds}  \into  (m(s)-k)_+^2 ds \leq T\, c_H\,  \|DV\|_\infty^2\int_Q m^2\mathds1_{\{m>k\}} 
\ee
Using \rife{stampa1} and \rife{stampa2}, and the fact that $m^2$ is estimated in  any $L^p$-space (by Step 1), we  readily get the $L^\infty$-bound;  this is well-known (see e.g. \cite{LSU}) but we detail the  steps for the reader's convenience. We define the level set 
$$
A_k:= \{(t,x)\in Q\,: m(t,x)>k\}
$$
and we obtain from \rife{stampa1}--\rife{stampa2}
\be\label{stampa3}
\int_Q  |D(m-k)_+|^2  + \sup_{t\in (0,T)} \left(\into  (m(t)-k)_+^2\right)  \leq  C\, (1+ T) \|DV\|_\infty^2 \|m\|_{L^p(Q)}^2 |A_k|^{1-\frac2p}\,. 
\ee
We recall   the interpolation  inequality (see e.g. \cite[Proposition 3.1, Chapter 1]{DiB})
\be\label{dibe}
 \| v\|_{L^\rho(Q)}^\rho \leq  c\,  \| v\|_{L^\infty((0,T);L^2(\Omega))}^{\frac{4}d}     \|Dv\|_{L^2(Q)}^2 \qquad \hbox{where $\rho=2 \,(\frac{d+2}d)$}
\ee
for any $v\in L^2((0,T);W^{1,2}(\Omega))\cap L^\infty((0,T);L^2(\T^d))$ such that $\into v(t) \, dx=0$ a.e. in $(0,T)$. Applying this inequality to $v= [(m -k)_+- |\Omega|^{-1}\into (m(t) -k)_+\, dx]$ and using \rife{stampa3} we get
\begin{align*}
\int_0^T\into (m -k)_+^{\frac{2(d+2)}d}  & \leq  
\left[\sup_{t\in (0,T)} \left(\into  (m(t)-k)_+^2\right) \right]^{\frac 2d} \left(\int_Q  |D(m-k)_+|^2 \right)
+ c \int_0^T \left(\into (m-k)_+\,  \right)^{\frac{2(d+2)}d}
\\
& \leq 
c\, |A_k|^{(1-\frac2p)(1+\frac 2d)}  + c \int_0^T  |A_k(t)|^{1+\frac 4d} \into (m(t)-k)_+^{\frac{2(d+2)}d}  
\end{align*}
where $A_k(t)= \{x\in \Omega\,:\, m(t,x)>k\}$ is the time-section of $A_k$. Using \rife{globp},   we know that $\sup_{t} \, |A_k(t)|$ is small for sufficiently large $k$. Hence  the last term can be absorbed into the left-hand side and we deduce that there exists $k_0>0$ such that
$$
\int_0^T\into (m -k)_+^{\frac{2(d+2)}d}dxdt  \leq C\, |A_k|^{(1-\frac2p)(1+\frac 2d)}  \qquad \forall k\geq k_0>0\,.
$$
We deduce from this inequality that
$$
\forall h>k\geq k_0\, \qquad  |A_h| (h-k)^{\frac{2(d+2)}d} \leq C\, |A_k|^{\beta } \quad \beta:= (1-\frac2p)(1+\frac 2d)\,.
$$
Choosing $p$ sufficiently large, we have $\beta= (1-\frac2p)(1+\frac 2d)>1$; in that case, by a  classical iteration lemma (see \cite[Lemma 4.1]{St}, or similar  arguments in \cite[Chapter 2]{LSU}) we have $|A_{\bar k}| =0$ for some $\bar k$ only depending on $k_0$ and $C$. 
This means that 
$$
\|m\|_\infty \leq K
$$
for some $K$ depending on $\|m_0\|_\infty, \|m_1\|_\infty, \|DV\|_\infty, T, c_0, r_0, H$.

In a similar way we localize the previous estimate.
%, without using any information on $m_0, m_1$. 
To this goal, we fix $t_0\in (0,T)$,  and $R<R_0:= \min(t_0, T-t_0)$; then,  for $\tau\in (0,R)$  let $\xi(t)$ be a smooth cut-off function such that 
$$
\begin{cases} 
\xi (t)= 1 & \hbox{if $t\in (t_0-\tau, t_0+\tau)$} \\
\xi(t)=0 & \hbox{if $|t-t_0 | \geq R$}  \\
|\xi'(t)|^2 + |\xi''(t)|\leq \frac c{(R-\tau)^2}
\end{cases}
$$
We denote 
$$
A_{k, \tau}:= \{ (t,x)\in (t_0-\tau, t_0+\tau)\times \Omega \,:\, m(t,x)>k\}\,.
%\quad A_{k, R}:= \{ (t,x)\in (t_0-R, t_0+R)\times \Omega \,:\, m(t,x)>k\}\,.
$$
Then we have from \rife{stampa0}
\be\label{stampaloc}
\begin{split}
\frac{d^2}{dt^2} \left( \xi^2 \into  (m(t)-k)_+^2  \right)  &  \geq    c_0\frac{\alpha_H}2 \into \xi^2 \, |D(m-k)_+|^2 - c_H\|DV\|_\infty^2 \into m^2|\xi|^2 \mathds1_{\{m>k\}}  
\\
& \qquad - \frac c{(R-\tau)^2}\mathds1_{\{|t-t_0|<R\}} \into  (m(t)-k)_+^2  + 4 \xi \xi' \frac{d}{dt} \into  (m(t)-k)_+^2  \,.
\end{split}
\ee
Integrating   in $(t,T)$ we get
\begin{align*}
\frac{d}{dt}  \left( \xi^2 \into  (m(t)-k)_+^2\right)  & \leq c_1 \int_{A_{k,R}} m^2 + \frac {c_2}{(R-\tau)^2} \int_{A_{k,R}}   (m-k)_+^2
%\\
+ 4 \xi \xi'  \into  (m(t)-k)_+^2  
\end{align*}
and one more integration yields
\begin{align*}
  \xi^2 \into  (m(t)-k)_+^2 &  \leq c_1\, T \int_{A_{k,R}} m^2 + \frac {c_2 T}{(R-\tau)^2} \int_{A_{k,R}}  (m-k)_+^2 + \int_Q 4 \xi \xi'   (m-k)_+^2  
\\ & \leq c_1\, T\int_{A_{k,R}} m^2 + \frac {c_3 T}{(R-\tau)^2} \int_{A_{k,R}}  (m-k)_+^2
\\
& \leq \left( c_1T + \frac {c_3 T}{(R-\tau)^2}\right) \|m\|_{L^p(A_{k,R})}^2 |A_{k,R}|^{1-\frac2p}\,.
\end{align*}
Similarly, integrating \rife{stampaloc} in $(0,T)$, we estimate
\begin{align*}
  \int_Q \xi^2 \, |D(m-k)_+|^2   & 
 \leq c \int_{A_{k,R}} m^2+ \frac {c\, T}{(R-\tau)^2} \int_{A_{k,R}}    (m-k)_+^2
 \\
 & \leq \left( c_1T + \frac {c_3 T}{(R-\tau)^2}\right) \|m\|_{L^p(A_{k,R})}^2 |A_{k,R}|^{1-\frac2p} 
\end{align*}
for $T\geq 1$. Hence, using as before the Gagliardo-Nirenberg inequality we get
\be\label{tloct}
\begin{split}
\int_Q (\xi (m -k)_+)^{\frac{2(d+2)}d}  & \leq  
\left[\sup_{t\in (0,T)} \left(\into  \xi^2(m(t)-k)_+^2\right) \right]^{\frac 2d} \left(\int_Q  \xi^2|D(m-k)_+|^2 \right)
\\ & \qquad\qquad + c \int_0^T \left(\into \xi (m-k)_+\,  \right)^{\frac{2(d+2)}d}
\\
& \leq 
\left[\left( c_1T + \frac {c_3 T}{(R-\tau)^2}\right) \|m\|_{L^p(A_{k,R})}^2\right]^{1+\frac2d} |A_{k, R}|^{(1-\frac2p)(1+\frac 2d)} 
\\
& \qquad  + c \int_0^T  |A_k(t)|^{1+\frac 4d}  \into (\xi(m(t)-k)_+)^{\frac{2(d+2)}d}  \,.
\end{split}
\ee
Notice that $\|m(t)\|_{p}$ can be  estimated from \rife{locp}, for any $p>1$; indeed we have 
%, by a constant  depending on $p, (t_0-R), (T-t_0-R)$. Moreover
$$
\hbox{for $t\in (t_0-R, t_0+R)$,}\quad |A_k(t)| \leq k^{-p} \into m(t)^p \leq C\, k^{-p} \left( (t_0-R)^{-d(p-1)}+ (T-t_0-R)^{-d(p-1)}\right)
$$
hence, if $k$ is large, $|A_k(t)|$ is uniformly small for $t\in (t_0-R, t_0+R)$. Absorbing last term in \rife{tloct} we deduce that
$$
\int_Q (\xi (m -k)_+)^{\frac{2(d+2)}d}  \leq  C_R \left( c_1T + \frac {c_3 T}{(R-\tau)^2}\right)^{1+\frac2d} |A_{k, R}|^{(1-\frac2p)(1+\frac 2d)} 
$$
which implies
$$
|A_{h, \tau}| (h-k)^{\frac{2(d+2)}d}\leq    C_R \left(    \frac {c(R^2+1) T}{(R-\tau)^2}\right)^{1+\frac2d} |A_{k, R}|^{(1-\frac2p)(1+\frac 2d)} 
$$
for every $h>k\geq k_0$ ($k_0$ possibly depending on $R$), and for every $\tau\in (0,R)$. As before, we can choose $p$ large so that $\beta:=(1-\frac2p)(1+\frac 2d)>1$
and we conclude with (a localized version of) the iteration lemma (see \cite[Lemma 5.1]{St}) that, for any $\sigma\in (0,1)$,  $m(t)$ is bounded for $t\in (t_0- (1-\sigma)R, t_0+ (1-\sigma)R)$. In particular, we have
$$
\|m(t_0)\|_\infty \leq C( \max[t_0^{-1}, (T-t_0)^{-1}], \|DV\|_\infty, T)\,.
$$ 
\qed

\begin{remark} Let us point out that the estimates of Proposition \ref{supm} apply to the  solutions of the penalized problem \rife{MFG-delta} as well. Moreover, since the estimates only depend on  $\|V\|_{W^{1,\infty}(\Omega)}$, on the uniform convexity (and upper bound) of the Hamiltonian and on the quite general condition \rife{coercf}, those bounds are inherited by weak solutions of general systems \rife{mfpp} which can be obtained in the limit as $\vep\to 0$. See Definition \ref{def-weak} and Theorem \ref{weaksol} in the next Section.
\end{remark}

\section{Existence results}\label{exi}

In this Section we collect all previous estimates and ingredients to deduce our main existence results.

\subsection{Smooth solutions}

Here we prove the existence of smooth solutions for the elliptic problem \rife{quasell}; equivalently, this yields the existence of smooth solutions to the general mean-field planning problem \rife{mfpp}.

\begin{theorem}\label{exi-smooth}
Let $\Omega$ be a  $C^3$, bounded, convex domain  in $\R^d$ and let $V\in W^{2,\infty}(\Omega), m_0,m_1\in W^{1,\infty}(\Omega)$ such that $m_0,m_1>0$ in $\oo$. 
Let $H$ be a $C^3$ function satisfying   \rife{hrad} and  $f \in C^2(0,\infty)$ be  a nondecreasing function. 
Assume  that at least one of the two following conditions is satisfied:
\vskip0.5em
(i) $H$ satisfies conditions \rife{Hpp}-\rife{hppp}.
\vskip0.3em
(ii) $H$ satisfies conditions \rife{Hppq}-\rife{hpppq} for some $q>1$ and $\varpi>0$, $r\mapsto f'(r)r$ is nondecreasing, and $V(x)$ is convex. 
\vskip0.5em
Then  there exists    $u\in C^{2,\alpha}(Q)\cap C^{1,\alpha}(\oQ)$, $m\in C^{1,\alpha}(Q)\cap C^{0,\alpha}(\oQ)$ such that   $(u,m)$ is a smooth solution to the problem \rife{mfpp}, and in addition $m>0$ in $\bar Q$. We also have that $u$ is  a solution to the elliptic problem \rife{quasell} (unique solution up to addition of a constant) and $m$ is the (unique) minimizer of the problem\footnote{A priori the minimizer can be considered in the class of $m, v \in L^\infty(Q)$, or even in the broader class of absolutely continuous curves $m\in C^0([0,T];\cP(\Omega))$, interpreting $v\in L^2(dm(x,t))$ as the metric derivative of $m$ (see \cite{AGS}) and extending the functional to general measures in a classical way, see e.g. \cite{LaSa}. However, the a priori setting is a minor point here, because the minimizer $(m,H_p(Du))$ turns out to be smooth.}
\be\label{varia}
\min \,\,  \int_0^T\!\! \into  L(v) dm+    \int_0^T\!\! \into  m ( \vep [\log (m) -1]+ V) dxdt +  \int_0^T\!\! \into F(m)dxdt \,,\quad (m,v)\,:\,\, \begin{cases} m_t - \dive(vm)=0 & \\ 
H_p(Du)\cdot \vec\nu \mathop{|}_{\partial\Omega}=0 & 
\\ 
m(0)=m_0\,, m(T)=m_1 & \end{cases} 
\ee
where $F'(r)=f(r)$ and  $L$ is the Fenchel conjugate of $H$. 
%\be\label{quasi}
%\begin{system}
%& -{\rm tr}\left({\mathcal A}({\mathcal D}u) \, {\mathcal D}^2u \right) + DV(x) \cdot H_p(Du)=0 &&\text{in } Q, \\%(0,1)\times \R^d, \\
%& \,\, \, H_p(Du)\cdot \vec\nu=0   &&\text{in } (0,T)\times \partial \Omega, \\%(0,1)\times \R^d
%& -u_t + H(Du)=   f^\vep (m_1)  &&\hbox{at $t=T$, $x\in \Omega$,}
%\\
%&   -u_t + H(Du) = f^\vep (m_0)  &&\hbox{at $t=0$, $x\in \Omega$,}
%\end{system}
%\ee
%where $\mA$ is given by \rife{aij}. 
\end{theorem}

\proof In a  first step, we consider a suitable approximation $f_n$ of $f$ in a way that condition \rife{fprime} is satisfied. If we are under condition (i), we simply take 
$f_n(r)= f(T_n(r))$, where $T_n(\cdot)$ is a $C^2$ function such that $T_n(r)=r$ for  $0\leq r\leq n$, $T_n(r)\leq 2n$ for every $r$, and $|T_n'(r)|r+ |T_n''(r)|r^2\leq n$. Then $|f_n''(r)|r^2 \leq C_n $ and condition \rife{fprime} is   satisfied (for some $\beta$ depending on $n$). If we are in condition (ii), we wish to preserve the nondecreasing character of $f'(r)r$. To this purpose, given that $(f'(\tau)\tau)'\geq 0$, we may take 
\be\label{fn2}
f_n(r):= f(1)+ f'(1)\ln(r) + \int_1^r \frac1s \int_1^s [(f'(\tau)\tau)' \wedge n]d\tau\,.
\ee
We observe that 
\be\label{fnsec0}
f_n' (r) r= f'(1)+  \int_1^r [(f'(\tau)\tau)' \wedge n]d\tau 
\ee
and 
\be\label{fnsec}
f_n''(r)r^2 = - f'(1) - \int_1^r [(f'(\tau)\tau)' \wedge n]d\tau+ r [(f'(r)r)' \wedge n]\,.
\ee
There is no loss of generality in  assuming that $(f'(\tau)\tau)'\geq \de_n>0$   (otherwise replace $f$ with $f(r)+ \frac1n r$); then we have, for every $r>2$,
%hence  $\int_1^r [(f'(\tau)\tau)' \wedge n]d\tau\geq \de_n\, r $ if $r$ is sufficiently large. We conclude that, for some $s_0>0$, we have
$$
%\forall r>s_0\,, \quad
 r [(f'(r)r)' \wedge n] \leq n\, r \leq  \frac {2n}{\de_n}\int_1^r [(f'(\tau)\tau)' \wedge n]d\tau \leq   \frac {2n}{\de_n} f_n'(r) r \,.
$$
Hence \rife{fnsec} implies  
$$
|f_n''(r)|r^2 \leq  \left(1+ \frac {2n}{\de_n}\right)f_n'(r) r \leq  \left(1+ \frac {2n}{\de_n}\right) (1+ f_n'(r) r )^{\frac32}
$$
which means that $f_n$ satisfies condition \rife{fprime}. At the same time, we have here that $f_n'(r) r$ is nondecreasing.  Notice also that $f_n$ is nondecreasing (from \rife{fnsec0}) and satisfies
\be\label{fnf}
f(1)+ f'(1) \ln (r) \leq f_n(r)\leq f(r)\quad \forall r\in (0,\infty). 
\ee
Now we take the solution $(u_\de, m_\de)$ of problem \rife{MFG-delta} corresponding to $f=f_n$, which is guaranteed by Theorem \ref{delta}. We define $\hat u_\de$ as in  Lemma \ref{norm}. Now we show that there exists a constant $M>0$, independent of $n$ and $\de$, such that $\frac1M\leq m_\de\leq M$.  First of all, we observe that under either conditions (i) or (ii), we have that $m_\de$ is  bounded above (independently of $n$). Indeed, if (ii) holds true, this follows directly from Lemma \ref{sup}, Lemma \ref{mound} and the definition of $m_\de$; and since $f_n\leq f$, the bound is independent of $n$.  If rather (i) holds true, we first observe that 
$$
\| m_\de(0)\|_\infty \leq  (f_n^\vep)^{-1} ( \de\|u_\de\|_\infty  + f_n^\vep(\|m_0\|_\infty) )
$$
and then, due to  Lemma \ref{sup}, we have that $\| m_\de(0)\|_\infty$ is bounded. Similarly we have for $\| m_\de(T)\|_\infty$. Therefore, under condition (i), we deduce that  $m_\de$ is uniformly bounded  by Proposition \ref{supm}.  Notice again that, due to \rife{fnf}, the bound is independent of $n$.  

Thanks to the bound of $m_\de$ and since $\log(m_0), \log(m_1)$ are Lipschitz continuous, we deduce from Lemma \ref{norm}
that  $\hat u_\de$ is uniformly bounded  and solves
$$
\begin{system}
& -{\rm tr}\left({\mathcal A}({\mathcal D}\hat u_\de) \, {\mathcal D}^2\hat u_\de \right)+ DV(x) \cdot H_p(D\hat u_\de) =0 &&\text{in } Q, \\%(0,1)\times \R^d, \\
& -( \hat u_\de)_t + H(D\hat u_\de)= \de \hat u_\de + c_\de+ f_n^\vep(m_1) + V(x) &&\hbox{at $t=T$, $x\in \Omega$,}
\\
&   -(\hat u_\de)_t + H(D\hat u_\de)+ \de \hat u_\de + c_\de= f_n^\vep(m_0) + V(x) &&\hbox{at $t=0$, $x\in \Omega$,}
\end{system}
$$
where $c_\de= \de\int u_\de(T)m_1\, dx$ is a  bounded sequence of real numbers.  Now we apply  Theorem \ref{bound},  and we get that  $\hat u_\de$ is bounded in Lipschitz norm by a constant which is independent of $n$. In fact, if $M$ is such that $m_\de\leq M$, then condition \rife{fprime} is only needed for $r\leq M$; and clearly  $f_n$ satisfies 
\rife{fprime} for $1\leq r\leq M$, with some constant $\beta$ independent of $n$. Thus, the Lipschitz bound of $\hat u_\de$ is independent of $n,\de$;  by definition of $m_\de$, this means that there exists  a constant $M$, independent of $n$, such that
$$
\frac1M\leq m_\de\leq M\qquad \forall \de>0\,.
$$
Now, if $f_n$ is defined by \rife{fn2},   since $(f'(\tau)\tau)'$ is continuous for $\tau \in [\frac1M,M]$,   for $n$ sufficiently large we have $f_n(m_\de)= f(m_\de)$. The same obviously holds if $f_n= f(T_n(r))$. This means that, in both cases (i) or (ii),  $(\hat u_\de, m_\de)$ are actually solutions of \rife{MFG-delta} with $f$.

As a next step,  from Theorem \ref{lieb}, we have that 
%\cite[Lemma 2.4]{lieb84NA}, the Lipschitz bound on $\hat u_\de$ implies that the sequence 
$\hat u_\de$ is actually bounded in $C^{1,\alpha}(\bar Q)$ for some $\alpha>0$.
In particular, we have that  $m_\de= (f^\vep)^{-1}(-\partial_t \hat u_\de + H(D\hat u_\de)-V(x))$ is bounded in $C^{0,\alpha}(\bar Q)$ and, up to subsequences, $(\hat u_\de, m_\de)$ converge to some $(u,m)$ and this convergence is uniform up to $t=0,T$. In particular, we deduce from \rife{margi} and from the increasing character of $f^\vep$ that $m(0)=m_0$ and $m(T)=m_1$. This also implies  that $c_\de\to 0$. Finally, the ellipticity implies that $u$ is actually $C^{2,\alpha}(Q)$, $m\in C^{1,\alpha}(Q)$ and the   system \rife{mfpp} is satisfied in a classical sense. Equivalently,  $u$ satisfies the quasilinear elliptic problem  \rife{quasell} in a classical sense.  Notice that $m>0$ in $\oQ$ because of the gradient bound on $u$. 

Finally, $m$ is a  minimizer of the variational problem \rife{varia}; this is  standard whenever $m>0$ and $(u,m)$ is a classical solution of \rife{mfpp}. Indeed, system \rife{mfpp} represents  the optimality conditions for the state-adjoint state of the optimization problem. In addition, setting $w= mv$, it is well known (see \cite{BB}, \cite{BCS}) that   \rife{varia} can be rephrased as a  convex optimization problem in terms of the couple $(m,w)$; hence  $m$ is a minimizer. Due to the log term, here the convexity is strict and $m$ is unique. 
\qed

\vskip1em 
The case when the Hamiltonian may  be degenerate, or singular, at $p=0$, e.g. in the model case $H(p)=|p|^q$,   corrresponds to having $\varpi=0$ in assumption \rife{Hppq}. In that situation, while the main gradient estimates  remain true (see Theorem \ref{bound}), the problem \rife{quasell}  lacks of uniform ellipticity and  $u$ cannot be proved to be smooth. The natural framework in that context is to use the divergence structure of the operator, namely of the continuity equation for the density. This yields the following result.

\begin{theorem}\label{exi-div}
Let $\Omega$ be a  $C^3$, bounded, convex domain  in $\R^d$ and let $ m_0,m_1\in W^{1,\infty}(\Omega)$ be such that $m_0,m_1>0$ in $\oo$. 
Assume that  $H\in C^3(\R^d\setminus \{0\})$   satisfies   \rife{hrad} and  \rife{Hppq}-\rife{hpppq}, and that $f\in C^2(0,\infty)$ is such that  $r\mapsto f'(r)r$ is nondecreasing  in $[0,\infty)$. Let $V\in     W^{2,\infty}(\Omega)$ be a convex function.

Then  there exists     $u\in W^{1,\infty}( Q) $, $m\in L^\infty( Q) $ such that   $(u,m)$ is a  solution to  \rife{mfpp}, where the continuity equation is taken in weak sense.  Moreover,  $u$ solves problem \rife{quasell} in weak sense, i.e.
\be\label{divu}
\into m_1\, \xi(T)\,dx - \into m_0\, \xi(0)\,dx- \int_0^T\into \vfi(-u_t + H( Du) - V(x)) \left[ \xi_t - H_p(Du)\cdot D\xi\right]\,dxdt =0 \quad \forall \xi\in C^1(\oQ)\,, 
\ee
where $\vfi(\cdot)= (f + \vep \log)^{-1}(\cdot)$.  We also have that $m>0$ in $\oQ$ and is the unique minimizer of \rife{varia},  and $u$ is the unique (up to addition of a constant) weak solution of \rife{divu}.
\end{theorem}

\proof We  approximate $H= h(|p|)$ with $H^n= h(\sqrt{\frac1n + |p|^2})$,  
%(possibly  regularizing $h$ with some $h_n$ satisfying \rife{Hppq}--\rife{hpppq} uniformly in $n$), 
so we can build a  solution $(u^n, m^n)$ according to Theorem \ref{exi-smooth}. Our assumptions allow us to apply Lemma \ref{mound}, so $m^n$ is bounded in $Q$.  From Lemma \ref{norm} we deduce that $\hat u^n$ (the normalization of $u^n$ as  in \rife{normou}) is uniformly bounded. Then we use Theorem \ref{bound} (since $m^n$ is bounded,   the assumption \rife{fprime} is only needed in a compact set  $[s_0, M]$, where it holds  because $f$ is $C^2$) and we get that $\hat u^n$   is bounded in Lipschitz norm, and then it is relatively compact in $C(\oQ)$.  We denote by $u\in W^{1,\infty}(Q)$ the uniform limit of (a subsequence of) $\hat u^n$.  We observe that $\hat u^n$ is a sequence of solutions to divergence form problems
$$
 -\partial_t \left( A_0(x,\mD \hat u^n ) \right) - \sum_{k=1}^d \frac{\partial}{\partial x_k}\left(A_k(x, \mD  \hat u^n )\right) =0
$$ 
where $A_0= -m^n= -(f^\vep)^{-1}(-  \hat u_t^n + H^n( D \hat u^n) - V(x))$, $A_k= (f^\vep)^{-1}(- \hat u_t^n + H^n( D \hat u^n) - V(x))H_{p_k}^n(D \hat u^n)$. We notice that the uniform gradient bound on $u^n$  implies that $m^n$ is uniformly bounded from below and from above.

Now we show the convergence of $\mD \hat u^n$ in $L^p(Q)$ for every $p<\infty$. To this purpose, we multiply the above equation by $( \hat u^n-u)$ and we get, using the divergence form structure,
\be\label{diven}
  \int_0^T \into \sum_{i=0}^d A_i(x,\mD  \hat u^n ) \mD_i ( \hat u^n-u)  \, dxdt = \into m_0 ( \hat u^n-u)(0)\, dx - \into m_1 ( \hat u^n-u)(T)\, dx\,.
\ee
The right-hand side goes to zero as $n\to \infty$ due to the uniform convergence of $ \hat u^n$. In the left-hand side we use the definition of $m^n$ and we get
\be\label{premn}
 - \int_0^T \into m^n ( \hat u^n-u)_t \, dxdt + \int_0^T \into m^n H_p^n(D \hat u^n) D( \hat u^n-u)  \, dxdt  \mathop{\rightarrow}^{n\to \infty} \,\, 0\,.
\ee
Recalling that $- \hat u^n_t = f^\vep (m^n)+ V- H^n(D \hat u^n)$, we obtain
\begin{align*}
  & \int_0^T \into m^n \, f^\vep(m^n) \, dxdt =    \int_0^T \into  m^n (-  u_t + H(Du)- V)\, dxdt   \\
  & \qquad \qquad +     \int_0^T \into m^n \left\{H(D \hat u^n)- H(Du)-  H_p(D \hat u^n) D( \hat u^n-u) \right\} \, dxdt  + o(1)_n 
  \end{align*}
where $o(1)_n$ denotes some quantity which vanishes as $n\to \infty$. Now we observe that the last integral is nonpositive, due to the convexity of $H$. If we denote by $m$ the weak$^*$ limit of $m^n$ in $L^\infty(Q)$, we get
\be\label{fbare}
\limsup_{n\to \infty} \int_0^T \into m^n \, f^\vep(m^n) \, dxdt\leq  \int_0^T \into  m (-u_t + H(Du)- V)\, dxdt\,.
\ee
But using again the convexity of $H$, the weak convergence of $ D  \hat u^n$ and denoting by $\bar f$ the weak$^*$ limit of $f^\vep (m^n)$ in $L^\infty(Q)$, we have
$$
- \hat u^n_t + H^n(D \hat u^n) = f^\vep (m^n)+ V  \quad \Rightarrow \quad -u_t + H(Du)  \leq \bar f + V\,.
$$
Therefore \rife{fbare} implies
$$
\limsup_{n\to \infty} \int_0^T \into m^n \, f^\vep(m^n) \, dxdt\leq  \int_0^T \into  m \, \bar f\, dxdt
$$
which yields
$$
\limsup_{n\to \infty} \int_0^T \into (m^n-m) \, (f^\vep(m^n)-f^\vep (m)) \, dxdt\leq 0\,.
$$
The strict monotonicity of $f^\vep$ implies that $m^n\to m$ in $L^1(Q)$, and then in $L^p(Q)$ for every $p<\infty$. Hence we deduce from \rife{premn} 
\begin{align*}
\int_0^T \into m^n (H_p(D \hat u^n)-H_p(Du)) D( \hat u^n-u)  \, dxdt & = -\int_0^T \into m^n  H_p(Du) D( \hat u^n-u)  \, dxdt
\\
& + \int_0^T \into m^n ( \hat u^n-u)_t \, dxdt + o(1)_n\,.
\end{align*}
All terms in the right-hand side converge to zero. Then,  using also the bound from below on $m^n$, we obtain
$$
\int_0^T \into   (H_p(D \hat u^n)-H_p(Du)) D( \hat u^n-u)  \, dxdt \mathop{\rightarrow}^{n\to \infty} \,\, 0\,.
$$
The strict monotonicity of $H_p(\cdot)$ implies that $D \hat u_n\to Du$ in $L^p(Q)$ for every $p<\infty$. The identity $- \hat u^n_t = f^\vep (m^n)+ V- H^n(D \hat u^n)$ allows us to conclude 
that $ \hat u^n_t\to u_t$ strongly in $L^p(Q)$ as well. Thus, we proved that $\mD \hat u^n\to \mD u$ in $L^p(Q)$ for every $p<\infty$. Passing to the limit in the equation, we conclude that $u$ is a weak solution of the limit problem in the sense of \rife{divu}. Accordingly, $(u,m)$ satisfy \rife{mfpp}, where $u$ satisfies the first equation almost everywhere in $Q$, and $m$ is a bounded weak solution of the continuity equation. We notice that $u$ is the unique Lipschitz solution (up to addition of a constant) of \rife{divu}; indeed, \rife{divu} can be rephrased as 
$$
 \int_0^T \into \sum_{i=0}^d A_i(x,\mD   u ) \mD_i \xi  \, dxdt =  \into m_0 \xi(0)\, dx - \into m_1 \, \xi(T) \, dx\quad \forall \xi\in C^1(\oQ)
$$
where ${\bf A}= (A_i)_{i=0,\dots d}$ is a strictly monotone vector field in $\R^{d+1}$.  Moreover, by a standard density argument, the weak formulation holds for all test functions  $\xi \in W^{1,\infty}(Q)$. Hence, if $u_1,u_2$ are two different solutions, using $\xi=u_1-u_2$ in the weak formulation of both equations, we readily conclude that $\mD u_1= \mD u_2$. So $u_1-u_2$ differ by a constant. Finally, the uniqueness of $m$ as a minimizer can be proved as in Theorem \ref{exi-smooth} by using the strict convexity of the associated functional.
%The only difference\footnote{rivedere} now lies inside Proposition \ref{boundaryHol} in the Appendix; indeed, after localizing near a portion of $\partial\Omega$ and after reflection, one has to use the divergence structure form for the equation of $u^*$ (as we did in \rife{cono} with $\psi$ in place of $u$). In fact, one can rewrite problem \rife{reflected} in divergence form, where a co-normal boundary condition is embedded into integration by parts, as is usual. Then one applies \cite[Thm 2]{lieb88} to deduce the $C^{1,\alpha}$ estimates up to the boundary. Once $\hat u_\vep$ is estimated in $C^{1,\alpha}(\oQ) $, then it is relatively compact in $C^1(\oQ)$, 
%and  a solution $u$ is found 
%satisfying  \rife{divu}. Accordingly, $(u,m)$ satisfy \rife{mfpp}.
\qed

\subsection{Transport with entropy in $\R^d$}\label{OT}

In this section we restrict to the   case of optimal transportation without any further congestion term (or mean field interaction). This means that we set $f=0$ and we focus on the very model problem
\be\label{OT}
\begin{cases}
-u_t + H( Du)  =   \vep \log m + V(x)  & \text{in }  Q:= (0,T)\times \Omega, \\
 m_t  - \dive (m\,H_p(Du))  = 0   &\text{in } Q,
\\
 \,\,  H_p(Du)\cdot \vec\nu=0  &\text{on $(0,T)\times \partial \Omega$,}
\\
 m(0,\cdot) = m_0, \; m(T,\cdot)  = m_1  &\text{in } \Omega\,,
\end{cases}
\ee
For the reader's convenience, we first  collect the main result of the previous section in the specific case of \rife{OT}. 

\begin{theorem}\label{thm1OT} Let $\Omega$ be a  $C^3$, bounded, convex domain  in $\R^d$ and let $V\in W^{2,\infty}(\Omega), m_0,m_1\in W^{1,\infty}(\Omega)$ such that $m_0,m_1>0$ in $\oo$. 
Let $H$ be a $C^3$ function satisfying   \rife{hrad} and  assume  that at least one of the two following conditions is satisfied:
\vskip0.5em
(i) $H$ satisfies conditions \rife{Hpp}-\rife{hppp}
\vskip0.3em
(ii) $H$ satisfies conditions \rife{Hppq}-\rife{hpppq} for some $q>1$ and $\varpi>0$,  and $V(x)$ is convex. 
\vskip0.5em
Then problem \rife{OT}  admits  a solution $(u,m)$ such that $u\in C^{2,\alpha}(Q)\cap C^{1,\alpha}(\oQ)$, $m\in C^{1,\alpha}(Q)\cap C^{0,\alpha}(\oQ)$, and this solution is unique (up to addition of a constant to $u$). Moreover we have  $m(t) >0 $ in $\oo$ for every $t\in [0,T]$. Finally, $m$ is the unique minimizer of the problem
$$
\min\,\,   \int_0^T\!\!\into   L(v)\,dm + \vep \int_0^T\!\!\into  \log\left( \frac {dm}{d\varrho}\right) dm\,,\quad (m,v): \, \begin{cases}
m_t - \dive(mv)=0\,, & \hbox{in $(0,T)\times \Omega$,}\\
v\cdot \vec\nu=0 & \hbox{on $(0,T)\times \partial \Omega$,} \\
 m(0)=m_0, \, m(T)=m_1 & \end{cases} 
$$ 
where $\varrho= e^{- \frac{V(x)}\vep}dx$ and $L$ is the Fenchel conjugate of $H$.
\end{theorem}

\vskip1em

We now devote the remaining part of this  subsection to a model result in the noncompact setting where $\Omega$ is replaced by 
the whole space $\R^d$. In this context, it is natural to consider  the entropy of $m$ with respect to Gaussian-type measures. 

Let us fix a reference measure $\nu:= e^{-V} dx$, where $V$ is twice differentiable and satisfies
%\footnote{penso che in setting Riemanniano dovrebbe leggersi come: $D^2 V + Ric \geq \gamma I$, che se non erro sarebbe  la condizione di Bakry-Emery?}
\be\label{dissip}
D^2 V \xi\cdot \xi \geq \gamma\, |\xi|^2 \qquad \gamma>0\,.
\ee
In particular, this implies that $V$ is a convex, coercive function at infinity, whose model case is given by $V= \gamma \frac{|x|^2}2$. 
 We consider the system
\be\label{OTrd}
\begin{system}
&-u_t + H(Du)  =\vep ( \log m + V(x) )&&\text{in }  (0,T)\times \R^d, \\
&m_t  - \dive (m\, H_p(Du) )  = 0  &&\text{in } (0,T)\times \R^d,\\
&m(0,\cdot) = m_0, \; m(T,\cdot)  = m_1 &&\text{in } \R^d\,,
\end{system}
\ee
where we assume that $m_0, m_1 $ satisfy
\be\label{mu01}
m_0 e^V, m_1 e^V \in W^{1,\infty}(\R^d)\qquad \hbox{and}\quad c_0 \, e^{-V(x)}\leq m_0\leq C_0 \, e^{-V(x)}\,,\qquad c_1 \, e^{-V(x)}\leq m_1\leq C_1 \, e^{-V(x)}
\ee
for some positive constants $c_i, C_i>0$, $i=0,1$. 
%This system corresponds to the minimization of the standard quadratic functional with additional entropy computed with respect to the measure $\nu$:
%$$
%\inf \{ \int_0^T\int_{\R^d}   m\, L(v) dxdt + \vep \int_0^T\int_{\R^d}  \log\left( \frac {dm}{d\nu}\right) dm\}\,,\qquad m_t = \dive(mv)=0\,,\,\,\, m(0)=m_0, \, m(1)=m_1
%$$
%Notice that, setting $ \mu:= m\, e^{V(x)}$, the above system can be rephrased as
%$$
%\begin{system}
%&-u_t + H(Du)   = \log \mu &&\text{in }  (0,1)\times \R^d, \\
%&\partial_t \mu  - \dive (\mu\, H_p(Du)) + \mu \, DV\cdot H_p(Du)  = 0  &&\text{in } (0,1)\times \R^d,
%%\\&m(0,\cdot) = m_0, \; m(1,\cdot)  = m_1 &&\text{in } \R^d\,,
%\end{system}
%$$
We can prove a similar result as above for the non compact case; for simplicity we restrict to  Hamiltonians with quadratic growth.

\begin{theorem}\label{OTspace} Assume that  $V\in W^{2,\infty}(\R^d)$ satisfies \rife{dissip}, and that $ m_0,m_1 $ satisfy \rife{mu01}.  
Let $H$ be a $C^3$ function satisfying   \rife{hrad}  and \rife{Hpp}-\rife{hppp}.

Then problem \rife{OTrd}  admits  a solution $(u,m)$ such that $u\in C^2(Q), m\in C^1(Q)$, and this solution is unique (up to addition of a constant to $u$). Moreover we have  that $m(t) e^{V(x)} $ is positive and bounded uniformly in $\R^d$ for every $t\in [0,T]$. Finally, $m$ is the unique minimizer of the problem
$$
\min \,  \int_0^T\into   m\, L(v)\, dxdt + \vep \int_0^T\into  \log\left( \frac {dm}{d\nu}\right) dm\,,  \quad (m,v)\,:\,\, \begin{cases} m_t - \dive(vm)=0 & \\ m(0)=m_0\,, m(T)=m_1 & \end{cases} 
$$ 
where $\nu= e^{- V(x) }dx$.
\end{theorem}

\proof  We define $B_n$ the $d$-dimensional ball of radius $n$, and we let 
%$$
%\Omega_n:= \{x\in \R^d\,: \, V(x) <n \}
%$$
%which is a $C^2$ open set, and we observe that $\Omega_n$ is bounded and strictly convex because of assumption \rife{dissip}. We let 
$(u_n,m_n)$
be the solution of the problem
$$
\begin{cases}
-u_t + H( Du)  =   \vep (\log m + V(x))  & \text{in }   (0,T)\times B_n, \\
 m_t  - \dive (m\,H_p(Du))  = 0   &\text{in } (0,T)\times B_n,
\\
 \,\,  H_p(Du)\cdot \vec\nu=0  &\text{on $(0,T)\times \partial B_n$,}
\\
 m(0,\cdot) = m_0, \; m(T,\cdot)  = m_1  &\text{in } B_n\,,
\end{cases}
$$
which exists by Theorem \ref{thm1OT}, with $u_n, m_n$ smooth since $m_0, m_1$ are Lipschitz and positive in $B_n$. 
Since $V$ is convex, we can   estimate  $-u_t + H(Du)$ in terms of its maximum at $t=0,T$ (see   Lemma \ref{mound}).  We obtain 
that
\be\label{munb}
 m_n e^V  \leq  \sup_{i=0,1} \, \|m_i e^V\|_\infty   
\ee
hence $m_n$ is bounded and uniformly decaying at infinity. 

Now we estimate $Du$ with a variant of Theorem \ref{bound}; in fact, we estimate $Du$ {\it independently} from the $L^\infty$-bound of $u$. Hereafter, we have $u=u_n$ (i.e. we avoid to write the index $n$).
We consider $w= H(Du)$ and we recall (see \rife{hdueq}) that $w$ satisfies
$$
-{\rm tr}\left({\mathcal A} \, {\mathcal D}^2 w\right)  + H_{pp}DV \cdot Dw + D^2V(x) H_p(Du)H_p( Du)   
  =       {\rm tr}\left([ \mD_\eta {\mathcal A} \cdot \mD w ] D^2 u \right)
-   \sum\limits_{k,\ell=1}^d  H_{p_kp_\ell}\left( \mA  \,\mD u_k \cdot \mD u_\ell\right)
$$
where, we recall,  $H_{pp}$ is computed on $Du$. This  implies, due to \rife{dissip},
$$
-{\rm tr}\left({\mathcal A} \, {\mathcal D}^2 w\right)  + H_{pp}DV \cdot Dw   
 -  {\rm tr}\left([ \mD_\eta {\mathcal A} \cdot \mD w ] D^2 u \right)  + \gamma \, |H_p(Du)|^2  \leq 0 \,.
$$
%
%
%and we can proceed as in the previous Section, i.e. reducing the system to a single elliptic equation for $u$:
%$$
%-u_{tt} + 2 Du_t \cdot Du - D^2u Du Du -\Delta u + DV \cdot Du =0\,.  
%$$
%Same additional terms can be added as in problem \rife{eps-de}. I omit them here in the next computation.  Let us observe that the previous equation can be re-written as
%$$
%u_{tt} + \Delta u =  2 \partial_t \left(\frac {|Du|^2}2\right)- D \left(\frac {|Du|^2}2\right)\cdot D u + DV \cdot Du
%$$ 
%If we set $w:=  \frac {|Du|^2}2 $, we compute\footnote{this is Lemma \ref{struc} but this is also, more simply, the Bochner formula in a Riemannian setting}
%\begin{align*}
%w_{tt} + \Delta w  & =  D \left(u_{tt} + \Delta u \right)\cdot Du + |D^2u|^2 + | Du_t|^2 \\
%& = 2 Dw_t \cdot Du - D^2 w Du \cdot Du- |Dw|^2  + DV \cdot Dw + D^2 V Du \cdot Du + |D^2u|^2 + | Du_t|^2 
%\end{align*}
%which implies
%$$
%-{\rm tr}\left({\mathcal A}(Du) \, {\mathcal D}^2w \right) - |Dw|^2 + DV \cdot Dw + D^2 V Du \cdot Du \leq 0\,.
%$$
Now we multiply $w$ by  an auxiliary function to take care of boundary conditions. We set
$$
z:= e^{\sigma\left(t-\frac{t^2}{T}\right)} \, w\,.
$$
Computing we find that $z$ solves 
\begin{align*}
-{\rm tr}\left({\mathcal A}(Du) \, {\mathcal D}^2z \right)  &      - e^{\sigma\left(t-\frac{t^2}{T}\right)} \sigma \, w (\frac 2T+ (1-2\frac tT)^2 \sigma) +  2 \sigma\, (1-2\frac tT) (z_t - H_p(Du) \cdot Dz) 
\\
 & \quad   + H_{pp}DV \cdot Dz   
 -  {\rm tr}\left([ \mD_\eta {\mathcal A} \cdot \mD z ] D^2 u \right)+ e^{\sigma\left(t-\frac{t^2}{T}\right)} \, \gamma \, |H_p(Du)|^2  \leq 0 
\end{align*}
where we used that $\mathcal A$ is independent of $u_t$, so $ e^{\sigma\left(t-\frac{t^2}{T}\right)}\mD_\eta {\mathcal A} \cdot \mD w= \mD_\eta {\mathcal A} \cdot \mD z$.  
Using \rife{coercq} (with $q=2$) we have $|H_p|^2 \geq \gamma_0\, H- \gamma_1$ for some constants $\gamma_0, \gamma_1$. Hence, we deduce 
$$
-{\rm tr}\left({\mathcal A}(Du) \, {\mathcal D}^2z \right) + z \left[ \gamma\, \gamma_0 - \sigma\, (\frac2T+ (1-2\frac tT)^2 \sigma)\right] + {\mathcal B}\cdot {\mathcal D} z\leq \gamma\, \gamma_1 \,e^{\sigma\left(t-\frac{t^2}{T}\right)}
$$
for some vector field $\mathcal B$. In particular, for $\sigma $ sufficiently small (only depending on $T, \gamma_0, \gamma$), we deduce that $z$ is uniformly bounded on any internal maximum point.   We can also exclude that $z$ has maximum on $\partial B_n$, due to the Neumann condition and the  convexity of $B_n$. On the time-boundary,  we have, using that $-u_t + H(Du)= \vep \log (m_ie^V)$,
\be\label{neu}
z_t= e^{\sigma\left(t-\frac{t^2}{T}\right)} \left( \sigma (1-2\frac tT) w+ H_p \cdot Du_t\right) = e^{\sigma\left(t-\frac{t^2}{T}\right)} \left( \sigma (1-2\frac tT) w- \vep H_p \cdot D \log(m_i e^V)  \right) + H_p \cdot Dz 
\ee
for $i=0,1$. 
%, where $\mu_0= m_0 e^{V(x)}$, $\mu_1= m_1 e^{V(x)}$. If we assume that 
%$$
%m_0 e^{V(x)}, m_1 e^{V(x)} \in W^{1,\infty}(\R^d)
%$$
By assumption \rife{mu01}, we have that $\log(m_i e^V) \in W^{1,\infty}(\R^d)$. Then it
  follows from \rife{neu} that $z$ cannot have a maximum at $t=0$ or $t=T$ such that
$|Du| > K $, for some $K$ depending on $\sigma, H, \vep  \max_i \, \|D\log(m_i e^V)\|_\infty$.  
We deduce from the above steps the a priori estimate
$$
\| Du\|_\infty \leq C\, (1+ \vep  \max \left(\|D\log(m_0 e^V)\|_\infty , \|D\log(m_1 e^V)\|_\infty\right) )\,,
$$
for some $C$ only depending on $\gamma, H, T$.  This estimate, together with assumption \rife{mu01},  imply that $|u_t|$ is bounded at $t=0, t=T$. Then, using Lemma \ref{ut} we also deduce a uniform bound for $u_t$.  Therefore we proved so far that
$$
\| \mathcal D u_n\|_\infty \leq C\,,\qquad  \alpha_0 \leq m_n e^{V(x)} \leq \alpha_1  \qquad \forall x\in B_n
$$
for some $\alpha_0, \alpha_1>0$.
At this stage  we can normalize $u$ as in Lemma \ref{norm} by setting 
$$
\hat u_n:= u_n - \int_{B_n} u_n(T) m_1\, dx\,.
$$
We stress that the proof of Lemma \ref{norm} can be adapted here to give a   bound 
on  $\hat u_n$ which is uniform in $B_n$. To this purpose one first needs to use \rife{munb} to obtain the estimate \rife{at0}. Secondly, when using   the Wasserstein geodesic in the proof of Lemma \ref{norm}, it is enough to notice that the two marginals $m_0,m_1$ have finite moments (of any order), and one can still obtain estimate \rife{above} using the geodesic between $m_1$ and measures $\tilde m$ with finite  second moments. Thus, following Lemma \ref{norm}, we obtain that $\hat u_n$ is bounded in $L^\infty(B_n)$.

We are only left with letting $n\to \infty$. To this goal, we can use the local $C^{1,\alpha} $ estimates (up to $t=0,t=T$, as in Proposition \ref{boundaryHol}) to get compactness of $\cD u_n$ and $m_n$. With a diagonal process on a sequence of balls invading $\R^d$, we build a solution $u$ (and $m$ in turn) of the problem in the whole space.
\qed

\subsection{Weak solutions, extensions and further comments}

We conclude with a  few further remarks.
%\section{Applications to mean field game systems}\label{MFG}

\vskip1em
\paragraph{A further estimate for optimal transport.}  

It is natural to ask what happens if we drop the strict positivity condition on the marginals $m_0, m_1$, together with letting the entropy term vanish, or even independently.
Alternatively, one can wonder whether the positivity estimates for $m$, or the gradient estimates for $u$, can be localized in time.  This seems to be an interesting question.  In this direction, we  only give the following estimate for the case of pure transport with entropy. It suggests, roughly speaking, that $u$ should remain smooth at least on the support of $m$, giving a precise quantitative estimate. Let us stress that this result only requires boundedness of the marginals $m_0, m_1$, and the estimate is independent of $\vep$. 

\begin{theorem}\label{thm2OT}  Under the assumptions of Theorem \ref{thm1OT}, let $(u,m)$ be the (smooth) solution of \rife{OT}. 
Then  there exists a constant $\theta>0$ such that 
\be\label{mduloc}
\theta H(Du(t)) + \vep \log (m(t)) \leq K_t \qquad \forall t\in (0,T)
%\exp\left(\frac\theta\vep H(Du(t))\right) m(t)\, e^{V(x) /\vep}\leq K_t \qquad \forall t\in (0,T)
\ee
where $K_t$ is a constant only depending on   $\|u\|_\infty, \|V\|_\infty, T$ and $\min(t, T-t)^{-1}$.
\end{theorem}

\proof By assumption \rife{Hppq}, there exists $\theta>0$ such that
\be\label{2theta}
H_p(p)\cdot p \geq  (1+ 2\theta) H(p) - c_0 \qquad \forall p\in \R^d \,.
\ee
We come back to the gradient estimate of Theorem \ref{bound}, which we aim at localizing in time. 
To this purpose  we set  
$$
z:= (1+\theta) H(Du) - u_t + \lambda \frac{u^2}2\,,\qquad \lambda= \frac \sigma{1+ \|u\|_\infty^2} 
$$
where $\sigma$ will be chosen  (sufficiently small) as we did in Theorem \ref{bound}.
By using Lemma \ref{struc}, we obtain (recall that $H_p, H_{pp}$ are computed on $Du$)
 \begin{align*}
& -{\rm tr}\left({\mathcal A} \, {\mathcal D}^2 z \right)+ H_{pp}DV \cdot Dz+ \rho z 
 +  \lambda u\,  \left(DV\cdot H_p - H_{pp}DV \cdot Du\right)  +  (1+\theta)  D^2V(x) H_p\cdot H_p 
\\
&    +   \lam\,\mA  \mD u\cdot \mD u - {\rm tr}\left(  [\mD_{\eta} {\mathcal A} \cdot \mD z] \cD^2 u \right)
 =   -  \lam \sum\limits_{i,j,\ell=1}^n u \partial_{\eta_\ell} (a_{ij}) \partial_\ell u \,u_{ij}
 -  (1+\theta)\sum\limits_{k,\ell=1}^d  H_{p_kp_\ell}\left( \mA  \,\mD u_k \cdot \mD u_\ell\right)
 \end{align*}
where we took advantage that, here, the matrix $\mathcal A$ does not depend on $(t,x)$ (and not either on $u_t$, which will be used below).
Using the coercivity of $\mathcal A$ and $H_{pp}$, as we did in the proof of Theorem \ref{bound}, we get that for $|Du|$ sufficiently large
\be\label{thetadu1} \begin{split}
& -{\rm tr}\left({\mathcal A} \, {\mathcal D}^2 z \right)+ H_{pp}DV \cdot Dz+ \rho z - {\rm tr}\left(  [\mD_{\eta} {\mathcal A} \cdot \mD z] \cD^2 u \right)
\\ & \qquad 
 +  \lambda u\,  \left(DV\cdot H_p - H_{pp}DV \cdot Du\right)  +  (1+\theta)  D^2V(x) H_p\cdot H_p  +   \lam\, \left( |H_p(Du)\cdot Du - u_t|^2 + \vep \,\gamma_H |Du|^q\right) 
\\
&     + (1+\theta ) \gamma_H|Du|^{q-2} \left(|Du_t-   D^2 uH_p(Du)|^2 + \vep \gamma_H\, |Du|^{q-2} |D^2 u|^2\right)
\leq 
%  {\rm tr}\left(  [\mD_{\eta_\ell} {\mathcal A} \cdot \mD z] \right)
     -  \lam \, u\sum\limits_{i,j,\ell=1}^n  \partial_{\eta_\ell} (a_{ij}) \partial_\ell u \,u_{ij}\,.
 \end{split}
 \ee
We estimate last term in a similar way as in Theorem \ref{bound}
\begin{align*}
    \lam\, u \sum\limits_{i,j,\ell}   \partial_{\eta_\ell} (a_{ij}) \partial_\ell u \, u_{ij} &  
    =  \lam u \left\{  2  \, H_{pp}Du[D^2u\, H_p(Du)-Du_t]  +  \vep \sum\limits_{i,j,k=1}^d   H_{p_ip_jp_k} \partial_k u\, u_{ij}\right\}
    \\ & \leq  \lam \, (\beta_H+ \vep) \,  | D u|^q  
%+ \frac\lam 2 |- v_t + H_p(Du)\cdot Dv|^2 
  + \lam\, u^2\, \beta_H\, |Du|^{q-2}\,|D^2u\, H_p(Du)-Du_t|^2 +    C_H\,  \lambda u^2\,   \vep\, |Du|^{2(q-2)}\,   |D^2u|^2    
  %  \\
  \\  & \leq  \lam \, (\beta_H+ \vep) \,  | D u|^q    
%+ \frac\lam 2 |- v_t + H_p(Du)\cdot Dv|^2 
  + \sigma \beta_H \, |Du|^{q-2}\, |D^2u\, H_p(Du)-Du_t|^2 +    \sigma\,  C_H\, \vep\, |Du|^{2(q-2)} \,  |D^2u|^2    
\end{align*}
where we used that $\lam u^2\leq \sigma$ due to the choice of $\lambda$. Notice that, for a sufficiently small $\sigma$ (independent of $\vep$!), last two terms can be absorbed by the left-hand-side in \rife{thetadu1}. 
% \begin{align*}
%& -{\rm tr}\left({\mathcal A} \, {\mathcal D}^2 z \right)+ H_{pp}DV \cdot Dz+ \rho z 
% + \lambda u\,  \left(DV\cdot H_p - H_{pp}DV \cdot Du\right)  +  (1+\theta)  D^2V(x) H_p(Du)H_p( Du) 
%\\
%&    +   \lam\, \left( |H_p\cdot Du - u_t|^2 + \vep \,\frac{\alpha_H}2 |Du|^2\right)  + \frac{(1+\theta )}2 \alpha_H \left(|Du_t- H_p D^2 u|^2 + \vep \alpha_H |D^2 u|^2\right)
%\\  & \qquad \leq 
%  {\rm tr}\left(  [\mD_{\eta_\ell} {\mathcal A} \cdot \mD z] \right)
% \end{align*}
Similarly we estimate, still using   $\lam\, u^2\leq \sigma$, 
\begin{align*}
 & \lambda u\,  \left(DV\cdot H_p - H_{pp}DV \cdot Du\right)  +  (1+\theta)  D^2V(x) H_p\cdot H_p 
\\ & \qquad \geq  - C (\lam+1)    |Du|^{2(q-1)} -  \sigma\, C_{V,H}     \,.
\end{align*}
Once we insert the above inequalities in \rife{thetadu1},  choosing $\sigma$ suitably small and dropping positive terms we get
\be\label{mdu0}
\begin{split}
 -{\rm tr}\left({\mathcal A} \, {\mathcal D}^2 z \right) & + H_{pp}DV \cdot Dz+ \rho z -  {\rm tr}\left(  [\mD_{\eta} {\mathcal A} \cdot \mD z] \cD^2 u\right)
  +   \lam\,   |H_p(Du)\cdot Du - u_t|^2   
  \\
  & \qquad \qquad  \leq C_0  +  \lam \, (\beta_H+ \vep) \,  | D u|^q+ K (1+ \lam )  |Du|^{2(q-1)}\,.
%+ \frac{(1+\theta )}2 \alpha_H \left(|Du_t- H_p D^2 u|^2 + \vep \alpha_H |D^2 u|^2\right)
\end{split}
 \ee
We point out that the choice of $\sigma$ is fixed, by now. We compare now the function $z$ with 
$$
\vfi(t):= L\left( \frac1{t^2}+ \frac1{(T-t)^2}\right)\,.
$$
Since $\vfi$ blows-up at $t=0, T$, we have that $z-\vfi$ admits a maximum point in $(0,T)\times \oo$.
As in  Step 2 in Theorem \ref{bound},  
we can  show 
%there exists $K_0$ such that  $Dz \cdot \vec\nu <0$ if $z-\vfi$ has a maximum point with $|Du|\geq K_0$. But we have
%\begin{align*}
%z-\vfi  & =  \theta\, H(Du) + H(Du)-u_t + \lam u^2 - \vfi 
%\\
%& = \theta\, H(Du)  + \vep \log m + V(x)+ \lam u^2  -\vfi
%\leq \theta\, H(Du) + \vep \log (\|m\|_\infty) + \| V\|_\infty + \sigma
%\end{align*}
%which yields 
%%(suppose $\vep, \sigma \leq 1$)
%$$
%\theta\, H(Du) \geq z-\vfi  - [\vep \log (\|m\|_\infty) + \| V\|_\infty+\sigma] \,.
%$$
%Hence,  there exists $L_0>0$ such that  if $\max (z-\vfi )\geq L_0$ then $|Du|\geq K_0$
that $Dz\cdot \vec \nu \leq 0$ on $\partial \Omega$, and  no maximum could occur on the boundary. We now analyze maximum points inside the domain. 
In this case, we have $Dz=0$ and $[\mD_{\eta} {\mathcal A} \cdot \mD z]=0$ (because $\mathcal A$ does not depend on $u_t$), and we also have 
$$
{\rm tr}\left({\mathcal A} \, {\mathcal D}^2 z \right) \leq {\rm tr}\left({\mathcal A} \, {\mathcal D}^2 \vfi \right) = 6L\left( \frac1{t^4}+ \frac1{(T-t)^4}\right).
$$
Moreover, if $\max (z-\vfi )\geq L_0$, then using \rife{2theta} we have
$$
H_p\cdot Du - u_t =  H_p\cdot Du  - (1+\theta) H(Du) -\lam \frac{u^2}2 + z \geq \theta H(Du) + \vfi  + L_0 - \sigma-c_0
$$
which implies, for  $L_0\geq \sigma+ c_0$,
$$
|H_p\cdot Du - u_t|^2\geq \theta^2 H(Du)^2+ \vfi^2 \geq c_H \, |Du|^{2q} + L^2 \left( \frac1{t^4}+ \frac1{(T-t)^4}\right)\,.
$$
Therefore \rife{mdu0} implies, on the maximum point,
\begin{align*}
\lam c_H \, |Du|^{2q} + (\lam L^2 -6 L)\left( \frac1{t^4}+ \frac1{(T-t)^4}\right) & \leq C_0  + K (1+ \lam )  |Du|^{2(q-1)} + (\beta_H+\vep)\lam |Du|^q
\\
&  \leq  C_0+   C\left(\lam,  \frac  1{\lam}\right)+ \lam c_H |Du|^{2q} 
\end{align*}
which cannot hold if $L$ is sufficiently large.
%, depending oe.g. if $(\lam L^2-6L) \frac{16}{T^2}> C_0+   \tilde K(\lam+ \frac  1{\lam})$. 
The conclusion is that, with a suitable choice of $L$, we have
$$
z\leq L\left( \frac1{t^2}+ \frac1{(T-t)^2}\right) + L_0 
$$
%which yields
%In particular, a closer look at $L_0, L$ yield the estimate
%$$
%z \leq C (1+ T)  (1+ \|u\|_\infty^2) \left( \frac1{t^2}+ \frac1{(T-t)^2}\right) + \vep \log (\|m\|_\infty) + \| V\|_\infty + C\,.
%$$
%We have proved so far that
%$$
%(1+\theta) H(Du) - u_t \leq C (1+T) (1+ \|u\|_\infty^2)  \left( \frac1{t^2}+ \frac1{(T-t)^2}\right) + \vep \log (\|m\|_\infty) + \| V\|_\infty+ C
%$$
which implies, using the definition of $z$ and $m$:
$$
\theta H(Du) + \vep \log m  + V(x) \leq L\left( \frac1{t^2}+ \frac1{(T-t)^2}\right) + L_0\,.
%C (1+T) (1+ \|u\|_\infty)^2 \left( \frac1{t^2}+ \frac1{(T-t)^2}\right) + \vep\log (\|m\|_\infty) +  \| V\|_\infty+ C
$$
This yields
%\footnote{Notice that, by using the equation of $u_t$, we also have:
%$$
%\exp\left(\frac\theta{(1+\theta)\vep} u_t\right) m\, e^{V(x)/\vep}    \leq C (1+T^2) (1+ \|u\|_\infty)^2 \left( \frac1{t^2}+ \frac1{(T-t)^2}\right) + \log (\|m\|_\infty) + \| V\|_\infty
%$$
%} 
\rife{mduloc}.
%$$
%\theta H(Du) + \vep \log \left(\frac m{\|m\|_\infty}\right) + V(x) \leq \| V\|_\infty+  \frac{K_T}{t^2(T-t)^2}
%$$
%$$
% \rightsquigarrow \quad       e^{\theta  H(Du)} m^\vep\, e^{V(x)}\leq  \|m\|_\infty^\vep   C(\| V\|_\infty) \exp[\frac {K_T}{t^2(T-t)^2} ]
%$$
%for some $K_T$ only depending on $T, \|u\|_\infty$.
\qed

\begin{remark} As we mentioned, the above estimate says that, once we drop the positivity condition on the marginals,  $u$ should remain smooth on the support of $m$.
But unfortunately, we do not have precise informations on the behavior of the support of the solution in that case (at least in dimension $d>1$).
% It is an open problem to understand the behavior of the support
%of $m(t)$ in the case that  $\vep>0$ (entropic problem) and $m_0, m_1$ have compact support. 
This issue is possibly  related to  the regularity
of solutions in critical sets for degenerate quasilinear elliptic problems. 
\end{remark}
\vskip2em

\paragraph{Convergence to weak solutions.} \qquad  Since \cite{Carda1}, \cite{CG}, \cite{CGPT}, a quite general theory of weak solutions is available for mean-field game systems  as 
\begin{equation}\label{mfg0}
\begin{system}
&-u_t + H( Du)  = f(m)  + V(x)  &&\text{in }  Q:= (0,T)\times \Omega, \\
&m_t  - \dive (m\,H_p(Du))  = 0  &&\text{in } Q,
\\
& H_p(Du)\cdot \vec\nu=0 &&\text{on $(0,T)\times \partial \Omega$,}
\\
&m(0,\cdot) = m_0, \; m(T,\cdot)  = m_1 &&\text{in } \Omega\,,
\end{system}
\end{equation}
whenever $H$ is convex and $f$ is nondecreasing.  The theory was actually developed when a terminal condition is prescribed on $u$ (a final pay-off in the cost functional), rather than for the  transport problem in which both marginals are imposed on $m$ (at $t=0,t=T$); but  the transport case  was also addressed (see \cite{OPS}, \cite{GMST}) relying on the variational interpretation of those solutions as relaxed minima of the corresponding functionals. 
Up to minor differences (related to possibly different growth conditions on $H$ and $f$), weak solutions are defined as follows. For simplicity, we consider here the case that the Hamiltonian has quadratic-like growth. Let us recall that $\cP(\Omega)$ denotes the space of probability measures, endowed with the  Kantorovich-Rubinstein distance 
$$
d_1(m,m')=\sup_{\phi} \into \phi \,d(m-m') \qquad \forall m,m'\in  \cP(\Omega),
$$
where the supremum is taken over all $1-$Lipschitz continuous maps $\phi:\Omega\to \R$.

\begin{definition}\label{def-weak} A pair  $(u,m)$ is a weak solution  of \rife{mfg0} if $m\in C^0([0,T];\cP(\Omega))\cap L^1(Q)$ with $m(0)= m_0$, $m(T)=m_1$,  $u \in L^2_{loc}((0,T);H^1(\Omega))$ and in addition  $m\, |Du|^2\in L^1(Q)$, 
$f(m)m \in L^1(Q)$ and $(u,m)$ satisfy:
\vskip0.3em
(i) $u$ is a weak sub-solution satisfying, in the sense of distributions,  
 $$
 -u_t + H( Du) \leq  f(m)  + V(x) \qquad \hbox{in $Q$}
 $$
\vskip0.3em
(ii) $m$ is a weak solution satisfying, in the sense of distributions, the continuity equation
$$
m_t  - \dive (m\,H_p(Du))  = 0  \qquad \hbox{in $Q$}
 $$
\vskip0.3em
(iii) $(u,m)$ satisfy the identity
\begin{equation}\label{en-id}
\begin{split}
 \into   m_0\, u(0)\, dx & - \into u(T)\,   m_1\,dx = \int_0^T\!\!\! \into  f(m) m\, dxdt
%\\  & 
+  \int_0^T\!\!\! \into  m\,  \left[   H_p(  D  u)\cdot Du-  H(  Du)\right]    dxdt
\end{split}
\end{equation}
where $u(0), u(T)$ are the one-sided traces of $u$ (which are well defined by properties of sub-solutions of Hamilton-Jacobi equations, see \cite[Section 5]{OPS}).
% which are assumed to satisfy $u(0)m_0, u(T)m_1\in L^1(\Omega)$.
\end{definition}

 One motivation  in the construction of smooth solutions (Theorem \ref{exi-smooth}) is that it provides a regularization for weak solutions of \rife{mfg0}, by adding a small entropy term in the coupling function $f$. In particular, this regularization allows one to justify 
the estimates proved  in Section \ref{disconv}, which  eventually holds for weak solutions in the limit as $\vep \to0$.  We give below a sample statement of this kind.
This is to be compared with the results in \cite{LaSa}, where similar estimates are justified for weak solutions by using time-discretization and the geodesic interpretation of $m(t)$ from optimal transport theory, whereas our approach is different and entirely  relies   on the PDE Eulerian approach  developed in Section \ref{disconv}. Let us stress that weak solutions, as defined above, coincides with minima of the corresponding functional  \rife{varia} (with $\vep=0$).

\begin{theorem}\label{weaksol}  Let $\Omega$ be a  $C^3$, bounded, convex domain  in $\R^d$ and let $V\in W^{1,\infty}(\Omega), m_0,m_1\in L^\infty(\Omega)\cap \cP(\Omega)$.
Let $H$ be a $C^2$ function satisfying   \rife{hrad} and \rife{Hpp}-\rife{hppp}.  Assume that     $f\in C^1([0,\infty))$ is an  increasing function satisfying condition \rife{coercf}. Then there is  a unique $m$ and a unique $u$ (up to $m$-negligible sets) such that   $(u,m)$ is a weak solution of problem \rife{mfg0}, with $\into u(T)m_1=0$. 

Moreover we have that $m\in L^\infty(Q)$, it is the unique minimum of the functional \rife{varia} (with $\vep=0$)  and satisfies estimates \rife{supglobm}--\rife{suplocm}. Finally $m$ is the limit as $\vep \to 0$ of the minima $m^\vep$ of regularized problems, obtained in Theorem \ref{exi-smooth}.
\end{theorem}

\proof We only sketch the argument. For $\vep>0$, we consider the solution $(u^\vep, m^\vep)$ of \rife{mfpp}, where $H,f,V$ are suitably regularized in order that Theorem \ref{exi-smooth} can be applied. We take for $u^\vep$ the normalized version given by Lemma \ref{norm},  i.e.  $\into u^\vep(T) m_1dx =0$. Applying 
Proposition \ref{supm}, we have that $\|m^\vep\|_\infty$ is uniformly bounded. Moreover, by estimate \rife{hatu2} in Lemma \ref{norm}, we also have that $u^\vep$ is locally uniformly bounded. We notice that $u^\vep$ satisfies, for some constant $C$,
\be\label{subeps}
-u_t^\vep + H(Du^\vep) \leq f(m^\vep)+ V + \vep C\,, 
\ee
where the right-hand side is bounded in $L^\infty(Q)$. Now we use some arguments from \cite{OPS} and \cite[Section 1.3.7]{CP-CIME}. Up to subsequences, we may assume that $u^\vep$ converges  to some function $u$ (at least weakly in $L^2((a,b); H^1(\Omega))$ for all $0<a<b<T$), and that, for a.e. $t\in (0,T)$,   $u(t)$ is the limit  of $u^\vep(t)$ in the weak-$*$ topology of $L^\infty(\Omega)$. By stability of sub-solutions (see e.g. \cite[Thm 5.11]{OPS}), we have that $u$ satisfies
\be\label{usub}
-u_t  + H(Du ) \leq \bar f+ V   
\ee
where $\bar f$ is the weak-$*$ limit of $f(m^\vep)$ in $L^\infty(Q)$. By property of subsolutions, $u$ has  one-sided limits as $t\to 0^+$, $t\to T^-$; since from \rife{subeps} $u^\vep$ satisfies
$$
\into u^\vep(t) m_1 \leq C (T-t)
$$
we also deduce, letting first $\vep \to 0$ and then $t\to T$, that $\into u(T)m_1\leq 0$. As for $m$, if $w$ denotes the weak limit (e.g. in $L^2(Q)$) of $m^\vep H_p(Du^\vep)$, then $m$ is a bounded weak solution of the continuity equation
\be\label{mcont}
\begin{cases}
m_t - \dive (w) =0 & \\
m\in C([0,T], \cP(\Omega))\,,\quad m(0)= m_0\,,\, m(T)= m_1\,. & 
\end{cases} 
\ee
The duality between \rife{usub} and \rife{mcont} was exploited in \cite[Section 5.3]{OPS}; roughly speaking, it is possible to cross multiply \rife{usub} and \rife{mcont}. 
Using that  $\into u(T)m_1\leq 0$, and reasoning as in \cite[Thm 1.15]{CP-CIME}, then one shows  
$$
\limsup_{\vep \to 0} \int_0^T \into f(m^\vep)m^\vep \leq \into \bar f \, m\, dxdt\,\,.
$$ 
A standard Minty-type argument implies that $\bar f= f(m)$, and since $f$ is increasing, one also deduces that $m^\vep\to m$ a.e. in $Q$, and therefore in $L^p(Q)$ for every $p<\infty$. Moreover, the identification of $\bar f$ goes together with the identification of $w= m\, H_p(Du)$, using the convexity of $(m,w)\mapsto L(w/m)m$. 
Hence one follows arguments in \cite{OPS} or \cite{CP-CIME}: $(u,m)$ satisfy the equality
$$
\int_0^T \into m  \left( H_p(Du  )Du  - H(Du)\right) + (f(m)+ V) m = \into m_0\, u(0)\, dx
$$
and $\into u(T)m_1\,dx $ is proved to be zero. It follows from the above equality that $m$ is a minimum of the limit functional and that $(u,m)$ is a weak solution of   \rife{mfg0} in the sense of Definition \ref{def-weak}.  The uniqueness of $m$ and of $u$ up to $m$-negligible sets can be proved   as in \cite[Thm 1.16]{CP-CIME}. 
\qed

\paragraph{The compact case.}  All results proved here holds in more generality if the state space is the flat torus ${\mathbb T}^d= \frac{\R^d}{\Z^d}$. In fact, requiring the Hamiltonian to be radial was only needed to handle the regularity near the boundary $\partial\Omega$. Therefore, all results stated remain  true in the torus removing the assumption \rife{hrad}.  In the same spirit, the results could be extended to a compact Riemannian manifold without boundary, although this requires to use Bochner's formula to handle the gradient estimates in the Riemannian setting.  We will exploit this case in future work.

%\section{Comments, extensions and open problems}

\section{Appendix:  existence of solutions to the elliptic problem}

Here we show the existence of solutions to problem \rife{eps-de} where, for simplicity, we fix $\vep=1$, $T=1$. In order to construct a solution, we follow the classical approach and we build a family of one-parameter problems. For $\tau\in [0,1]$, we set
\be\label{mtau}
m^\tau := \vfi^\tau( -u_t+ H(Du)- \tau V(x))\,,\qquad \vfi^\tau:= (\tau\, f(\cdot)+ \log(\cdot) )^{-1}
\ee
and we consider the problem
\begin{equation}\label{LU}
\begin{cases}
   -{\rm tr} ( {\mathcal A}^\tau(x, \cD u) \cD^2 u) 
%-u_{tt} + 2  H_p(Du) Du_t -   D^2 u H_p(Du) \cdot H_p(Du) - ( 1+ \tau\, m^\tau\, f'(m^\tau ) ){\rm tr}\left(H_{pp}(Du)D^2u \right)   
+ \tau DV(x) \cdot H_p(Du) + \rho u= 0
& 
\hbox{in  $(0,1)\times \Omega$,} 
 \\
  - u_t + H(Du)- \tau (f(m_1)+\tilde V(x)) - \de u = \psi_1 (x)   
%f^\vep(m_1)  
 &
\hbox{at $t=1$, $x\in \Omega$,}
\\
    -u_t   + H(Du)- \tau (f(m_0)+\tilde V(x)) + \de u=   \psi_0 (x)
 %f^\vep (m_0)  
  &
 \hbox{at $t=0$, $x\in \Omega$,}
\\
  H_p(Du)\cdot \vec\nu=0  &
\text{on $(0,1)\times \partial \Omega$,}
\end{cases}
\end{equation}
 where 
 $$
{\mathcal A}^\tau= \left(a_{ij}^\tau((t,x),\mD u)\right) : = 
\begin{pmatrix} 1 & - H_p(Du) \\ - H_p(Du) &   H_p(Du)\otimes H_p(Du)\end{pmatrix}
+( 1+ \tau\, m^\tau\, f'(m^\tau ) ) \begin{pmatrix} 0 & 0 \\ 0 &  
H_{pp}\end{pmatrix}\,
$$
and where  $\tilde V, \psi_0, \psi_1\in C^{1,\alpha}(\oo)$. Problem \rife{eps-de} actually corresponds to $\tau =1$, $\tilde V= V$, $\psi_0= \log(m_0), \psi_1= \log(m_1)$.
%It is however convenient to 

 %$f^\vep(s)= f(s) + \vep \log s$. 

To this problem we are going to apply Theorem 1.2 (Chapter X, page 462) in \cite{LU}. Unfortunately, we can not apply directly this result because our domain does not look regular in the time-space environment; otherwise said,  our nonlinear Neumann boundary  condition  is  only piecewisely defined at the $(t,x)$ boundary. 
To overcome this technical issue,  we will use the radial structure of the Hamiltonian $H$ (see \rife{hrad}).
%
%To overcome this problem, here we assume that $H$ is a radial function, namely
%\be\label{radial}
%\exists \quad  \hbox{$h\in C^2([0,\infty))$, increasing and convex: }\quad H(p)= h(|p|)\,.
%\ee 
%In particular, $H$ satisfies \rife{Hpp} if $h'(0)=0$ and $\alpha\leq h''\leq \beta$.  
This  assumption  simplifies  the Neumann condition on $\partial\Omega$ and allows us to follow the classical steps up to using an extra reflection argument, in order to infer the second order estimates. A similar result for general Hamiltonians (without requiring \rife{hrad}) would need  to develop different  technical  tools to handle the $C^1$ estimates for piecewisely defined nonlinear  boundary value problems; this is beyond our present scopes.
\vskip0.5em
We recall that this kind of quasilinear problems (with nonlinear boundary conditions) can be solved provided one shows (uniformly in $\tau$):
 \begin{itemize}
 
 \item A global Lipschitz bound for (smooth) solutions of \rife{LU}
 
 \item A  global $C^{1,\alpha}$ estimate for (smooth) solutions of \rife{LU}
 
 \item Solvability and $C^{2,\alpha}$ estimates for the linearized version of \rife{LU}

 \end{itemize} 

We split those three tasks in the next propositions.

\begin{proposition}\label{boundK} Assume that  $V\in C^{2,\alpha}(\overline \Omega)$,  that $H\in C^{3,1}(\R^d)$ satisfies \rife{hrad} and \rife{Hppq}-\rife{hpppq} for some $q>1$ and $\varpi>0$, and that $f\in C^{2,1}(0,\infty)$ satisfies
\be\label{fprime-mod}
\exists\,\, \beta, s_0>0\,:\, | f''(s)|\, s^2 \leq \beta \,    (1+ f'(s)s)  \qquad \hbox{ for all $s\geq s_0$. }
\ee
Let  $m_0, m_1, \tilde V, \psi_0,\psi_1 \in C^{1,\alpha}(\oo)$. Then   there exists $M$, independent of $\tau\in [0,1]$, such that any $u\in C^{2,\alpha}(\oQ) $ which is a solution of \rife{LU} satisfies
\be\label{tauM}
\|u\|_\infty + \| \mD u\|_\infty \leq M
\ee
where $M$ depends on $\de, \|m_0\|_{W^{1,\infty}(\Omega)}, \|m_1\|_{W^{1,\infty}(\Omega)}, \|\psi_0\|_{W^{1,\infty}(\Omega)}, \|\psi_1\|_{W^{1,\infty}(\Omega)}, \|\tilde V\|_{W^{1,\infty}(\Omega)}, \|V\|_{W^{2,\infty}(\Omega)}$ and on $H, f$. 
\end{proposition}

\proof This is the result of  Theorem \ref{bound} (together with Lemma \ref{sup}), up to minor remarks. The vector field $V$ is here replaced by $\tau V(x)$ in the interior and by $\tau \tilde V$ on the boundary, but this makes no difference  in the proof.  The condition \rife{fprime-mod} is slightly more restrictive than  \rife{fprime} in Theorem \ref{bound}; this way  we can assert that $\tau f(\cdot)$ satisfies the same condition as $f$, with the same constant.  Hence, the gradient estimates are independent of $\tau$. We point out that the regularity required on $V, H, f$ guarantee  that any $C^2$ solution $u$ actually belongs to $C^{3}(Q)$, so that the gradient bounds can be classically derived.  Finally, we notice that the bound on $u$ depends on $\de$ from maximum principle, as in Lemma \ref{sup}. 
\qed

\vskip1em

\begin{proposition}\label{boundaryHol} Under the same assumptions of Proposition  \ref{boundK},  there exists $C$, independent of $\tau\in [0,1]$, such that any $u\in C^{2,\alpha}(\overline Q)$ which is a solution of \rife{LU} satisfies
\be\label{gradhol}
\|u\|_{C^{1,\alpha}(\oQ)} \leq C
\ee
where $C$ depends on the constant $M$ in \rife{tauM} (as well as  on the bounds of $H_p, H_{pp}, H_{ppp}$   for $|p|\leq M$) and again on $\|m_0\|_{W^{1,\infty}(\Omega)}, \|m_1\|_{W^{1,\infty}(\Omega)}, \|\psi_0\|_{W^{1,\infty}(\Omega)}, \|\psi_1\|_{W^{1,\infty}(\Omega)}, \|\tilde V\|_{W^{1,\infty}(\Omega)}$.
\end{proposition}

\proof  The elliptic problem in  \rife{LU} can be rephrased as
\be\label{ellip}
\begin{cases}
L_\tau (u):= -a_{ij}^\tau(\eta,\mD u) u_{\eta_i\eta_j} + \tau DV(x) \cdot H_p(Du) + \rho u =0
&  \hbox{ in $Q$,} \\
N^\tau(u)= (\psi_1, \psi_0)  & \hbox{on $\Sigma$,}
\\
 H_p(Du)\cdot \vec\nu=0  &
\text{on $(0,1)\times \partial \Omega$,}
\end{cases}
\ee
where $\eta=(t,x)$, $\Sigma= (\{1\}\times \Omega) \cup (\{0\}\times \Omega) $  and the boundary operator $N^\tau$ is defined  as
%$$
%N_\tau (u):= \begin{cases} 
%N_1(x,u, \mD u) = - f^\vep(m_1)   &\hbox{at $t=1$, $x\in \Omega$,}
%\\
% N_0(x,u, \mD u)  =   f^\vep (m_0)  &\hbox{at $t=0$, $x\in \Omega$,}
%\\
%  D_xu\cdot \vec\nu(x)  =0 &\hbox{on $(0,1)\times \partial \Omega$,}
%\end{cases}
%$$
%where, we recall, $\eta=(t,x)$.
\be\label{ntau}
N_\tau (u):= \begin{cases} 
 - u_t   + H(Du) - \tau (f(m_1)+\tilde V(x)) -\de u   &\hbox{at $t=1$, $x\in \Omega$,}
\\
  -u_t  +H(Du) - \tau (f(m_0)+\tilde V(x)) + \de u    &\hbox{at $t=0$, $x\in \Omega$.}
\end{cases}
\ee
Due to the estimate \rife{tauM} in  Proposition \ref{boundK}, we have that the coefficients $a_{ij}^\tau$ need only be considered in the compact set $ \mathcal M:= \{\eta\in \oQ, |u|\leq M, |\mD u|\leq M \}$.  
Assuming $V$ of class $C^{2,\alpha}$, $f$ of class $C^{2,1}$  and $H$ of class $C^{3,1}$, we get that the coefficients $a_{ij}^\tau$, as well as the first order coefficients, are $C^{1,\alpha}$ in their arguments,  with $\frac{\partial a_{ij}^\tau(y,p)}{\partial p}$ being Lipschitz in $p$.  In particular, for the interior regularity we can apply the classical results (see e.g. \cite[Chapter 6, Thm 1.1]{LU})  and we deduce the $C^{1,\alpha}$ interior estimate for  $u$.

We are only concerned now with the boundary regularity. 
%{Similarly, the boundary functions $N_0, N_1$ will be $C^{2,\alpha}$ with $\frac{\partial^2 a_{ij}^\tau(y,u,p)}{\partial p^2}$ being Lipschitz in $(u,p)$.}
Without loss of generality (because $H_{pp}$ is bounded below),  we assume here that 
\be\label{rad}
H(p)= h(|p|^2)\,,
\ee 
for some $C^2$ function $h$ such that $h'>0$.  With this notation we have  $H_p= 2h'\,p$ and $H_{pp}= 2h' I_d+4h''(p\otimes p)$.
\vskip1em
Now we use a (local) reflection argument in the $x$ variable.  To this purpose, recall that for $x\in \Omega$, we denote by $d(x)$ the distance of  $x$ to the boundary. It is well known (see e.g. \cite{GT}) that if $\Omega$ is of class $C^k$, $d(x)$ is a $C^k$ function in a neighborhood of the boundary.  More precisely,  there exists a positive number $d_0>0$ such that, if we denote $ \Omega_{d_0}:= \{x\in \Omega \,:\, d(x) < d_0\}$, then any point $x\in \Omega_{d_0}$ can be represented as
\be\label{proj}
 x=  - \delta\, \vec{n}(\bar x) + \bar x
\ee
where $\delta= d(x)$ and $\bar x$ is the  projection of $x$ onto $\partial \Omega$ ($\vec n(\bar x)$ is the outward unit vector at $\bar x$). To ease notations, we denote hereafter by $N$ the dimension of the euclidean space (i.e. $\Omega\subset \R^N$). Therefore, if we represent $\bar x\in \partial \Omega$ in a system of coordinates of some local chart, we may assume that $\bar x= (y', g(y'))$  for some $C^k$ function $g$ defined on some open subset $U\subset \R^{N-1}$. Hence \rife{proj} yields a local diffeomorphism  between $\Omega_{d_0}$ and $U\times (0,d_0)$ defined as
$$
T(x):= (y', \delta)  \quad \iff \quad x=  - \delta\, \vec n(\bar x) + \bar x\,,\quad \bar x= (y', g(y')) \in \partial\Omega.
$$
Let us now set $y_N:=\delta$ so that $y=(y',y_N)$ and $y= T(x)$ will shortly denote the above change of variables. We notice that  assuming  $ \Omega$   of class  $C^3$ guarantees that the mapping $T$ is of class $C^2$.  The advantage of this change of coordinates lies in the property that
\be\label{ortho}
\nabla T_k(x) \cdot \nabla T_N(x)=\nabla T_k(x) \cdot \nabla d(x)= 0 \quad \forall k=1,\ldots, N-1
\ee
because $T_N(x)= d(x)$ (the distance function) and the projection coordinates remain constant if we move along the normal direction (i.e. $\nabla T_k(x)\cdot \nabla d(x)=0$ for all $N-1$-coordinates $y_k=T_k(x)$). 

Now we set
$$
\tilde u (t,y):= u(t, T^{-1}y) \,\qquad   \rightsquigarrow \qquad u(t,x)= \tilde u(t, T(x))\,
$$
and  we compute:
\be\label{uxi}
u_{x_i}=  \sum_{k=1}^N\tilde u_{y_k} \frac{\partial T_k}{\partial x_i} \quad ; \quad  u_{x_ix_j}= \sum_{k,\ell=1}^N\tilde u_{y_ky_\ell} \frac{\partial T_k}{\partial x_i}\frac{\partial T_\ell}{\partial x_j} + \sum_{k=1}^N\tilde u_{y_k}\frac{\partial^2 T_k}{\partial x_i\partial x_j}
\ee
so the equation \rife{ellip} is transformed into
\be\label{tilde-eq}
\tilde u_{tt} - 2    \sum_{k } \left( DT_k\cdot  H_p(Du)\right)  \tilde u_{t y_k} +     \sum_{k,\ell=1}^N  \tilde a^\tau_{kl}  \tilde u_{y_ky_\ell} + \sum_{k=1}^N\tilde  b_k^\tau\, \tilde u_{y_k} 
%+ \tau DV(x) \cdot H_p(Du) 
= \rho \tilde u
\ee
where
%$$
%\sum_{i,j} \tilde a_{ij}^\tau((t,y) , \tilde u, \mD \tilde u )\sum_{k,\ell=1}^N   \tilde u_{y_ky_\ell} \frac{\partial T_k}{\partial x_i}\frac{\partial T_\ell}{\partial x_j} + \sum_{k=1}^N\tilde u_{y_k}\frac{\partial^2 T_k}{\partial x_i\partial x_j}
%+  \sum_i  b_i(T^{-1}y) \sum_{k=1}^N\tilde u_{y_k}\frac{\partial T_k}{\partial x_i}=0
%$$
%where
\be\label{akl}
\tilde a^\tau_{kl} := \sum_{i,j}  a_{ij}^\tau\, \frac{\partial T_k}{\partial x_i}\frac{\partial T_\ell}{\partial x_j} \,; \qquad \tilde b_k^\tau:= \sum_{i,j}a_{ij}^\tau\,\frac{\partial^2 T_k}{\partial x_i\partial x_j} 
- \tau 2h' DV\cdot DT_k  
\ee
and where the coefficients $a_{ij}^\tau$ are computed on $Du$; recall that   $Du= M D\tilde u$ for $M=(m_{ij})= \left(\frac{\partial T_j}{\partial x_i}\right)$.

We extend now  $\tilde u$ for $y_N<0$ by setting
$$
u^*(t,y):= \begin{cases} \tilde u(t,y',y_N) & \hbox{if $y'\in U, y_N\in (0,d_0)$}
\\
\tilde u(t,y', -y_N) & \hbox{if $y'\in U, y_N\in (-d_0,0)$}
\end{cases}
$$
Notice that \rife{ortho} implies
$$
D u(t,x)\cdot D d(x)= \sum_{k=1}^N \tilde u_{y_k} (DT_k\cdot Dd(x))= \tilde u_{y_N}\, |Dd(x)|^2= \tilde u_{y_N}\,\,,\qquad x\in \Omega_{d_0}
$$
so the Neumann condition $H_p(Du)\cdot \vec\nu=0$ together with \rife{rad} translate into  $\tilde u_{y_N}=0$; this implies that $u^*$ is a $C^1$ function. 
Now we look at the equation satisfied by $u^*$;  this involves the even reflection (across $y_N=0$) of coefficients $\tilde a^\tau_{kl}$ for $k,\ell \neq N$ or $k=\ell=N$ and the odd reflection of coefficients $\tilde a^\tau_{kl}$ for $k\neq \ell$ and $k$ or $\ell=N$. Let us check that the structure condition \rife{rad} and the property \rife{ortho} yield a nice equation for $u^*$. First of all we observe that, due to \rife{ortho}, 
\begin{align*}
|Du|^2 & = \sum_{i} u_{x_i}^2= \sum_{i}     \left( \tilde u_{y_N}  \frac{\partial T_N}{\partial x_i}     + \sum_{p=1}^{N-1}       \tilde u_{y_p}   \frac{\partial T_p}{\partial x_i} \right)^2
\\
%& =  \sum_{i}     \left\{ \tilde u_{y_N}^2  \left|\frac{\partial T_N}{\partial x_i} \right|^2    +  \left( \sum_{p=1}^{N-1}       \tilde u_{y_p}   \frac{\partial T_p}{\partial x_i} \right)^2  +2 \sum_{p=1}^{N-1}    \tilde u_{y_N}  \frac{\partial T_N}{\partial x_i}    \tilde u_{y_p}   \frac{\partial T_p}{\partial x_i}\right\}
%\\
& = \tilde u_{y_N}^2 + \sum_{i}  \left( \sum_{p=1}^{N-1}       \tilde u_{y_p}   \frac{\partial T_p}{\partial x_i} \right)^2  
 =  |u^*_{y_N}|^2+ |M'D'u^*|^2 \\
&\qquad \hbox{for $ M'_{ip} =  ( \frac{\partial T_p}{\partial x_i})\,,\,\, D'u^*= (u^*_{y_p})_{p=1,\ldots,N-1}$}\quad 
\rightsquigarrow \quad    (M' D'u^*)_i: =  \sum_{p=1}^{N-1}       u^*_{y_p}   \frac{\partial T_p}{\partial x_i}
\end{align*}
In particular, by definition of $m^\tau$ in \rife{mtau},  we have
\be\label{kastar}
( 1+ \tau\, m^\tau\, f'(m^\tau ) ) = {\mathcal \kappa}^\tau (V^*(y), \cD u^*)
\ee
for some smooth function ${\mathcal \kappa}^\tau$, where   $V^*(y)$ is the even reflection of $V(T^{-1}(y))$.
 
Then we compute from \rife{akl}, using the definition of $a_{ij}^\tau$ and the radial structure of $H$, and due to condition \rife{ortho},
\begin{align*}
\forall k\neq N,\qquad  \tilde a^\tau_{kN} & = \sum_{i,j}   \left[ H_{p_i}\otimes H_{p_j} + H_{p_ip_j}  ( 1+ \tau\, m^\tau\, f'(m^\tau ) ) \right] \frac{\partial T_k}{\partial x_i}\frac{\partial T_N}{\partial x_j}  
\\
& = \sum_{i,j}   \left( 4(h')^2+ 4 h''  {\mathcal \kappa}^\tau (V^*(y), \cD u^*) \right)   u_{x_i}u_{x_j} \frac{\partial T_k}{\partial x_i}\frac{\partial T_N}{\partial x_j}  
%\\
%& = \sum_{i,j,p,q}  \tau^2  \tilde u_{y_p} \tilde u_{y_q} \frac{\partial T_p}{\partial x_i}\frac{\partial T_q}{\partial x_j}    \frac{\partial T_k}{\partial x_i}\frac{\partial T_N}{\partial x_j}  
%\\
%& = \sum_{i,j,p\neq N}  \tau^2  \tilde u_{y_p} \tilde u_{y_N} \frac{\partial T_p}{\partial x_i}\frac{\partial T_N}{\partial x_j}    \frac{\partial T_k}{\partial x_i}\frac{\partial T_N}{\partial x_j}  
\\
& =  \left(4(h')^2+ 4 h''  {\mathcal \kappa}^\tau (V^*(y), \cD u^*) \right) \sum_{p=1}^{N-1}   DT_p\cdot DT_k\tilde u_{y_p} \tilde u_{y_N}
\end{align*}
where we used \rife{uxi} (and once more \rife{ortho}) in the last step. In particular, we notice that $\tilde a^\tau_{kN} $ is an odd coefficient with respect to reflection; hence, using the above notation, 
$$
\forall k\neq N,\qquad  \tilde a^\tau_{kN} \tilde u_{y_ky_N}= \left(4(h')^2+ 4 h''  {\mathcal \kappa}^\tau (V^*(y), \cD u^*) \right) (DT_k\cdot MD'u^*) u^*_{y_N} u^*_{y_ky_N}\,.
$$
With similar computations we find
\begin{align*}
  \tilde a^\tau_{NN} & = \left(4(h')^2+ 4 h''  {\mathcal \kappa}^\tau (V^*(y), \cD u^*) \right) (u^*_{y_N})^2 + 2h'\,  {\mathcal \kappa}^\tau (V^*(y), \cD u^*)\,.
\end{align*}
Still using \rife{rad}, \rife{uxi} and \rife{ortho}, we observe that
\begin{align*}
 \sum_{k} \left( DT_k\cdot  H_p(Du)\right)  \tilde u_{t y_k} & = 2h'  ( \tilde u_{y_N} \tilde u_{ty_N} + \sum_{k,\ell=1}^{N-1} \left( DT_k\cdot DT_\ell )  \tilde u_{y_\ell}   \tilde u_{t y_k}\right)
 \\
 & = 2h'( u^*_{y_N}   u^*_{ty_N} + \sum_{k,\ell=1}^{N-1} \left( DT_k\cdot DT_\ell\right)    u^*_{y_\ell}     u^*_{t y_k})
\end{align*}
Therefore, $u^*(t,y)$ solves  the following elliptic equation for $t\in (0,1), y\in B$, where $B\subset \R^N$ is any  open set  which is contained in $U\times (-d_0,d_0)$:
\be\label{reflected}
 u_{tt}^* - 4 h' ( u^*_{y_N}   u^*_{ty_N} + \sum_{k,\ell=1}^{N-1} \left( DT_k\cdot DT_\ell\right)    u^*_{y_\ell}     u^*_{t y_k}) +     \sum_{k,\ell=1}^N   a^*_{kl}   u^*_{y_ky_\ell} 
%+ \sum_{k=1}^N\tilde  b^*_k(y)\, u^*_{y_k} 
+ b^*(t,y)= \rho   u^*
\ee
where
\begin{align*}
& {\small a^*_{kl} =  \left(4(h')^2+ 4 h''  {\mathcal \kappa}^\tau (V^*(y), \cD u^*) \right)(DT_k\cdot M'D'u^*) (DT_\ell\cdot M'D'u^*) +2h'\,  {\mathcal \kappa}^\tau (V^*(y), \cD u^*)(DT_k\cdot DT_\ell)  }  \quad \hbox{if $k, \ell\neq N$,}
\\
& a^*_{kN} =  \left(4(h')^2+ 4 h''  {\mathcal \kappa}^\tau (V^*(y), \cD u^*) \right) \left(DT_k\cdot M'D'u^*\right)  u^*_{y_N}  \quad \hbox{if $k\neq N$,} \\
&  a^*_{NN} =   \left(4(h')^2+ 4 h''  {\mathcal \kappa}^\tau (V^*(y), \cD u^*) \right) |u^*_{y_N}|^2+ 2h'\,  {\mathcal \kappa}^\tau (V^*(y), \cD u^*)
%\\
%& M_{ip} =  ( \frac{\partial T_p}{\partial x_i})\,,\,\, D'u^*= (u^*_{y_i})_{i=1,\ldots,N-1}\quad \rightsquigarrow \quad    (M D'u^*)_i: =  \sum_{p=1}^{N-1}       u^*_{y_p}   \frac{\partial T_p}{\partial x_i}
\end{align*}
%and we denoted by $V^*(y)$ the even reflection of $  V(T^{-1}(y))$, 
and  where $b^*(t,y) $ is the even extension (through $y_N=0$) of the term $\tilde b^\tau\cdot D\tilde u$. 
%defined as
%$$
%b^*(t,y)= \begin{cases}
%\sum_{k=1}^{N-1}  b_k^\tau(t, T^{-1}(y)) \,  u^*_{y_k} + b_N^\tau(t, T^{-1}(y)) u^*_{y_N} & \hbox{if $y_N\geq 0$}, 
%\\
%\sum_{k=1}^{N-1}  b_k^\tau(t, T^{-1}(y)) \,  u^*_{y_k} - b_N^\tau(t, T^{-1}(y)) u^*_{y_N} & \hbox{if $y_N<0$}
%\end{cases}
%$$
%$$
%b_k^*:= \begin{cases} \tilde b_k^\tau (t,y',y_N) & \hbox{if $y_N>0$}
%\end{cases}
%\sum_{i,j}a_{ij}^\tau((t,T^{-1}y) , \tilde u, \mD  u)\frac{\partial^2 T_k}{\partial x_i\partial x_j}  + (\vep + \tau h(-u^*_t+ \frac\tau2  \left(|u^*_{y_N}|^2+ |MD'u^*|^2\right)-\tau V^*(y)) \Delta T_k(y)
%+     \sum_i  b_i^*(y)  \frac{\partial T_k}{\partial x_i}
%$$
\vskip0.5em
Now we apply to equation \rife{reflected} the $C^{1,\alpha}$-estimates (up to the boundary)  for elliptic equations with nonlinear first order condition at the boundary, see e.g. \cite[Lemma 2.3]{lieb84NA}\footnote{We warn the reader that, literally, the result in \cite[Lemma 2.3]{lieb84NA} assumes the coefficients $a_{ij}(y, \mD u)$ to be $C^1$ with respect to $y$; however this is not necessary, as can be readily checked by inspecting the proof as well as by noticing that the estimate only requires $\frac{\partial a_{ij}}{\partial y}$ to be bounded. In fact, the Lipschitz character of $a_{ij}$ with respect to $y$ is enough, which is here preserved by the  reflection argument.}. 
Let us point out that the boundary  here is represented by the hypersurfaces $(\Sigma_0\cap B) \cup(\Sigma_1\cap B)$; the solution $u^*$ satisfies 
a uniform bound $|u^*| + |\mD u^*|\leq K$ (by Proposition \ref{boundK}) and the 
coefficients $a_{ij}^*((t,y), \mD u^*)$  are  $C^1$ with respect to  $\mD u^*$  and $\frac{\partial a_{ij}^*}{\partial q}$ as well as $\frac{\partial a_{ij}^*}{\partial y}$ are uniformly bounded (because $H$ is $C^{3,1}$ in $p$ and $V$ is Lipschitz). Note that $a_{ij}^*$ depends on $y$ only through (the even reflections of) $V(y)$ and $T^{-1}(y)$. %Ellipticity is verified (uniformly in $\tau$) because, setting $\eta=(s,\xi)$, we have
%\begin{align*}   
%a_{ij}^* \eta_i \cdot \eta_j   &= \left\{ s -2 \tau  \left( u^*_{y_N} \xi_N+ \sum_{k=1}^{N-1} \left(DT_k\cdot M D'u^*\right)\xi_k \right) \right\}^2  + \vep |\xi|^2
%\\ &  \geq \lambda_{\vep, K} (s^2+ |\xi|^2) 
%\end{align*}
%for some $\lambda_{\vep, K}$ depending on $\vep, K$ and $\|DT\|_\infty$. 
As for the coefficient $b^*(t,y) $, it is actually continuous (thanks to the fact that $u^*_{y_N}=0$) and, what only matters,  is uniformly bounded by some constant depending on $K, \|DV\|_\infty$ and $\|D^2T\|_\infty$. Finally, the boundary operator $N^\tau$ is Lipschitz  with respect to $y$ and $C^2$ with respect to the gradient variable. According to  \cite[Lemma 2.3]{lieb84NA}, $\mD u^*$ satisfies a bound in the H\"older norm, which is uniform with respect to $\tau$. With a standard localization argument, based on a partition of unity, 
we conclude the bound \rife{gradhol}.
\qed
\vskip1em

\begin{remark}\label{plaplacian}  We stress that a $C^{1,\alpha}$ estimate (up to the boundary) can also be obtained  if we look at the elliptic equation in the divergence form, which  is clearly inherited from the continuity equation  in \rife{mfpp}. In the notations used in the previous proof, we observe that the function $u^*$, obtained after local change of coordinates and reflection, satisfies the equation
\begin{align*}
 & \partial_t  m^\tau   - \sum_{k=1}^N \frac{\partial}{\partial y_k}\left(m^\tau 2h'(|MDu^*|^2) DT_k\cdot (M Du^*)\right) + \rho \, u^* \frac{m^\tau}{1+ \tau f'(m^\tau)m^\tau} =0 
 \\
 &\qquad  \hbox{where $m^\tau = \vfi^\tau\left( -u_t + h( |MDu^*|^2 )- \tau V \right)$,}
\end{align*}
where we recall that $\vfi^\tau$ is defined in \rife{mtau} and $M= (m_{ij})= (\frac{\partial T_j}{\partial x_i})$, $|MDu^*|^2= |u_{y_N}^*|^2+ \sum_{i=1}^N\left(\sum_{\ell=1}^{N-1}\frac{\partial T_\ell}{\partial x_i}u^*_{y_\ell}\right)^2$.
Hence the equation can be rewritten as a nonlinear divergence form equation on $u^*$:
$$
\partial_t \left( A_0(y,\mD u^*) \right) + \sum_{k=1}^N \frac{\partial}{\partial y_k}\left(A_k(y, \mD u^* )\right) = \beta(y, u^*, \mD u^*)\qquad (t,y)\in (0,1)\times B
$$
complemented with a co-normal boundary condition at $t=0,t=1$:
\begin{align*}
{\bf A}\cdot {\bf n} & = \left( A_0(y,\mD u^*) \mathop{|}_{t=1}  , -A_0(y,\mD u^*) \mathop{|}_{t=0}\right)
\\ &  = \left( \vfi^\tau( \psi_1^* + \tau (f(m_1)^*+  \tilde V^*-V^*)+ \de u^*) ,  - \vfi^\tau( \psi_0^* + \tau (f(m_0)^*+  \tilde V^*-V^*)- \de u^*)\right) 
\end{align*}
where $(\cdot)^*$ denotes the even extension of the various functions. 
Notice that when $\tilde V=V$, $\de=0$  and $\psi_i=\log(m_i)$, then we have ${\bf A}\cdot {\bf n}=(m_1^*, -m_0^*)$. 

Since $|\beta(y, u^*, \mD u^*)|\leq C$ (because $\mD u^*$ is bounded), and due to the regularity of the boundary terms, if $H_{pp}$ is nondegenerate the $C^{1,\alpha}$ estimate for $u^*$ can then be deduced by \cite[Thm 2]{lieb88}.

We mention this alternative approach because it could be exploited if one aims at generalizing  the $C^{1,\alpha}$ estimate to possibly degenerate Hamiltonians, e.g. satisfying \rife{Hppq} with $\varpi=0$. 
\end{remark}

\vskip1em
Last ingredient is the well-posedness of the linearized problem around one solution $u$ of \rife{LU} and, correspondingly, a compactness property  for sequences of solutions of \rife{LU}.  In order to use the $C^{2,\alpha}$ regularity, here we need to use a compatibility condition at the boundary  $\partial\Omega$. To this purpose, we denote  $C^{1,\alpha}_N(\oo)$ the space of functions $\phi \in C^{1,\alpha}(\oo)$ such that $D\phi\cdot \vec \nu=0$ on $\partial \Omega$. We also denote $a_{ij}^\tau= a_{ij}^\tau(\eta,q)$.

\begin{proposition}\label{linearized}  In addition to the  conditions  of Proposition \ref{boundK}, assume that $\partial\Omega$ is of class $C^4$, and that the functions $f(m_0)$, $f(m_1)$, $\tilde V(x)$, $\psi_0,\psi_1$  belong to  $C^{1,\alpha}_N(\oo)$. Then the  solution  $u$ of \rife{ellip} belongs to $C^{2,\alpha}(\oQ)$, and  for every $\vartheta\in C^{0,\alpha}(\oQ), \zeta= (\zeta_1, \zeta_0) \in C^{1,\alpha}_N(\oo)$ there exists a unique $\phi\in C^{2,\alpha}(\oQ)$ which is a solution of  the linear problem
\be\label{linea}
\begin{cases}
-a_{ij}^\tau((t,x),\mD u) \phi_{ij} - \frac{\partial a^\tau_{ij}}{\partial q_k}((t,x), \mD u) u_{ij} \phi_k + \tau H_{pp}(Du) DV \cdot D\phi  + \rho \phi = \vartheta  
& \hbox{in $Q$ }
\\
- \phi_t + H_p(Du) \cdot D\phi = \de \phi +   \zeta_1 & \hbox{on $\Sigma_T$,}
\\
   -  \phi_t + H_p(Du) \cdot D\phi + \de \phi=  \zeta_0  & \hbox{on $\Sigma_0$,}
\\
  D\phi\cdot \vec\nu=0 & \text{on $(0,T)\times \partial \Omega$.}
\end{cases}
\ee
Moreover, if $\psi_{1m}, \psi_{0m}$ are  sequences which  converge in $C^{1,\alpha}_N(\oo)$, then the  corresponding solutions $u_m$ of  \rife{ellip}  are  relatively compact in $C^{2,\alpha}(\oQ)$.
\end{proposition}

\proof  The first assertion follows from Proposition \ref{boundaryHol} if we  come back to the equation \rife{reflected} complemented  with the boundary condition satisfied by $u^*$, which reads as 
$$
\begin{cases}
  -u_t^* + h( |MDu^*|^2 )- \de u^* = \tau (f(m_1)^*+\tilde V^*) + \psi_1^* & \quad \hbox{at $t=1$,}\\
 -u_t^* + h( |MDu^*|^2 )+ \de u^* = \tau (f(m_0)^*+\tilde V^*) + \psi_0^* & \quad \hbox{at $t=0$.}
\end{cases}
$$
Since the boundary data are assumed to  belong to  $C^{1,\alpha}_N(\oo)$, they are  reflected into $C^{1,\alpha}$ functions. This implies  that $u^*_{t}-h( |MDu^*|^2 ) \in C^{1,\alpha}(U)$. One can therefore apply the results of \cite{lieb82} to invoke that $u^*\in C^{2,\alpha}$, which leads to the $C^{2,\alpha}$ regularity of $u$ up to the boundary.

Let us now consider the linearized problem \rife{linea}. Since $u\in C^{2,\alpha}(\oQ)$, \rife{linea} is a linear problem with H\"older coefficients, and Schauder's theory applies to get, at least, the interior estimates of $\phi$. Let us only check the boundary regularity.  With the notations of Proposition \ref{boundaryHol},   we localize near the boundary with the change of variable $y=T(x)$, and then we extend the problem by reflection, setting  $\phi^*(t,y)$ the (even) reflection of $\tilde \phi(t,y):= \phi(t,T^{-1}(y))$ through the hyperplane $\{y_N=0\}$. 
From the computations of Proposition \ref{boundaryHol}, we know that the matrix ${\mathcal A}^\tau$ is transformed into a new matrix ${\mathcal A^*}$ which depends on $(y, \cD u^*)$ as follows:
$$
{\mathcal A^*}= \left(a_{ij}^*\right) = 
\begin{pmatrix} 1 & -2h' \Lambda_u \\ - 2h' \, \Lambda_u & [4(h')^2] \Lambda_u\otimes \Lambda_u\end{pmatrix}
+{\mathcal \kappa}^\tau (V^*(y), \cD u^*) \begin{pmatrix} 0 & 0 \\ 0 &  
 (2h') \left(DT_i\cdot DT_j\right)+ (4h'') \Lambda_u\otimes \Lambda_u  \end{pmatrix}
$$
with $\Lambda_u:= \left( (DT_k\cdot M'D'u^*)_{k=1,\dots N-1}, u^*_{y_N}\right) $ and ${\mathcal \kappa}^\tau (V^*(y), \cD u^*)$ defined as in \rife{kastar} .

We look now at the first order terms. By definition of $a_{ij}^\tau$, we have 
%Let us recall that 
%$$
%{\mathcal A}= \left(a_{ij}^\tau((t,x),\mD u)\right) = 
%\begin{pmatrix} 1 & - H_p(Du) \\ - H_p(Du) &   H_p(Du)\otimes H_p(Du)\end{pmatrix}
%+( 1+ \tau\, m^\tau\, f'(m^\tau ) ) \begin{pmatrix} 0 & 0 \\ 0 &  
%H_{pp}\end{pmatrix}
%$$
%where $m^\tau$ is defined in \rife{mtau}.  Hence if $\eta= (s, \xi)$ and $g(s)= sf'(s)$, we have
%\begin{align*}
%& \frac{\partial a^\tau_{ij}}{\partial s}((t,x), \mD u) u_{ij}=  -   \tau \, g' (m^\tau)(\vfi^\tau)' H_{p_ip_j}u_{ij}\\
%& 
%\frac{\partial a^\tau_{ij}}{\partial \xi_k}((t,x), \mD u) u_{ij} =  -2H_{p_ip_k}u_{tx_i} + 2 H_{p_ip_k}H_{p_j}u_{ij} + \tau  \, g' (m^\tau) (\vfi^\tau)' H_{p_k}  H_{p_ip_j}u_{ij}+ (1+\tau \, g(m^\tau)) (H_{p_ip_jp_k} u_{ij})
%\end{align*}
%which yields
\be\label{newterm}\begin{split}
\frac{\partial a^\tau_{ij}}{\partial q_k}((t,x), \mD u) u_{ij} \phi_k & = -2 H_{pp}Du_t \cdot D\phi+ 2 D^2u H_p \cdot H_{pp}D\phi + 
\tau \, g' (m^\tau)(\vfi^\tau)' \,  [-\phi_t +   H_p D\phi]\,  {\rm tr}(H_{pp}D^2u)
\\ & \qquad \quad +   (1+\tau \, g(m^\tau)){\rm tr} (H_{ppp}D\phi D^2 u)
\,
\end{split}
\ee 
where $g(s)= sf'(s)$. 
Using the radial structure of $H$ (\rife{rad}) and the chosen reference frame where \rife{ortho} holds true, one  can check (with long but routine computations, similarly as in Proposition \ref{boundaryHol}) that all terms in \rife{newterm} can be transformed into functions of $\phi^*$ and $u^*$ except for  a few terms which are only first order in $u$ and terms which involve second derivatives of the map $T$. 
Then, after the change of coordinates and the reflection,  $\phi^*$ satisfies a linear  elliptic  problem of the following form in the open set $(0,T)\times B$:
\be\label{eqpsi*}
a_{ij}^* \phi^*_{ij}  + B^*_i \phi^*_i + c^*+ \rho \phi^* = \vartheta^* \qquad \hbox{$(t,y)\in (0,T)\times B$}\,.
\ee
Now we wish to apply to $\phi^*$ the H\"older estimates up to the boundary in order to conclude that $\phi^* \in C^{2,\alpha}([0,T]\times U \times (-d_0, d_0))$ for some open set $U\subset \R^{N-1}$.  The regularity of $\phi^*$ will yield the $C^{2,\alpha}$ regularity of $\phi$ in $[0,T]\times ({\mathcal B}\cap \overline\Omega)$ for any ball $\mathcal B$ such that $\partial \Omega \cap {\mathcal B}$ is the portion of a smooth graph.

In \rife{eqpsi*}, the ($d+1$-dimensional) vector field $B^*$ depends on $\cD^2u^*$ (which is H\"older continuous) and on the (even reflection) of $DT, DV, D\tilde u$, which are  Lipschitz continuous. Therefore,  overall $B^*$ belongs to $C^{0,\alpha}$ for some $\alpha\in (0,1)$. Conversely, the term $c^*$ in \rife{eqpsi*} is the even reflection of a function depending on $D\phi$, $D\tilde u$ and $D^2 T$; in this case we cannot say  that $c^*$ is reflected into a H\"older continuous function until we establish the H\"older regularity for $D\phi$.
%\begin{align*}
%B^*_0& = - \tau \, h' (m^*)  \left(   \sum_{p,q} u^*_{y_py_q} (DT_p\cdot DT_q ) + \sum_{p}u_{y_p} \Delta  T_p\right)
%\\
%\hbox{for $\ell=1,\dots N$,}\,\, \quad B^*_\ell & = -2 (DT)^* Du_t^* DT_\ell   + 2 \sum_{p,q,k} u_{y_py_q}^* (DT_p\cdot DT_k)(DT_q\cdot DT_\ell)u_{y_k}^*
%\\ & \qquad 
%+   \tau \, h' (m^*)  \left( \tau \sum_{p,q,k} u^*_{y_py_q} (DT_p\cdot DT_q )   u^*_{y_k}  (DT_k\cdot DT_\ell) +   \tau \sum_{p,k}u_{y_p} \Delta  T_pu^*_{y_k}  (DT_k\cdot DT_\ell)\right)
%\end{align*}
%and
%$$
%c^*(t,y)= \begin{cases}
%\sum_{\ell=1}^{N-1}  b_\ell^\tau(t,y',y_N) \,  \psi^*_{y_\ell} + b_N^\tau (t,y',y_N)\, \psi^*_{y_N} & \hbox{if $y_N\geq 0$}, 
%\\
%\sum_{\ell=1}^{N-1}  b_\ell(t,y',-y_N) \,  \psi^*_{y_\ell} - b_N^\tau(t,y',-y_N)\, \psi^*_{y_N} & \hbox{if $y_N<0$}
%\end{cases}
%$$
%where 
%$$
%b_\ell^\tau=  2 \sum_{p,k} (D^2T_p DT_k\cdot DT_\ell ) u_{y_p} u_{y_k}+ \sum_{i,j}a_{ij}^\tau(\mD  u )\frac{\partial^2 T_\ell}{\partial x_i\partial x_j} +   \tau DV\cdot DT_\ell 
%$$
%all functions here computed at $(t,x)$, i.e. $(t,T^{-1}(y))$.  

Therefore, since  the coefficient $c^*$  does not belong to $C^{0,\alpha}$ a priori, 
 we cannot apply the Schauder's estimates in one shot to infer the $C^{2,\alpha}$ regularity of $\phi^*$. However, the coefficients $a_{ij}^*$ belong to $C^{0,1}$ (because $u^*$ is $C^{2,\alpha}$, while $V, DT$ are $C^1$ and so they are reflected into a Lipschitz function); hence, by rewriting the equation in divergence form, $\phi^*$ is, in particular,  a bounded weak solution to the equation 
\be\label{cono}
(a_{ij}^* \phi^*_{j})_i  + H(y, \phi^*, \phi^*_i)=0\qquad (t,y)\in (0,T)\times U\times (-d_0, d_0)
\ee 
for some function $H(y, \phi^*, \phi^*_i)$ which satisfies $|H(y, \phi^*, \phi^*_i)|\leq C_0 (1+ |D\phi^*|)$. In addition, $\phi^*$ satisfies  a co-normal boundary condition at $\Sigma_0, \Sigma_T$. By  regularity for divergence form equations (which is true for general nonlinear problems, see e.g. \cite[Thm 2]{lieb88}), we  deduce that  $\phi^*\in C^{1,\alpha}$ and 
$$
\|\phi^*\|_{C^{1,\alpha}}\leq K\qquad  \hbox{for some $K= K(\|\phi\|_\infty, \|u\|_{C^{2,\alpha}}, \|T\|_{2}, \|DV\|_{\infty}, \|\vartheta\|_\infty, \|\zeta\|_{\alpha}, Q)$.}
$$ 
We can now use this information, which implies that the term $c^*$ in \rife{eqpsi*} belongs to $C^{0,\alpha}$. Since $a^*_{ij}$, $B^*$, $\vartheta^*$ also belong to $C^{0, \alpha}$, and since  the boundary  data $\zeta^* \in C^{1,\alpha}$ (because $\zeta$ satisfies a Neumann condition),   the classical Schauder's estimates (see e.g. \cite[Lemma 1]{lieb82}) imply that $\phi^*\in C^{2,\alpha}$ and 
$$
\|\phi^*\|_{C^{2,\alpha}}\leq K\qquad  \hbox{for some $K= K(\|\phi\|_\infty, \|u\|_{C^{2,\alpha}}, \|T\|_{C^{2,\alpha}}, \|V\|_{C^{1,\alpha}}, \|\vartheta\|_{\alpha}, \|\zeta\|_{C^{1,\alpha}}, Q)$.}
$$ 
We notice that $\|\phi\|_\infty$ is estimated uniformly by maximum principle.  We also recall that, requiring $\partial\Omega$ of class $C^4$, we have that $T$ is of class $ C^3$. Finally, we proved that all solutions to the linear problem \rife{linea} satisfy the estimate
$$
\|\phi\|_{C^{2,\alpha}}\leq C \left( \|\vartheta\|_{\alpha}+  \|\zeta\|_{C^{1,\alpha}}\right) 
%\qquad \forall \tau\in [0,1]\,.
$$
for some $C= C( \|u\|_{C^{2,\alpha}}, \partial \Omega, \|V\|_{C^{1,\alpha}})$. By linear theory, this implies that problem \rife{linea}
 is uniquely solvable for every $\tau\in [0,1]$.
 
Now we prove the last assertion of the Proposition.  To this goal, we first observe that, since $\psi_m=(\psi_{0m}, \psi_{1m})$ is convergent in $C^{1,\alpha}_N(\oo)$, by Proposition \ref{boundaryHol} the sequence $u_m$ is bounded in $C^{1,\alpha}(\oQ)$, and therefore it is relatively compact in $C^1(\oQ)$ by Ascoli-Arzel\`a\ theorem. We observe that $w^{m,n}:= u_m-u_n $ solves 
$$
\begin{system}
& -a_{ij}^\tau(x,\mD u_m) w^{m,n}_{ij}  + B_k^{n,m} w^{m,n}_k+ \tau \, \Gamma^{n,m}DV(x) \cdot Dw^{m,n} + \rho w^{m,n}= 0
 &&\text{in }  (0,1)\times \Omega, \\
& -\partial_t w^{m,n} +   \beta^{n,m}\cdot Dw^{m,n} - \de w^{m,n} = \psi_{1n}-\psi_{1m}   
%f^\vep(m_1)  
&&\hbox{at $t=1$, $x\in \Omega$,}
\\
&   -\partial_t w^{m,n} +   \beta^{n,m}\cdot D w^{m,n}   + \de w^{m,n}=  \psi_{0n}-\psi_{0m} 
 %f^\vep (m_0)  
 &&\hbox{at $t=0$, $x\in \Omega$,}
\\
& Dw^{m,n}\cdot \vec\nu=0 &&\text{on $(0,1)\times \partial \Omega$,}
\end{system}
$$
 where
\begin{align*}
&  B_k^{n,m}= \sum_{ij}(u_n)_{ij} \int_0^1 \frac{a_{ij}^\tau}{\partial q_k} (x, \lambda Du_m+ (1-\lambda)Du_n) d\lambda
\\
& \beta^{n,m}= \int_0^1 H_p(\lambda Du_m+ (1-\lambda)Du_n) d\lambda\,,\qquad \Gamma^{n,m}= \int_0^1 H_{pp}(\lambda Du_m+ (1-\lambda)Du_n) d\lambda \,.
\end{align*}
This is  a linear problem for $w^{m,n}$ which is of the same kind as \rife{linea}. Therefore, we proceed as we did before, by localizing near the boundary and reflecting the solution through (the straightened part of)  $\partial \Omega$. Then we apply to the reflected function $(w^{n,m})^*$ the $C^{2,\alpha}$ estimates as before (\cite[Lemma 1]{lieb82}). We only stress that, differently than it was for the function $\phi^*$ above, now we already know that $w^{n,m}$ is bounded in $C^{1,\alpha}(\oQ)$; this allows us to apply directly the Schauder estimates because (in the notation used before) $a_{ij}^* (w^{n,m})^*_{ij} $ belongs to (and can be estimated in) $C^{0,\alpha}(\overline Q)$. Applying the precise form of  \cite[Lemma 1]{lieb82}, we conclude that
\be\label{wnm}
 \|w^{n,m}\|_{C^{2,\alpha}} \leq C_0 \|w^{n,m}\|_{C^1} (1+ \|u_n\|_{C^{2,\alpha}} + \|u_m\|_{C^{2,\alpha}} )+ C_1 \|w^{n,m}\|_{C^{1,\alpha}} (\|D^2u_n\|_{\infty} + \|D^2u_m\|_{\infty} )
\,.
\ee
From this estimate we conclude as in \cite[Lemma 2]{lieb82}; one shows first that $\|u_n\|_{C^{2,\alpha}}$ is bounded (using the smallness of $ \|w^{n,m}\|_{C^1}$ in \rife{wnm}) and this implies the compactness of $u_n$ in $C^{C^{1,\alpha}}$. Hence \rife{wnm}  yields the convergence of $u_n$  in $C^{2,\alpha}(\overline Q)$.
\qed

\vskip0.4em
We are ready to prove the existence of solutions to the elliptic problem \rife{eps-de}.

\vskip0.6em

{\bf Proof of Theorem \ref{lieb}.} \quad  We start by assuming that the data are more regular, namely that $\partial \Omega$ is of class $C^4$, $V\in C^{2,\alpha}(\overline \Omega)$,  $f\in C^{2,1}(0,\infty)$ satisfies \rife{fprime-mod} and $H\in C^{3,1}(\R^d)$ satisfies \rife{rad}. We also assume that the boundary data $f(m_0)$, $f(m_1)$, $\tilde V(x)$, $\psi_0,\psi_1$  belong to  $C^{1,\alpha}_N(\oo)$.
%, in order to ensure the compatibility condition at $\partial\Omega$ for the boundary operator $N^\tau$. 
Under those conditions, we can use Propositions \ref{boundK}--\ref{linearized} and we establish the existence of a solution with a continuity method,  following the same proof as in \cite[Chapter X, Thm 1.1 \& Thm 1.2]{LU}. Let us sketch some detail: we set 
$$
X:= \{v\in C^{2,\alpha}(\oQ): Dv\cdot \vec\nu=0 \,\,\hbox{on}\, (0,T)\times \partial\Omega\}
$$
and $X':= C^{0,\alpha}(\oQ)\times C^{1,\alpha}(\oo)^2$, and we define a family of mappings from $X $ into $X' $ by setting $\Phi(u,\tau):= (L^\tau(u), N^\tau_0(u), N^\tau_1(u))$, where $N^\tau_0, N^\tau_1$ are the restrictions to $t=0$, $t=1$ respectively, of the boundary operator $N^\tau$ in \rife{ntau}. We also set 
$X'_0$ the closed subset of $X'$ made of elements $(0, \psi_0, \psi_1)$, with $(\psi_0, \psi_1)\in C^{1,\alpha}_N(\oo)^2$. Thanks to Proposition \ref{linearized}, the mapping $\Phi$ is locally invertible in a neighborhood of points $u\in X$ such that $\Phi(u,\tau)\in X_0'$, and in addition a sequence $u_k$ is compact in $X$ if  $\Phi(u_k,\tau)$ is converging in $X'_0$. Finally, for $\tau=0$ and $\psi_0=\psi_1=0$, the problem \rife{ellip} has the unique solution $u=0$. Therefore,   as in  \cite[Chapter X, Thm 1.1]{LU} one concludes that the problem with $\tau=1$ is uniquely solvable. 

In a second step, we approximate the data with smoother sequences.  In particular, we approximate the domain $\Omega$ with a sequence of $C^4$ domains $\Omega_n$, and the functions $m_0, m_1$ with sequences $m_{0n}, m_{1n}\in C^{1,\alpha}(\oo)$ such that $\|m_{0n}\|_{W^{1,\infty}(\Omega)}, \|m_{1n}\|_{W^{1,\infty}(\Omega)}$ are bounded, as well as $\|\log(m_{0n})\|_{W^{1,\infty}(\Omega)}$ and $\|\log(m_{1n})\|_{W^{1,\infty}(\Omega)}$. It is also possible to build $m_{0n}, m_{1n}$ in a way that those functions be constant along the normal in a small neighborhood of $\partial\Omega$, so that they belong to $C^{1,\alpha}_N(\oo)$. Similarly we take a sequence $\tilde V_n\in C^{1,\alpha}_N(\oo)$ which converges uniformly to $V$ and such that  $\|\tilde V_n\|_{W^{1,\infty}(\Omega)}$ is bounded. Finally, we take another sequence $V_n\in C^{2,\alpha}(\oo)$ converging to $V$ and such that $\|V_n\|_{W^{2,\infty}(\Omega)}$ is bounded. As for the function $f$, we approximate it with 
$f_n(m)= f(\frac m{1+ \frac1n m})$ which satisfies the stronger condition \rife{fprime-mod} (for some constant $\beta$ possibly dependent on $n$) but also satisfies the 
weaker condition \rife{fprime} for some constant which is uniform with respect to $n$.
As for the function $H$, by assumption it is $C^3$ and satisfies \rife{hrad} and \rife{Hppq}-\rife{hpppq} for some $\varpi>0$; hence, up to replacing $h$ with $h(\sqrt{n^{-1}+ t^2})$, and up to a  further regularization if needed, we can assume that $H$ is approximated by some $H^n$ which is  $C^{3,1}(\R^d)$ and satisfies \rife{rad}, and in a way  that \rife{Hppq}-\rife{hpppq} hold uniformly for $H^n$.  We build so far an approximating problem
$$
\begin{system}
& -{\rm tr}\left({\mathcal A}_n(x,{\mathcal D}u_n) \, {\mathcal D}^2u_n \right)+ \rho u_n + DV_n(x)\cdot H_p^n(Du_n) =0 &&\text{in } Q, \\%(0,1)\times \R^d, \\
& -\partial_t u_n + H^n(Du_n)= \de u_n+ f^\vep (m_{1n}) + \tilde V_n(x) &&\hbox{at $t=1$, $x\in \Omega$,}
\\
&   -\partial_t u_n + H^n(Du_n)+ \de u_n= f^\vep(m_{0n})+ \tilde V_n(x)   &&\hbox{at $t=0$, $x\in \Omega$,}
\\
& \,\, H_p^n(Du_n)\cdot \vec\nu=0 &&\text{on $(0,1)\times \partial \Omega$.}
%\\%(0,1)\times \R^d,\\
%&m(0,\cdot) = m_0, \; m(1,\cdot)  = m_1 &&\text{in } \R^d.
\end{system}
$$
which admits a solution $u_n\in C^{2,\alpha}(\oQ)$ by what we proved in the first part. This solution $u_n$ also belongs to $C^3(Q)$ and we can apply Lemma \ref{sup} and Theorem \ref{bound} to infer a uniform estimate for $\|u_n\|_{W^{1,\infty}(Q)}$. By Proposition \ref{boundaryHol}, we also deduce that $u_n$ is uniformly bounded in $C^{1,\alpha}(\oQ)$, because the $C^{1,\alpha}$ estimate only depends on the bound for  $\|u_n\|_{W^{1,\infty}(Q)}$ and on the Lipschitz bounds of the boundary terms $f^\vep(m_{in})$, $ \tilde V_n$. Hence $u_n$ is relatively compact in $C^1(\oQ)$. By interior regularity, we also have that $u_n$ is bounded in $C^{2,\alpha}(Q)$, so it is relatively compact in $C^2$ if restricted to  compact subsets in the interior. This is enough to pass to the limit and conclude that, up to subsequences,  $u_n$ converges towards some $u\in C^{2,\alpha}(Q)\cap C^{1,\alpha}(\oQ)$ which solves problem \rife{eps-de} and satisfies estimate \rife{lastpray}.
\qed

 \vskip2em
 
{\bf Acknowledgement.} I warmly thank Giuseppe Savar\'e for stimulating my interest in this problem and for sharing with me several opinions and hints on the topics
of the paper.

I also thank Filippo Santambrogio for pointing me out reference \cite{SaWa}.

\end{document}